\documentclass[11pt]{article}
\usepackage{amssymb,latexsym,amsmath,amsbsy,amsthm,amsxtra,amsgen}
\numberwithin{equation}{section}

\newfont{\msbm}{msbm10 at 11pt}

\newcommand {\Z} {\mbox{\msbm Z}}

\newcommand {\N} {\mbox{\msbm N}}

\def\ni{\noindent}

\def\be{\begin{equation}}
\def\ba{\begin{align}}
\def\ea{\end{align}}
\def\mn{\bigskip\noindent}
\def\Var{\textup{Var}}
\def\Cov{\textup{Cov}}
\def\ln2b{O \bigg( \frac{1}{(\log N)^2} \bigg)}
\def\upat{\vskip-0.2in\noindent}

\def\Box{\vcenter{\vbox{\hrule height .4pt
  \hbox{\vrule width .4pt height 7pt \kern 7pt
        \vrule width .4pt} \hrule height .4pt}}}

\newtheorem{Theo}{Theorem}[section]
\newtheorem{Lemma}[Theo]{Lemma}

\newtheorem{Prop}[Theo]{Proposition}

\newtheorem{Dfn}[Theo]{Definition}

\begin{document}
\title{Random partitions approximating the coalescence \\ of lineages during a selective sweep}
\author{by Jason
Schweinsberg\thanks{Supported by an NSF Postdoctoral Fellowship while the author was at Cornell 2001--2004.}
\kern.5em and Rick Durrett\thanks{Partially supported by NSF grants from
the probability program (0202935) and from a joint DMS/NIGMS initiative to support
research in mathematical biology (0201037).} \\ \\
University of California at San Diego and Cornell University}
\maketitle

\begin{abstract}
When a beneficial mutation occurs in a population, the new, favored allele
may spread to the entire population.  This process is known as a selective
sweep.  Suppose we sample $n$ individuals at the end of a selective sweep.
If we focus on a site on the chromosome that is close to the location of
the beneficial mutation, then many of the lineages will likely
be descended from the individual that had the beneficial mutation, while others
will be descended from a different individual because of recombination
between the two sites.  We introduce two approximations for the effect of
a selective sweep. The first one is simple but not very accurate:
flip $n$ independent coins with probability $p$ of heads and say that the
lineages whose coins come up heads are those that are descended from the
individual with the beneficial mutation. A second approximation, which is
related to Kingman's paintbox construction, replaces the coin flips by
integer-valued random variables and leads to very accurate results. 
\end{abstract}

\bigskip
\noindent Running head: Approximating a selective sweep.


\clearpage

\section{Introduction}
A classical continuous-time model for a population with overlapping
generations is the Moran model, which was introduced by Moran (1958).
Thinking of $N$ diploid individuals, 
we assume the population size is fixed at $2N$.  
However under the assumption that each individual is a random union of gametes,
the dynamics are the same as for a population of $2N$ haploid individuals,
so we will do our computation for that case.  In the 
simplest version of the Moran model, each individual
independently lives for a time that is exponentially distributed with
mean $1$ and then is replaced by a new individual.  The parent of the new
individual is chosen at random from the $2N$ individuals, including the
individual being replaced.

Here we will consider a variation of the Moran model that involves two loci,
one subject to natural selection, the other neutral, and with a probability
$r$ in each generation of recombination between the two loci. To begin
to explain the last sentence, we assume that at the selected locus there are two alleles,
$B$ and $b$, and that the relative fitnesses of the two alleles are $1$ and
$1-s$. The population then evolves with the same
rules as before, except that a replacement of an individual with a $B$
allele by an individual with a $b$ allele is rejected with probability
$s$.  Consequently, if at some time there are $k$ individuals with the
$B$ allele and $2N - k$ with the $b$ allele, then the rate of transitions
that increase the number of $B$ individuals from $k$ to $k+1$ is
$k(2N - k)/(2N)$, but the rate of transitions that reduce the number of
$B$ individuals to $k-1$ is $k(2N - k)(1 - s)/(2N)$.  See 
chapter 3 of Durrett (2002) for a summary of some work with this model.

We assume that the process starts at time zero with $2N-1$
individuals having the $b$ allele and one individual having the advantageous
$B$ allele.  We think of the individual with the $B$ allele as having had
a beneficial mutation at time zero.  There is
a positive probability that eventually all $2N$ individuals will have
the favorable allele.  When this happens, we say that a selective sweep
occurs, because the favorable allele has swept through the entire population.

If we assume that the entire chromosome containing the selected locus 
is passed down from one generation to the next, as is the case for the
$Y$ chromosome or mitochondrial DNA, then all $2N$
chromosomes at the end of the selective sweep
will have come from the one individual that
had the beneficial mutation at the beginning of the sweep.
However, non-sex chromosomes in diploid individuals are typically not
an identical copy of one of their parents' chromosomes. Instead, because of
a process called recombination, each chromosome that an individual inherits
consists of pieces of each of a parent's two chromosomes. In this case, if 
we are interested in the origin of a second neutrally evolving locus
on the chromosome
and a selective sweep occurs because of an advantageous mutation at a
nearby site, then some of the lineages will be traced back to the chromosome
that had the favorable allele at the beginning of the sweep but other
lineages will be traced back to different individuals because of
recombination between the neutral and selected loci.
When a lineage can be traced back to an individual other than the one with the
beneficial mutation, we say that the lineage
escapes from the selective sweep.

The combined effects of recombination and selective sweeps have been studied
extensively.  Maynard Smith and Haigh (1974) observed that selective
sweeps can alter the frequencies of alleles at sites nearby the site at
which the selective sweep occurred.  They referred to this as the
``hitchhiking effect''.  They considered a situation with a neutral 
locus with alleles $A$ and $a$ and a second locus where allele $B$ 
has a fitness of $1+s$ relative to $b$. Suppose $p_0$ is the initial 
frequency of the $B$ allele, and $Q_n$ and $R_n$ are the frequencies in 
generation $n$ of the $A$ allele on chromosomes containing $B$ and $b$
respectively.
If $Q_0=0$ (i.e., the advantageous mutation arises on a chromosome with
the $a$ allele) 
and the recombination probability in each generation is $r$, Maynard Smith 
and Haigh (1974) showed (see (8) on page 25) that the frequency of the 
$A$ allele after the selective sweep is reduced from $R_0$ to
\begin{displaymath}
\lim_{n\to\infty} Q_n = R_0 \sum_{n=0}^\infty r (1-r)^n \cdot
\frac{1-p_0}{1 - p_0 + p_0(1+s)^{n+1}}.
\end{displaymath}

In the calculation of Maynard Smith and Haigh, the number of individuals
with the $B$ locus grows deterministically. Kaplan, Hudson, and Langley (1989)
used a model involving an initial phase in which the number
of $B$'s is a supercritical branching process, a middle deterministic piece
where the fraction $p$ of $B$'s follows the logistic differential equation
\be
\frac{dp}{dt} = s p(1-p),
\label{logde}
\end{equation}
and a final random piece where the number of $b$'s follows a subcritical
branching process. This process is too difficult to study analytically
so they resorted to simulation.

Stephan, Wiehe, and Lenz (1992) further simplified this approach
by ignoring the random first and third phases and modeling the change
in the frequency of $B$'s by the logistic differential equation (\ref{logde}),
which has solution
\begin{displaymath}
p(t) = \frac{p(0)}{p(0) + (1-p(0)) e^{-st}}.
\end{displaymath}
This approach has been popular with biologists in simulation studies
(see, for example, Simonsen, Churchill, and Aquadro (1995) and
Przeworski (2002)).  However, as results in Barton (1998) and
Durrett and Schweinsberg (2004a) show, this can introduce substantial errors,
so rather than using this approximation for
our analysis, we will consider a modification of the Moran model that
allows for recombination as well as beneficial mutations.

We consider two sites on each chromosome.  At one site,
each of the $2N$ chromosomes has either the advantageous $B$ allele
or a $b$ allele.  Our interest, however, is in the
genealogy at another neutral site, at which all alleles
have the same fitness. As before, we assume that each individual lives for
an exponential time with mean $1$ and is replaced by a new individual
whose parent is chosen at random from the population, except that
we disregard disadvantageous replacements of a $B$ chromosome by a
$b$ chromosome with probability $s$.  We will also now assume that
when a new individual is born, it inherits alleles at both sites
from the same individual with probability $1 - r$.  With probability $r$,
there is recombination between the two sites, and the
individual inherits the allele at the neutral site from its
parent's other chromosome.  Since a parent's two chromosomes are 
considered to be two distinct individuals in the population, we model this
by saying that the new individual inherits the two alleles from
two ancestors chosen independently at random from the $2N$ individuals 
in the population.

Suppose we sample $n$ chromosomes at the end of a
selective sweep and follow their
ancestral lines back until the beginning of the sweep.  We will describe
the genealogy of the sample by a marked partition of $\{1, \dots, n\}$,
which we define to be a partition of $\{1, \dots, n\}$ in which one
block of the partition may be designed as a ``marked'' block.  We define the
marked partition $\Theta$ of $\{1, \dots, n\}$ as follows.  We say that
two integers
$i$ and $j$ are in the same block of $\Theta$, denoted $i \sim_{\Theta} j$,
if and only if the alleles at the neutral site on the $i$th and $j$th
chromosomes in the sample have the same ancestor at the beginning of the
sweep.  Thus, if we are following the lineages associated with the allele at
the neutral site, we have $i \sim_{\Theta} j$ if and only if the $i$th
and $j$th lineages coalesce during the selective sweep.  We also
``mark'' the block of $\Theta$ containing the integers $i$ for which the
$i$th individual is descended from the individual that had the beneficial
mutation at the beginning of the sweep.  Thus, to understand
how a selective sweep affects the genealogy of a sample of size $n$, we
need to understand the distribution of the random marked partition $\Theta$.

In this paper, we study two approximations to the distribution
of $\Theta$.  The approximations were introduced, and studied by
simulation, in Durrett and Schweinsberg (2004a).  Here we
provide precise bounds on the error in the approximations.
The idea behind the first approximation is that a large number of
lineages will inherit their allele at the neutral site from the
individual that had the beneficial mutation at the beginning of
the sweep, and the corresponding integers $i$ will be in the
marked block of $\Theta$.  With high probability, the lineages
that escape the selective sweep do not coalesce with one another,
so the corresponding integers are in singleton blocks of $\Theta$.

Before stating the first approximation precisely, we need a definition.
Let $p \in [0,1]$.  Let $\xi_1, \dots, \xi_n$ be independent random
variables such that, for $i = 1, \dots, n$, we have $P(\xi_i = 1) = p$ and
$P(\xi_1 = 0) = 1-p$.  We call the random
marked partition of $\{1, \dots, n\}$ such that one marked
block consists of $\{i \in \{1, \dots, n\}: \xi_i = 1\}$ and
the remaining blocks are singletons a $p$-partition of $\{1, \dots, n\}$.
Let $Q_p$ denote the distribution of a $p$-partition of $\{1, \dots, n\}$.

Theorem \ref{logthm} below shows that the distribution of $\Theta$ can be
approximated by the distribution of a $p$-partition.  For this result,
and throughout the rest of the paper, we assume that the selective
advantage $s$ is a fixed constant that does not depend on the population
size $N$.  However, the recombination probability $r$ is allowed to depend
on $N$, even though we have not recorded this dependence in the notation.
We will assume throughout the paper that $r \leq C_0/(\log N)$ for some
positive constant $C_0$.  We denote
by ${\cal P}_n$ the set of marked partitions of $\{1, \dots, n\}$.

\begin{Theo}
Fix $n \in \N$.  Let $\alpha = r \log(2N)/s$.  Let $p = e^{-\alpha}$.
Then, there exists a positive constant $C$ such that
$|P(\Theta = \pi) - Q_p(\pi)| \leq C/(\log N)$
for all $N$ and all $\pi \in {\cal P}_n$.
\label{logthm}
\end{Theo}

\noindent In this theorem, and throughout
the rest of the paper, $C$ denotes a positive
constant that may depend on $s$ but does not depend on $r$ or $N$.
The value of $C$ may change from line to line. 

A consequence of Theorem \ref{logthm} is that if
$\lim_{N \rightarrow \infty} r \log(2N)/s = \alpha$ for some
$\alpha \in (0, \infty)$
and $p = e^{-\alpha}$, then the distribution of $\Theta$ converges to $Q_p$
as $N \rightarrow \infty$.  However, the rate of convergence 
that the theorem gives is rather slow, and simulation results
of Barton (1998) and Durrett and Schweinsberg (2004a) show that the
approximation is not very accurate for realistic values of $N$.
Consequently, it is necessary to look for a better approximation.
Theorem \ref{log2thm} below gives an approximation with an error term
that is of order $1/(\log N)^2$ rather than $1/\log N$.  It follows from
the improved approximation that the error in Theorem \ref{logthm} is
actually of order $1/\log N$.

The motivation for the second approximation comes from the observation
that, at the beginning of the selective sweep, the number of $B$'s
can be approximated by a continuous-time branching process in which
each individual gives birth at rate $1$ and dies at rate $1-s$.
Some individials in this supercritical branching process will have
an infinite line of descent, meaning that they have descendants alive
in the population at all future times.  As we will show later, the
individuals with an infinite line of descent can be approximated
by a Yule process, a continuous-time branching process in which each
individual splits into two at a constant rate $s$.  Since our sample,
taken at the end of the selective sweep, comes from lineages that
have survived a long time, we can get a good approximation to the
genealogy by considering only individuals with an infinite line of
descent.  We will also show that, during the time when there are
exactly $k \geq 2$ lineages with an infinite line of descent,
the expected number of recombinations along these lineages is $r/s$.
For simplicity, we assume that the number of such recombinations is
always either $0$ or $1$.  Such a recombination causes individuals
descended from the lineage with the recombination to be traced back
to an ancestor at time zero different from descendants of the other $k-1$
lineages (and therefore to belong to a
different block of $\Theta$).  Well-known facts about the Yule process
(see e.g. Joyce and Tavar\'e, 1987) imply that when there are $k$
lineages, the fraction of individuals at the end of the sweep that are
descendants of a given lineage has approximately a beta($1,k-1$) distribution.
Furthermore, we will show that with probability $r(1-s)/(r(1-s) + s)$,
there is a recombination when there is only one individual with an
infinite line of descent, in which case none of the sampled lineages
will get traced back to the individual with the $B$ allele at time zero.

These observations motivate the definition of a class of marked
partitions of $\{1, \dots, n\}$, which we will use to approximate
the distribution of $\Theta$.  The construction resembles the
paintbox construction of exchangeable random partitions due to Kingman (1978).
To start the construction,
assume $0 < r < s$, and let $L$ be a positive integer.
Then let $(W_k)_{k=2}^L$ be independent random variables such that
$W_k$ has a Beta distribution with parameters $1$ and $k-1$.
Let $(\zeta_k)_{k=2}^L$ be a sequence of independent random variables such
that $P(\zeta_k = 1) = r/s$ and $P(\zeta_k = 0) = 1 - r/s$ for all
$k$.  As the reader might guess from the probabilities, $\zeta_k =1$
corresponds to a recombination when there are $k$ lineages with an
infinite line of descent.
For $k = 2, 3, \dots, L$, let $V_k = \zeta_k W_k$, and let
$Y_k = V_k \prod_{j=k+1}^L (1 - V_j)$ be the fraction of individuals
carried away by recombination. 
Let $Y_1 = \prod_{j=2}^L (1 - V_j)$.
Note that $\sum_{k=1}^L Y_k = 1$.  Finally, let $Q_{r,s,L}$ be the
distribution of the random
marked partition $\Pi$ of $\{1, \dots, n\}$ constructed
in the following way.  Define random variables $Z_1, \dots, Z_n$ to be
conditionally independent given $(Y_k)_{k=1}^L$ such that for
$i = 1, \dots, n$ and $j = 1, \dots, L$, we have
$P(Z_i = j|(Y_k)_{k=1}^L) = Y_j$.  Here the integers $i$ such
that $Z_i = k$ correspond to lineages that recombine when there
are $k$ members of the $B$ population with an infinite line of descent.
Then define $\Pi$ such that $i \sim_{\Pi} j$ if and only if $Z_i = Z_j$.
Independently of $(Z_i)_{i=1}^n$, we mark the block $\{i: Z_i = 1\}$ with
probability $s/(r(1-s) + s)$ and, with probability $r(1-s)/(r(1-s) + s)$,
we mark no block. When the block is
marked, the integers $i$ such that $Z_i = 1$ correspond to the lineages
that do not recombine and therefore can be traced back to the individual
that had the beneficial mutation at time zero; otherwise, they correspond to the lineages that recombine
when there is only one member of the $B$ population with an infinite line
of descent.

We are now ready to state our main approximation theorem,
which says that the distribution of $\Theta$ can be approximated
well by the distribution $Q_{r,s,L}$, where $L = \lfloor 2Ns \rfloor$,
and $\lfloor m \rfloor$ denotes the greatest integer less than or equal to $m$.
The choice of $L$ comes from the fact that
in a continuous-time branching process with births at rate $1$ and
deaths at rate $1-s$, each individual has an infinite line of descent
with probability
$s$.  Therefore, the number of such individuals at the end of the selective
sweep is approximately $L$.  

\begin{Theo}
Fix $n \in \N$ and let $L = \lfloor 2Ns \rfloor$.  Then, there exists a positive constant $C$ such that 
for all $N$ and all $\pi \in {\cal P}_n$
$$
|P(\Theta = \pi) - Q_{r,s,L}(\pi)| \leq C/(\log N)^2
$$
\label{log2thm}
\end{Theo}

Consider for concreteness $N=10,000$, a number commonly used for the
``effective size" of the human population. To explain the term in
quotes, we note that although there are now 6 billion humans,
our exponential population growth is fairly recent, so for many measures
of genetic variability the human population is the same as a
homogeneously mixing population of constant size 10,000.
When $N=10,000$, $\log N = 9.214$ and $(\log N)^2 = 84.8$, so
Theorem \ref{log2thm} may not appear at first glance to be a big improvement.
Two concrete examples however show that the improvement is dramatic.
In each case $N=10^4$ and $s=0.1$.  More extensive simulation results
comparing the two approximations are given in Durrett and Schweinsberg (2004a).

\mn
\begin{center}
\begin{tabular}{llllll}
& & $pinb$ & $p2inb$ & $p2cinb$ & $p1B1b$ \\
\hline
& Theorem \ref{logthm} & 0.1 &  0.01 & 0 & 0.18 \\
$r = .00106$ & Moran &  0.08203 & 0.00620 & 0.01826 & 0.11513 \\
& Theorem \ref{log2thm} & 0.08235 & 0.00627 & 0.01765 & 0.11687 \\
\hline
& Theorem \ref{logthm} & 0.4  & 0.16  &  0  &  0.48 \\
$r = .00516$& Moran & 0.33656 & 0.10567 & 0.05488 & 0.35201 \\
&Theorem \ref{log2thm} & 0.34065 & 0.10911 & 0.05100 & 0.36112 \\
\hline
\end{tabular}
\end{center}

\bigskip
Here $pinb$ is the probability that a lineage escapes the selective sweep.
The remaining three columns pertain to two lineages:
$p2inb$ is the probability that two lineages both escape the sweep but do
not coalesce, $p2cinb$ is the probability both lineages escape but coalesce
along the way, and $p1B1b$ is the probability one lineage escapes the sweep
but the other does not.  The remaining possibility is that neither lineage
escapes the sweep, but this probability can be computed by subtracting the
sum of the other three probabilities from one.  The first row in each
group gives the probabilities obtained from the approximation in
Theorem \ref{logthm}, and the third row gives the probabilities obtained
from the approximation in Theorem \ref{log2thm}.  The second row gives
the average of 10,000 simulation runs of the Moran model described earlier.
The values of the recombination rate $r$ were chosen in the two examples
to make the approximations to $pinb$ given by Theorem \ref{logthm}
equal to 0.1 and 0.4 respectively.  It is easy to see from the table
that the approximation from Theorem \ref{log2thm} is substantially more
accurate.  In particular, note that in the approximation given by Theorem
\ref{logthm}, two lineages never coalesce unless both can be traced back to
the individual with the beneficial mutation.  Consequently,
$p2cinb$ would be zero if this approximation were correct.  However, in
simulations, a significant percentage of pairs of lineages both coalesced
and escaped from the sweep, and this probability is approximated
very well by Theorem \ref{log2thm} in both examples.

The results in this paper are a first step in studying situations
in which, as proposed by Gillespie (2000),
selective sweeps occur at times of a Poisson process in a single locus
or distributed along a chromosome at different distances from the neutral
locus at which data have been collected.
It is well-known that in the Moran model when there are no advantageous
mutations, if we sample $n$ individuals and follow their ancestors backwards
in time, then when time is sped up by $2N$,
we get the coalescent process introduced by Kingman (1982).
It is known (see Durrett, 2002) that
selective sweeps require an average amount of time $(2/s) \log N$, so
when time is sped up by $2N$, the selective sweep occurs almost
instantaneously.  Durrett and Schweinsberg (2004b) show that Theorem
\ref{logthm} implies that if advantageous mutations
occur at times of a Poisson process then the ancestral processes 
converge as $N \rightarrow \infty$ to a coalescent with multiple collisions
of the type introduced by Pitman (1999) and Sagitov (1999).  At times of a 
Poisson process, multiple lineages coalesce simultaneously into one.
The more accurate approximation
in Theorem \ref{log2thm} suggests that a better approximation to the ancestral
process can be given by a coalescent with simultaneous multiple collisions.
These coalescent processes were studied by M{\"o}hle and Sagitov (2001) and
Schweinsberg (2000).  

Finally, it is important to emphasize that the results in this
paper are for the case of ``strong selection'', where the selective
advantage $s$ is $O(1)$.  There has also been considerable interest
in weak selection, where $Ns$ is assumed to converge to a limit as
$N \rightarrow \infty$, which means $s$ is $O(1/N)$.  In this case,
there is a diffusion limit as
$N \rightarrow \infty$.  For work in this direction that incorporates
the effect of recombination, see Donnelly and Kurtz (1999) and
Barton, Etheridge, and Sturm (2004).  Recently, Etheridge,
Pfaffelhuber, and Wakolbinger (2005) have shown that many of the
results in this paper carry over to the diffusion setting.
They assume that $Ns \rightarrow \alpha$ as $N \rightarrow \infty$,
so that they can work with a diffusion limit, and then obtain an
approximation to the distribution of the ancestral partition $\Theta$
which has an error of order $1/(\log \alpha)^2$ as $\alpha \rightarrow \infty$,
by using approximations to the genealogy similar to those used in the present paper.

\section{Overview of the Proofs}

Since the proofs of Theorems \ref{logthm} and \ref{log2thm} are rather long,
we outline the proofs in this section.  A precise definition of the
genealogy is given in subsection 2.1.  The proof of Theorem 1.1
is outlined in subsection 2.2.  In subsection 2.3, we describe the
coupling with a supercritical branching process and outline the
proof of Theorem 1.2.

\subsection{Precise definition of the genealogy}

We now define more precisely our
model of a selective sweep.  We construct a process
$M = (M_t)_{t=0}^{\infty}$.  The vector $M_t = (M_t(1), \dots, M_t(2N))$
contains the information about the population at the time of the $t$th
proposed replacement, and $M_t(i) = (A_t^0(i), \dots, A_t^{t-1}(i), B_t(i))$
contains the information about the ancestors of
the $i$th individual at time $t$.
For $0 \leq u \leq t - 1$, $A_t^u(i)$ is the individual at time $u$
that is the ancestor of the $i$th individual at time $t$, when we
consider the genealogy at the neutral locus. The final 
coordinate $B_t(i) = 1$ if the $i$th individual at time $t$
has the $B$ allele, and $B_t(i) = 0$  if this individual has the $b$ allele.
Note that this is a discrete-time process, but one can easily recover the
continuous-time description by replacing discrete time steps with independent
holding times, each having an exponential distribution with mean $1/2N$.

At time zero, only one of the chromosomes will have the $B$ allele.
We define a random variable $U$, which is uniform on the set
$\{1, \dots, 2N\}$, and we let $B_0(U) = 1$ and $B_0(i) = 0$ for
$i \neq U$.  We now define a collection of independent random variables
$I_{t,j}$ for $t \in \N$ and $j \in \{1, \dots, 5\}$.
For $j \in \{1, 2, 3\}$, the random variable
$I_{t,j}$ is uniform on $\{1, \dots, 2N\}$.  

\begin{itemize}

\item
$I_{t,1}$ will be the individual that dies at time $t$. 

\item
$I_{t,2}$ will be the parent of the new individual at time $t$.

\item
$I_{t,3}$ will be the other parent from whom the new chromosome will inherit 
its allele at the neutral locus if there is recombination.  

\item
$I_{t,4}$ will be an
indicator for whether a proposed disadvantageous change will be rejected,
so $P(I_{t,4} = 1) = s$ and $P(I_{t,4} = 0) = 1-s$.  

\item
$I_{t,5}$ will determine whether there is
recombination at time $t$, so $P(I_{t,5} = 1) = r$ and $P(I_{t,5} = 0) = 1-r$.

\end{itemize}
\clearpage

\begin{center}
\begin{picture}(320,250)
\put(16,25){\line(1,0){264}}
\put(16,225){\line(1,0){264}}
\put(15,25){\line(0,1){200}}
\put(280,25){\line(0,1){200}}
\put(15,12){0}
\put(275,12){$\tau$}
\put(210,45){$B$ population}
\put(25,210){$b$ population}
\put(285,147){$i$}
\put(280,150){\line(-1,0){82}}
\put(197,147){.}\put(197,144){.}\put(197,141){.}\put(197,138){.}\put(197,135){.}
\put(197,132){.}\put(197,129){.}\put(197,126){.}
\put(200,22){\line(0,1){6}}
\put(192,10){$G(i,j)$}
\put(285,123){$j$}
\put(280,126){\line(-1,0){127}}
\put(155,160){\line(-1,0){140}} 
\put(85,158){$\bullet$}
\put(80,168){$A^t_\tau(j)$}
\put(152,157){.}\put(152,154){.}\put(152,151){.}\put(152,148){.}\put(152,145){.}
\put(152,142){.}\put(152,139){.}\put(152,136){.}\put(152,133){.}\put(152,130){.}
\put(152,128){.}
\put(155,22){\line(0,1){6}}
\put(145,10){$R(i)$}
\put(285,73){$k$}
\put(280,76){\line(-1,0){151}}
\put(127,102){.}\put(127,99){.}\put(127,96){.}\put(127,93){.}\put(127,90){.}
\put(127,87){.}\put(127,86){.}\put(127,83){.}\put(127,80){.}\put(127,78){.}
\put(130,105){\line(-1,0){51}}
\put(77,102){.}\put(77,99){.}\put(77,96){.}\put(77,93){.}\put(77,90){.}\put(77,87){.}
\put(77,84){.}\put(77,81){.}\put(77,78){.}\put(77,75){.}\put(77,72){.}\put(77,69){.}
\put(77,66){.}\put(77,63){.}\put(77,60){.}\put(77,57){.}\put(77,54){.}\put(77,51){.}
\put(77,48){.}\put(77,45){.}\put(77,42){.}\put(77,39){.}\put(77,36){.}\put(77,33){.}
\put(77,30){.}
\put(80,27){\line(-1,0){65}}
\put(3,56){$J$}
\put(16,58){\line(1,0){91}}
\put(107,58){\line(0,-1){33}}
\put(105,10){$\tau_J$}
\put(15,26){.}
\put(16,26){.}\put(17,26){.}\put(18,26){.}\put(19,26){.}\put(20,26){.}
\put(21,26){.}\put(22,26){.}\put(23,26){.}\put(24,26){.}\put(25,26){.}
\put(26,26){.}\put(27,26){.}\put(28,26){.}\put(29,26){.}\put(30,26){.}
\put(31,26){.}\put(32,26){.}\put(33,27){.}\put(34,27){.}\put(35,27){.}
\put(36,27){.}\put(37,27){.}\put(38,27){.}\put(39,27){.}\put(40,27){.}
\put(41,27){.}\put(42,27){.}\put(43,28){.}\put(44,28){.}\put(45,28){.}
\put(46,28){.}\put(47,28){.}\put(48,28){.}\put(49,28){.}\put(50,28){.}
\put(51,29){.}\put(52,29){.}\put(53,29){.}\put(54,29){.}\put(55,29){.}
\put(56,30){.}\put(57,30){.}\put(58,30){.}\put(59,30){.}\put(60,30){.}
\put(61,31){.}\put(62,31){.}\put(63,31){.}\put(64,31){.}\put(65,32){.}
\put(66,32){.}\put(67,32){.}\put(68,33){.}\put(69,33){.}\put(70,33){.}
\put(71,34){.}\put(72,34){.}\put(73,34){.}\put(74,35){.}\put(75,35){.}
\put(76,35){.}\put(77,36){.}\put(78,36){.}\put(79,37){.}\put(80,37){.}
\put(81,38){.}\put(82,38){.}\put(83,39){.}\put(84,39){.}\put(85,40){.}
\put(86,40){.}\put(87,41){.}\put(88,42){.}\put(89,42){.}\put(90,43){.}
\put(91,44){.}\put(92,44){.}\put(93,45){.}\put(94,46){.}\put(95,46){.}
\put(96,47){.}\put(97,48){.}\put(98,49){.}\put(99,50){.}\put(100,51){.}
\put(101,52){.}\put(102,53){.}\put(103,54){.}\put(104,55){.}\put(105,56){.}
\put(106,57){.}\put(107,58){.}\put(108,59){.}\put(109,60){.}\put(110,61){.}
\put(111,62){.}\put(112,64){.}\put(113,65){.}\put(114,66){.}\put(115,68){.}
\put(116,69){.}\put(117,70){.}\put(118,72){.}\put(119,73){.}\put(120,75){.}
\put(121,76){.}\put(122,78){.}\put(123,79){.}\put(124,81){.}\put(125,83){.}
\put(126,84){.}\put(127,86){.}\put(128,88){.}\put(129,89){.}\put(130,91){.}
\put(131,93){.}\put(132,95){.}\put(133,97){.}\put(134,98){.}\put(135,100){.}
\put(136,102){.}\put(137,104){.}\put(138,106){.}\put(139,108){.}\put(140,110){.}
\put(141,112){.}\put(142,114){.}\put(143,116){.}\put(144,118){.}\put(145,120){.}
\put(146,122){.}\put(147,124){.}\put(148,126){.}\put(149,128){.}\put(150,130){.}
\put(151,132){.}\put(152,134){.}\put(153,136){.}\put(154,138){.}\put(155,140){.}
\put(156,142){.}\put(157,144){.}\put(158,146){.}\put(159,147){.}\put(160,149){.}
\put(161,151){.}\put(162,153){.}\put(163,155){.}\put(164,157){.}\put(165,158){.}
\put(166,160){.}\put(167,162){.}\put(168,164){.}\put(169,165){.}\put(170,167){.}
\put(171,169){.}\put(172,170){.}\put(173,172){.}\put(174,173){.}\put(175,175){.}
\put(176,176){.}\put(177,178){.}\put(178,179){.}\put(179,181){.}\put(180,182){.}
\put(181,183){.}\put(182,185){.}\put(183,186){.}\put(184,187){.}\put(185,188){.}
\put(186,189){.}\put(187,191){.}\put(188,192){.}\put(189,193){.}\put(190,194){.}
\put(191,195){.}\put(192,196){.}\put(193,197){.}\put(194,198){.}\put(195,199){.}
\put(196,200){.}\put(197,200){.}\put(198,201){.}\put(199,202){.}\put(200,203){.}
\put(201,204){.}\put(202,204){.}\put(203,205){.}\put(204,206){.}\put(205,206){.}
\put(206,207){.}\put(207,208){.}\put(208,208){.}\put(209,209){.}\put(210,209){.}
\put(211,210){.}\put(212,211){.}\put(213,211){.}\put(214,212){.}\put(215,212){.}
\put(216,212){.}\put(217,213){.}\put(218,213){.}\put(219,214){.}\put(220,214){.}
\put(221,215){.}\put(222,215){.}\put(223,215){.}\put(224,216){.}\put(225,216){.}
\put(226,216){.}\put(227,217){.}\put(228,217){.}\put(229,217){.}\put(230,217){.}
\put(231,218){.}\put(232,218){.}\put(233,218){.}\put(234,218){.}\put(235,219){.}
\put(236,219){.}\put(237,219){.}\put(238,219){.}\put(239,220){.}\put(240,220){.}
\put(241,220){.}\put(242,220){.}\put(243,220){.}\put(244,220){.}\put(245,221){.}
\put(246,221){.}\put(247,221){.}\put(248,221){.}\put(249,221){.}\put(250,221){.}
\put(251,221){.}\put(252,222){.}\put(253,222){.}\put(254,222){.}\put(255,222){.}
\put(256,222){.}\put(257,222){.}\put(258,222){.}\put(259,222){.}\put(260,222){.}
\put(261,222){.}\put(262,222){.}\put(263,223){.}\put(264,223){.}\put(265,223){.}
\put(266,223){.}\put(267,223){.}\put(268,223){.}\put(269,223){.}\put(270,223){.}
\put(271,223){.}\put(272,223){.}\put(273,223){.}\put(274,223){.}\put(275,223){.}
\put(276,223){.}\put(277,223){.}\put(278,223){.}\put(279,223){.}
\end{picture}
\end{center}

\mn
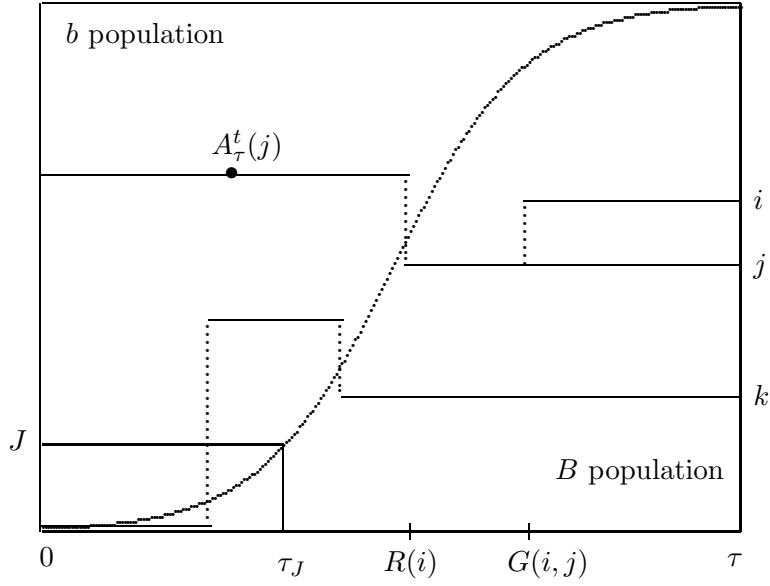
Figure 1. A picture to explain our notation. The lineages jump around as we move backwards
in time, but for simplicity we have only indicated the recombination events. Here as we work backwards in time 
$i$ and $j$ coalesce and then recombine into the $b$ population. Proposition
2.4 shows that this event has probability at most $C/\log N$.
Proposition 2.1 estimates the
probability of two recobminations as shown in lineage $k$.

\bigskip
Using these random variables we can construct the process in the obvious way.
Refer to Figure 1 for help with the notation.

\begin{enumerate}
\item If $B_{t-1}(I_{t,1}) = 1$, $B_{t-1}(I_{t,2}) = 0$, and
$I_{t,4} = 1$, then the population will be the same at time $t$ as
at time $t-1$ because the proposed replacement of a $B$ chromosome by a
$b$ chromosome is rejected.  In this case, for all $i = 1, \dots, 2N$ 
we define $B_t(i) = B_{t-1}(i)$, $A_t^{t-1}(i) = i$, and 
$A_t^u(i) = A_{t-1}^u(i)$ for $u = 0, \dots, t-2$.  

\item If we are not in the previous case and $I_{t,5} = 0$, then there
is no recombination at time $t$.  So, the individual $I_{t,1}$ dies,
and the new individual gets its alleles at both sites from $I_{t,2}$.
For $i \neq I_{t,1}$, define $B_t(i) = B_{t-1}(i)$,
$A_t^u(i) = A_{t-1}^u(i)$ for $u = 0, \dots, t-2$, and 
$A_t^{t-1}(i) = i$.  Let $B_t(I_{t,1}) = B_{t-1}(I_{t,2})$, 
$A_t^u(I_{t,1}) = A_{t-1}^u(I_{t,2})$ for $u = 0, \dots, t-2$,
and $A_t^{t-1}(I_{t,1}) = I_{t,2}$.

\item If we are not in either of the previous two cases, then there
is recombination at time $t$.  This means that the new individual
labeled $I_{t,1}$ gets a $B$ or $b$ allele from $I_{t,2}$ but gets its
allele at the neutral locus from $I_{t,3}$.  For $i \neq I_{t,1}$, define
$B_t(i) = B_{t-1}(i)$, $A_t^u(i) = A_{t-1}^u(i)$
for $u = 0, \dots, t-2$, and $A_t^{t-1}(i) = i$.
Let $B_t(I_{t,1}) = B_{t-1}(I_{t,2})$,
$A_t^u(I_{t,1}) = A_{t-1}^u(I_{t,3})$ for $u = 0, \dots, t-2$,
and $A_t^{t-1}(I_{t,1}) = I_{t,3}$.
\end{enumerate}

It will also be useful to have notation for the number of individuals with
the favorable allele.  For nonnegative integers $t$, define
$X_t = \# \{i: B_t(i) = 1\}$, where $\# S$ denotes the cardinality of the
set $S$.  For $J = 1, 2, \dots, 2N$, let $\tau_J = \inf\{t: X_t \geq J\}$
be the first time at which the number of $B$'s in the population
reaches $J$.  Let $\tau = \inf\{t: X_t \in \{0, 2N\}\}$ be the time
at which the $B$ allele becomes fixed in the population
(in which case $X_{\tau} = 2N$) or disappears (in which case
$X_{\tau} = 0$).
Since our main interest is in studying a selective sweep,
$P'$ and $E'$ will denote probabilities and expectations under
the unconditional law of $M$, and $P$ and $E$ will denote probabilities
and expectations under the conditional law of $M$ given $X_{\tau} = 2N$.
Likewise, $\Var$, and $\Cov$ will always refer to conditional
variances and covariances given $X_{\tau} = 2N$.

To sample $n$ individuals from the population at the time $\tau$ when the
selective sweep ends, we may take the individuals $1, \dots, n$ because
the distribution of genealogy of $n$ individuals does not depend on which
$n$ individuals are chosen.  Therefore, we can define $\Theta$ to be the
random marked partition of $\{1, \dots, n\}$ such that $i \sim_{\Theta} j$
if and only if the $i$th and $j$th individuals at time $\tau$ get their
allele at the neutral site from the same ancestor at time $0$, with the
marked block corresponding to the individuals descended from the
individual $U$, which had the beneficial mutation at time zero.
More formally, we have $i \sim_{\Theta} j$ if and only if
$A^0_{\tau}(i) = A^0_{\tau}(j)$ with the marked block being
$\{i: A^0_{\tau}(i) = U\}$ or, equivalently, $\{i: B_0(A^0_{\tau}(i)) = 1\}$.

\subsection{The first approximation}

Recall that Theorem 1.1 says that we can approximate $\Theta$
by flipping independent coins for each lineage, which come up heads
with probability $p$, to determine which lineages fail to escape
the selective sweep.  These lineages are then in one block of the
partition, because they are descended from the ancestor with the
beneficial mutation at time zero, while the other lineages do not
coalesce and correspond to singleton blocks of the partition.

The first step in establishing this picture is to calculate the
probability that one lineage escapes the selective sweep.  In the
notation above, we need to find $P(B_0(A^0_{\tau}(i)) = 0)$.
Define $R(i) = \sup\{t \geq 0: B_t(A^t_{\tau}(i)) = 0\}$, where
$\sup\emptyset = -\infty$.
If we work backwards in time, $R(i)$ is the first moment that the lineage of 
the neutral locus resides in the $b$ population. Note that it is 
possible to have $R(i) \geq 0$ and $B_0(A^0_{\tau}(i)) = 1$
if a lineage is affected by two recombinations, one taking it from the
$B$ population to the $b$ population, and another taking it back into
the $B$ population.  The next result shows that
the probability of this is small.

\begin{Prop}
$P( B_t(A_{\tau}^t(i)) = 1 \mbox{ for some }t \leq R(i) ) \leq C/(\log N)^2.$
\label{2rprop}
\end{Prop}

Proposition \ref{2rprop} implies that in the proofs of Theorems 1 and 2,
the probability that a lineage
escapes the selective sweep can be approximated by $P(R(i) \geq 0)$.  It will
also be useful to have an approximation of $P(R(i) \geq \tau_J)$ for
$J \geq 1$, which is the probability that
a given lineage escapes into the $b$ population after the time
when the number of $B$'s in the population reaches $J$.
The next result gives such an approximation.

\begin{Prop}
If $q_J = 1 - \exp \left( -\frac{r}{s} \sum_{k=J+1}^{2N}
\frac{1}{k} \right)$  then
\begin{displaymath}
P(R(i) \geq \tau_J) = q_J + O \bigg( \frac{1}{(\log N)^2} +
\frac{1}{(\log N) \sqrt{J}} \bigg).
\end{displaymath}
\label{1linprop}
\end{Prop}

\upat
Propositions \ref{2rprop} and \ref{1linprop} will be proved in
Section 3.

The next step is to consider two lineages.  We now need to 
consider not only recombination but also the possibility that the lineages
may coalesce, meaning that the alleles at the neutral site on the two
lineages are descended from the same ancestor at the beginning of the sweep.
Let $G(i,j)$ be the time that the $i$th and $j$th lineages coalesce.
More precisely, we define $G(i,j) = \sup\{t: A^t_{\tau}(i) = A^t_{\tau}(j)\}$
with $\sup\emptyset=-\infty$. Our first result regarding coalescence
shows that it is unlikely for two
lineages to coalesce at a given time unless both alleles at the
neutral site are descended from a chromosome with the $B$ allele at that time.

\begin{Prop}
$P( G(i,j) \geq 0 \mbox{ and }B_{G(i,j) + 1}(A_{\tau}^{G(i,j) + 1}
(i)) = 0 ) \leq C(\log N)/N$.
\label{bcoalprop}
\end{Prop}

\ni Next, we bound the probability that, if we trace two lineages back through
the selective sweep, the lineages coalesce and then escape from the sweep.

\begin{Prop}
$P(0 \leq R(i) \leq G(i,j)) \leq C/(\log N)$.
\label{crprop}
\end{Prop}

\ni Note that Proposition \ref{bcoalprop} says that, with high probability,
only lineages in the $B$ population merge, while Proposition
\ref{crprop} says that, in the first-order approximation, lineages that have
merged do not escape into the $b$ population.  Together, these results
will justify the approximation of $\Theta$ by a random partition in
which the only non-singleton block corresponds to lineages that
fail to escape the selective sweep.
The next result bounds the probability that two lineages coalesce
after time $\tau_J$. 

\begin{Prop}
Let $C' > 0$. If $J \leq C' N/(\log N)$, then $P(G(i,j) \geq \tau_J) \leq C/J$.
\label{coaltJprop}
\end{Prop}

\ni We prove Propositions \ref{bcoalprop}, \ref{crprop}, and
\ref{coaltJprop} in Section 4.

We now consider $n$ lineages.
To prove Theorem \ref{logthm}, we will need to show that the events
$\{R(1) \geq 0\}, \dots, \{R(n) \geq 0\}$ are approximately independent.
Let $K_t = \#\{i \in \{1, \dots, n\}: R(i) \geq t\}$.
If the events that the $n$ lineages escape the selective sweep after
time $t$ are approximately independent, then $K_t$ should have
approximately a binomial distribution.  The following proposition, which
we prove in Section 5, provides a binomial approximation to the distribution of
$K_{\tau_J}$.  Since $\tau_1 = 0$, the $J = 1$ case will be used in the
proof of Theorem \ref{logthm}, while the general case will help to prepare
us for the proof of Theorem \ref{log2thm}.

\begin{Prop}
Define $q_J$ as in Proposition \ref{1linprop}. If $J \leq C' N/(\log N)$, then 
$$\bigg|P(K_{\tau_J} = d) - \binom{n}{d} q_J^d (1-q_J)^{n-d}\bigg| \leq
\min \bigg\{ \frac{C}{\log N}, \frac{C}{J} \bigg\} +
\frac{C}{(\log N)^2}
\quad\hbox{for $d = 0, 1, \dots, n$.}
$$
\label{indprop}
\end{Prop}

\begin{proof}[Proof of Theorem \ref{logthm}]
Define a new partition $\Theta'$ of $\{1, \dots, n\}$ such that
$i \sim_{\Theta'} j$ if and only if $R(i) = R(j) = -\infty$.
We mark the block of $\Theta'$ consisting of $\{i: R(i) = -\infty\}$.
In words, only the lineages that recombine and hence stay in the $B$ population
are trapped by the sweep. To do this we observe

\begin{itemize}

\item
Proposition 2.1 implies that the probaility of two recombinations
affecting a lineage can be ignored.

\item
Proposition 2.3 says that we can ignore coalescence in the $b$ population.

\item
Propostion 2.4 says that the probability two lineages coalesce and then escape 
has small probaility.

\end{itemize}

The results above imply that $P(\Theta \neq \Theta') \leq C/(\log N)$.
Therefore, to prove Theorem \ref{logthm}, it suffices to show that
$|P(\Theta' = \pi) - Q_p(\pi)| \leq C/(\log N)$ for all marked partitions
$\pi$ of $\{1, \dots, n\}$.  It follows from Proposition \ref{indprop}
with $J = 1$ and the exchangeability of $\Theta'$ that
$|P(\Theta' = \pi) - Q_{1 - q_1}(\pi)| \leq C/(\log N)$
for all $\pi \in {\cal P}_n$. Using the definition of $q_1$ and $| \frac{d}{dx} e^{-x} | \le 1$
for $x \ge 0$ gives
$$|(1 - q_1) - p| = \bigg| \exp \bigg( -\frac{r}{s} \sum_{k = 2}^{2N}
\frac{1}{k} \bigg) - \exp \bigg( -\frac{r}{s} \log(2N) \bigg) \bigg|
\leq \frac{r}{s} \bigg| \sum_{k=2}^{2N} \frac{1}{k} - \log(2N) \bigg|
\leq \frac{C}{\log N},$$
and the theorem follows.
\end{proof}

\subsection{Branching process coupling and the second approximation}

We now work towards improving the approximation to the distribution of
$\Theta$ so that we can prove Theorem \ref{log2thm}.  To do this, we will
break the selective sweep into two stages.
Let $J = \lfloor (\log N)^a \rfloor$, where $a > 4$ is a fixed constant.
We will consider separately the time intervals
$[0,\tau_J)$ and $[\tau_J, \tau]$.  

\mn
{\bf Part 1.} $\Theta \approx \Theta_1 \approx \Theta_2$.

\medskip
We first establish that we can ignore coalescence involving
a lineage that escapes the sweep after time $\tau_J$.
Define a random marked partition
$\Theta_1$ of $\{1, \dots, n\}$ such that $i \sim_{\Theta_1} j$ if and only if 
$R(i) < \tau_J$, $R(j) < \tau_J$, and 
$A_{\tau}^0(i) = A_{\tau}^0(j)$.  Mark the block of $\Theta_1$ consisting of
$\{i: R(i) < \tau_J \mbox{ and } B_0(A_{\tau}^0(i)) = 1\}$.
Note that $\Theta_1 = \Theta$ unless, for some $i$ and $j$, we have
$R(i) \geq \tau_J$ and either $i \sim_{\Theta} j$ or
$B_0(A_{\tau}^0(i)) = 1$.  It follows from 
Propositions \ref{2rprop}, \ref{bcoalprop}, and \ref{coaltJprop} that
$P(\Theta \neq \Theta_1) \leq C/(\log N)^2$.  Thus, we may now
work with $\Theta_1$.

The next step is to approximate the distribution of $\Theta_1$. Let
$K_t = \{i \in \{1, \dots, n\}:
R(i) \geq t\}$, as defined before the statement of Proposition \ref{indprop}.
Define $m = n - \#K_{\tau_J}$ to be the number of lineages in
the $B$ population at time $\tau_J$. Proposition \ref{coaltJprop}
shows that lineages are unlikely to coalesce in $[\tau_J,\tau]$.
Relabel the lineages using an arbitrary bijective function
$f$ from $\{1, \dots, n\} \setminus
K_{\tau_J}$ to $\{1, \dots, m\}$.  

To describe the first stage of the selective sweep precisely,
we define, for each $m \leq J$, a new marked partition $\Psi_m$ of
$\{1, \dots, m\}$.  Let $\sigma_m$ be a random
injective map from $\{1, \dots, m\}$ to $\{i: B_{\tau_J}(i) = 1\}$
such that all $(J)_m = (J)(J-1) \dots (J-m+1)$ maps are equally likely.
Thus, $\sigma_m(1), \dots, \sigma_m(m)$ is a random sample from the $J$
individuals with the $B$ allele at time $\tau_J$.  Then define $\Psi_m$
such that that $i \sim_{\Psi_m} j$ if and only if
$A_{\tau_J}^0(\sigma_m(i)) = A_{\tau_J}^0(\sigma_m(j))$.
This means $i$ and $j$ are in the same block of $\Psi_m$ if and only if
the $\sigma_m(i)$th and $\sigma_m(j)$th individuals at time
$\tau_J$ inherited their allele
at the neutral locus from the same individual at the beginning of the sweep.
The block $\{i: B_0(A_{\tau_J}^0(\sigma_m(i))) = 1\}$ is marked.

Define $\Theta_2$ to be the marked
partition of $\{1, \dots, n\}$ such that $i \sim_{\Theta_2} j$ if and only if
$R(i) < \tau_J$, $R(j) < \tau_J$, and $f(i) \sim_{\Psi_m} f(j)$.
Let the marked block of $\Theta_2$ consist of all $i$ such that
$R(i) < \tau_J$ and $f(i)$ is in the marked block of $\Psi_m$. 
To compare $\Theta_1$ and $\Theta_2$, note that
$f(i) \sim_{\Psi_m} f(j)$ if and only if
$A_{\tau_J}^0(\sigma_m(f(i))) = A_{\tau_J}^0(\sigma_m(f(j)))$.
On the other hand, 
$A_{\tau}^0(i) = A_{\tau}^0(j)$ if and only if $A_{\tau_J}^0(A_{\tau}^{\tau_J}
(i)) = A_{\tau_J}^0(A_{\tau}^{\tau_J}(j))$.
For $i \neq j$, we have $P(A_{\tau}^{\tau_J}(i) = A_{\tau}^{\tau_J}(j))
\leq C/(\log N)^4$ by Proposition \ref{coaltJprop}.
By the strong Markov property, the genealogy of the process up to
time $\tau_J$ is independent of $K_{\tau_J}$.  From these observations
and the exchangeability of the model, 
it follows that for all $\pi \in {\cal P}_n$, we have
$|P(\Theta_1 = \pi) = P(\Theta_2 = \pi)| \leq C/(\log N)^4$.

\mn
{\bf Part 2.} $\Psi_m \approx \Upsilon_m \approx Q_{r,s, \lfloor Js \rfloor}(\pi)$

\medskip
Our next step is to understand
the distribution of $\Psi_m$.  The first step is to show that the
first stage of a selective sweep can be approximated by a branching process.
Recall that when the number of individuals with the favorable $B$ allele
is $k \ll 2N$, the rate of transitions that increase the number of
$B$ individuals from $k$ to $k + 1$ is $k(2N-k)/2N \approx k$, while the
rate of transitions that decrease the number of $B$ individuals from
$k$ to $k-1$ is $k(2N-k)(1-s)/2N \approx k(1-s)$.  Therefore, the
individuals with the $B$ allele follow approximately a continuous-time
branching process in which each individual gives birth at rate one and
dies at rate $1-s$.  Also, each new individual born with the $B$ allele
inherits the allele at the neutral site from its parent with 
probability $1-r$.  We can model this recombination by considering
a multi-type branching
process starting from one individual in which each new individual is the
same type as its parent with probability $1-r$ and is a new type, different
from any other member of the current population, with probability $r$.

Say that an individual in the branching process at time $t$ has an infinite
line of descent if it has a descendant in the population at time $u$ for
all $u > t$.  Otherwise, say the individual has a finite line of descent.
It is well-known that the process consisting only of the individuals with
an infinite line of descent is also a branching process.  This is discussed,
for example, in Athreya and Ney (1972).  For more recent work in this
direction, see O'Connell (1993) and Gadag and Rajarshi (1987, 1992). 
In Section 6, we will show that when the original branching process is a
continuous-time branching process with births at rate 1 and deaths at rate
$1-s$, the process consisting only of the individuals with an infinite
line of descent is a continuous-time branching process with no deaths in
which each individual gives birth at rate $s$.  That is, this process is
a Yule process with births at rate $s$.  The probability that a randomly
chosen individual has an infinite line of descent is $s$, so when the original
branching process has $J$ individuals, there are approximately $Js$
individuals with an infinite line of descent.  Furthermore, since the
past and future are independent by the Markov property, the genealogy
of a sample will not be affected if we sample only from the individuals
with infinite lines of descent.

In section 6, we justify these approximations.  This will lead to
a proof of the following proposition, which explains how the genealogy
of the first phase of a selective sweep can be approximated by the
genealogy of a continuous-time branching process.

\begin{Prop}
Consider a continuous-time multi-type branching process started with one
individual at time zero such that each individual gives birth at rate one
and dies at rate $1-s$.  Assume that each individual born has the same
type as its parent with probability $1-r$ and a new type with probability
$r$.  Condition this process to survive forever.  At the first
time at which there are $\lfloor Js \rfloor$ individuals with an infinite
line of descent, sample $m$ of the $\lfloor Js \rfloor$ individuals with
an infinite line of descent.  Define $\Upsilon_m$ to be the marked partition of
$\{1, \dots, m\}$ such that $i \sim_{\Upsilon_m} j$ if and only if the
$i$th and $j$th individuals in the sample have the same type, and the
marked block consists of the individuals with the same type as the
individual at time zero.  Then for all $\pi \in {\cal P}_m$, we have
$|P(\Psi_m = \pi) - P(\Upsilon_m = \pi)| \leq C/(\log N)^2$.
\label{branchprop}
\end{Prop}

Recall that in the introduction we constructed a random marked partition $\Pi$
with distribution $Q_{r, s, L}$, where $L = \lfloor 2Ns \rfloor$.
To compare this partition with $\Theta$, we will consider the construction
in two stages, just as we considered two stages of the selective sweep.
The first stage of the construction will involve the integers $i$ such that
$Z_i \leq \lfloor Js \rfloor$, and the second stage involves the integers $i$
such that $Z_i > \lfloor Js \rfloor$.  We think of $Z_i = k$ as meaning
that the $i$th lineage escapes the selective sweep at a time when there
are $k$ individuals in the Yule process (or, equivalently, $k$ lineages
in the branching process with an infinite line of descent).  
We use $\lfloor Js \rfloor$ as the boundary between the two stages
because, when the population size of the branching process reaches $J$,
there are approximately $Js$ individuals with an infinite line of
descent.

The next result compares the first stage of a selective sweep
to the random variables $Z_i$ such that $Z_i \leq \lfloor Js \rfloor$.

\begin{Prop}
There is a positive constant $C$ such that for all partitions $\pi$
of $\{1, \dots, n\}$, we have $|P(\Upsilon_n = \pi) -
Q_{r,s, \lfloor Js \rfloor}(\pi)| \leq C/(\log N)^2$.
\label{stage1prop}
\end{Prop}

\mn
{\bf Part 3.} $\Theta_2 \approx Q_{r, s, \lfloor Js \rfloor, q_J} \approx Q_{r,s,L}$ 

\medskip
Proposition 2.6 shows that the number of lineages that escape the sweep 
during $[\tau_J,\tau]$ has approximately a binomial distribution with 
success probability $q_J$. This motivates the following:

\begin{Dfn}
Let $r$, $s$, and $q$ be in $(0,1)$, and let $L$ be a positive integer.
Let $Q_{r, s, L, q}$ be the distribution of the random marked partition
$\Pi'$ of $\{1, \dots, n\}$
obtained as follows.  First, let $\Pi$ be a random marked partition
of $\{1, \dots, n\}$
with distribution $Q_{r, s, L}$.  Let $\xi_1, \dots, \xi_n$ be i.i.d.
random variables such that $P(\xi_i = 1) = q$ and $P(\xi_i = 0) = 1-q$.
Then say that $i \sim_{\Pi'} j$ if and only if $i \sim_{\Pi} j$ and
$\xi_i = \xi_j = 0$.  Mark the block of $\Pi'$ consisting of all integers
$i$ in the marked block of $\Pi$ such that $\xi_i = 0$.
\end{Dfn}

The next two propositions establish
the connection between the second stage of the construction of $\Pi$ and the
second stage of the selective sweep.  Proposition \ref{stage2prop1}
shows that it is unlikely to have $Z_i = Z_j$ if both are at least
$\lfloor Js \rfloor$, just as Proposition \ref{coaltJprop} shows that 
lineages are unlikely to coalesce 
during the second stage of a selective sweep.
Likewise, Proposition \ref{stage2prop2} shows that
the number of $Z_i$ greater than $\lfloor Js \rfloor$ has
approximately a binomial distribution, just as Proposition \ref{indprop}
shows that the number of lineages that escape the selective sweep
during the second stage has approximately a binomial distribution.

\begin{Prop} $P(Z_i = Z_j > \lfloor Js
\rfloor) \leq C/(\log N)^5$ for all $i \neq j$. 
\label{stage2prop1}
\end{Prop}

\begin{Prop}
Let $D = \#\{i: Z_i > \lfloor Js \rfloor\}$, and 
define $q_J$ as in Proposition 2.2.  Then
$$
\bigg|P(D = d) - \binom{n}{d} q_J^d (1 - q_J)^{n-d} \bigg| \leq 
\frac{C}{(\log N)^5} \qquad\hbox{for $d = 0, 1, \dots, n$}.
$$
\label{stage2prop2}
\end{Prop}

\ni
Propositions \ref{stage1prop}, \ref{stage2prop1}, and \ref{stage2prop2}
are proved in Section 7.  The proofs of Propositions \ref{stage2prop1} and
\ref{stage2prop2} are straightforward, but the proof of
Proposition \ref{stage1prop} is more difficult.  It involves considering
marked partitions $\pi$ with different numbers of blocks and doing
combinatorial calculations in each case.  

\begin{proof}[Proof of Theorem \ref{log2thm}]

By Propositions \ref{branchprop} and \ref{stage1prop}, we have
$|P(\Psi_n = \pi) - Q_{r,s, \lfloor Js \rfloor}(\pi)| \leq C/(\log N)^2$
for all $\pi \in {\cal P}_n$.  It follows from this fact, Proposition 
\ref{indprop}, and the construction of $\Theta_2$ that
$|P(\Theta_2 = \pi) - Q_{r,s, \lfloor Js \rfloor, q_J}(\pi)| \leq C/(\log N)^2$
for all $\pi \in {\cal P}_n$.  Also, by defining
$\xi_i = 1_{\{Z_i > \lfloor Js \rfloor\}}$ and applying Propositions
\ref{stage2prop1} and \ref{stage2prop2}, we see that
$|Q_{r,s, \lfloor Js \rfloor, q_J}(\pi) - Q_{r,s, L}(\pi)| \leq
C/(\log N)^5$ for all $\pi \in {\cal P}_n$.  This observation,
combined with the discussion in Part 1 of this subsection, completes the
proof of Theorem \ref{log2thm}.
\end{proof}

\section{Recombination of one lineage}

Our goal in this section is to prove Propositions \ref{2rprop} and
\ref{1linprop}, which pertain to the recombination probabilities
for a single lineage.  The strategy will be to study the process
$X = (X_t)_{t=0}^{\tau}$, which describes how the number of individuals with
the $B$ allele evolves during the selective sweep,
and then calculate recombination probabilities
conditional on the process $X$.  In subsection 3.1, we show that
the process $X$ behaves like an asymmetric random walk, and work out
some calculations that will be needed later.
We prove Proposition \ref{2rprop} in subsection 3.2 and Proposition
\ref{1linprop} in subsection 3.3.

\subsection{Random walk calculations}

Suppose $1 \leq X_{t-1} \leq 2N - 1$.  Then $X_t = X_{t-1} + 1$
if and only if $B_{t-1}(I_{t,1}) = 0$ and $B_{t-1}(I_{t,2}) = 1$.
Also, $X_t = X_{t-1} - 1$ if and only if 
$B_{t-1}(I_{t,1}) = 0$, $B_{t-1}(I_{t,2}) = 1$, and $I_{t,4} = 0$.
Otherwise, $X_t = X_{t-1}$.  It follows that, for $1 \leq k \leq 2N - 1$,
\begin{align}
&P'(X_t = X_{t-1} + 1|X_{t-1} = k) = \bigg( \frac{2N - k}{2N} \bigg)
\bigg( \frac{k}{2N} \bigg), \label{xt1} \\
&P'(X_t = X_{t-1} - 1|X_{t-1} = k) = \bigg( \frac{2N - k}{2N} \bigg)
\bigg( \frac{k}{2N} \bigg) (1-s), \label{xt2} \\
&P'(X_t = X_{t-1}|X_{t-1} = k) = 1 - \frac{(2-s)k(2N - k)}{(2N)^2}. \label{xt3}
\end{align}
Let $S_0 = 0$, and, for $m \geq 1$, let $S_m = \inf\{t > S_{m-1}:
X_t \neq X_{S_{m-1}}\}$ be the time of the $m$th jump.  
It follows from (\ref{xt1}) and (\ref{xt2})
that the process $(X_{S_m})_{m=0}^{\infty}$ is a random walk on
$\{0, 1, \ldots, 2N\}$ that starts at $1$, at each step moves to the
right with probability $1/(2-s)$ and to the left with probability
$(1-s)/(2-s)$, and is absorbed when it first reaches $0$ or $2N$.
A standard calculation for random walks (see e.g., section 3.1
of Durrett (2002)) gives the following result.

\begin{Lemma} Let $p(a,b,k) = P'(\inf\{s > t: X_s = b\} <
\inf\{s > t: X_s = a\}|X_t = k)$ be the probability that if the
number of $B$'s is $k$, then the number of $B$'s will reach 
$b$ before $a$.
For $0 \leq a < k < b \leq 2N$,  
$$
p(a,b,k) = \frac{1 - (1-s)^{k-a}}{1 - (1-s)^{b-a}}
\quad\hbox{and}\quad
P(X_{\tau} = 2N) = p(0,2N,1) = \frac{s}{1 - (1-s)^{2N}}.
$$
\label{survive}
\end{Lemma}

Given $1 \leq j \leq 2N-1$ and $1 \leq k \leq 2N - 1$, we define the
following quantities:
\begin{align}
\hbox{up jumps}\qquad & U_{k,j} = \#\{t \geq \tau_j: X_t = k \mbox{ and }X_{t+1} = k+1\} \nonumber \\
\hbox{down jumps}\qquad & D_{k,j} = \#\{t \geq \tau_j: X_t = k \mbox{ and }X_{t+1} = k-1\}  \nonumber \\
\hbox{holds}\qquad & H_{k,j} = \#\{t \geq \tau_j: X_t = k \mbox{ and }X_{t+1} = k\}  \nonumber \\
\hbox{total}\qquad & T_{k,j} = U_{k,j} + D_{k,j} + H_{k,j}  \nonumber
\end{align}
Also, let $U_k = U_{k,1}$, $D_k = D_{k,1}$, $H_k = H_{k,1}$, and
$T_k = T_{k,1}$.  The expected values
of these quantities are given in the following lemma.

\begin{Lemma}
Suppose $1 \leq j \leq 2N - 1$ and $1 \leq k \leq 2N - 1$.  Define
$$
q_k = \frac{p(k, 2N, k + 1)}{p(0, 2N, k+1)} =
\frac{s}{(1 - (1-s)^{2N - k})} \cdot \frac{(1 - (1 - s)^{2N})}{(1 - (1-s)^{k+1})} \ge s.
$$
Also, define $q_0 = 1$. Define $r_{k,j} = 1$ for $j \leq k$, and let
$r_{0,j} = 0$.  If $j > k$, let
$$
r_{k,j} = 1 - \frac{p(k, 2N, j)}{p(0, 2N, j)} =
1 - \frac{(1 - (1-s)^{j-k})}{(1 - (1-s)^{2N - k})}
\cdot \frac{(1 - (1-s)^{2N})}{(1 - (1-s)^j)} \le (1-s)^{j-k}.
$$  
Then $E[U_{k,j}] = r_{k,j}/q_k$.  Also, $E[D_{k,j}] = (1/q_{k-1}) - 1$
for $k > j$ and $E[D_{k,j}] = r_{k-1, j}/q_{k-1}$ for $k \leq j$.
Furthermore,
\begin{equation}
E[H_{k,j}] = E[U_{k,j} + D_{k,j}] \bigg( \frac{1}{2-s} \bigg) \frac{1}{\beta_k}
\leq \frac{\min\{(1-s)^{j-k}, 1\}}{s \beta_k},
\label{ehkj}
\end{equation}
where $\beta_k = k(2N-k)/ ( k^2 + (2N-k)^2 + sk(2N-k))$.
\label{udhlem}
\end{Lemma}

\noindent {\it Proof.}
First, suppose $k \geq j$.  On the event $\{X_{\tau} = 2N\}$,
we have $X_t = k$ and $X_{t+1} = k+1$ for some $t \geq \tau_j$.
Note that $P'(X_s > k \mbox{ for all }s > t|X_t = k+1) = p(k, 2N, k+1)$
for all $t$, so $P(X_s > k \mbox{ for all }s > t|X_t = k+1) =
p(k, 2N, k+1)/p(0, 2N, k+1) = q_k$.  It follows that the distribution
of $U_{k,j}$ is Geometric($q_k$), so $E[U_{k,j}] = 1/q_k$.
If instead $k < j$, then $P(X_t > k \mbox{ for all }t > \tau_j) =
p(k, 2N, j)/p(0, 2N,j)$.  Therefore, $P(T_{k,j} \geq 1) = r_{k,j}$.
It follows from the strong Markov property that, conditional on
$T_{k,j} \geq 1$, the distribution of $U_{k,j}$ is Geometric($q_k$),
so $E[U_{k,j}] = r_{k,j}/q_k$.

To find $E[D_{k,j}]$, note that if $k > j$ then $X$ takes a downward step
from $k$ to $k-1$ after each step from $k-1$ to $k$ except the last one,
so $D_{k,j} = U_{k-1,j} - 1$.  If $k \leq j$, then the number of steps after
$\tau_j$ from $k$ to $k-1$ is the same as the number of steps from $k-1$ 
to $k$, so $D_{k,j} = U_{k-1, j}$.  The formulas for $E[D_{k,j}]$ follow
immediately from these observations.

Let $p_k = P(X_t \neq X_{t-1}|X_{t-1} = k)$.
To prove (\ref{ehkj}), note that (\ref{xt3}) gives
$$
p_k = \frac{k(2N-k)(2-s)}{(2N)^2}.
$$
It follows that the conditional distribution of $T_{k,j}$ given
$U_{k,j}$ and $D_{k,j}$ is the same as the distribution of the sum of
$U_{k,j} + D_{k,j}$ independent
random variables with a Geometric($p_k$) distribution.
Therefore,
\begin{equation}
E[H_{k,j}] = E[T_{k,j}] - E[U_{k,j}] - E[D_{k,j}] = E[U_{k,j} + D_{k,j}] 
\bigg( \frac{1}{p_k} - 1 \bigg). \nonumber
\end{equation}
Straightforward algebraic manipulations give
$1/p_k - 1 = 1/[\beta_k(2-s)]$,
which implies the equality in (\ref{ehkj}).  To check the inequality in
(\ref{ehkj}), note that if $k > j$ then $$E[U_{k,j} + D_{k,j}] =
\frac{1}{q_k} + \frac{1}{q_{k-1}} - 1 \leq \frac{1}{s} + \frac{1}{s} - 1
= \frac{2-s}{s}$$ and if $k \leq j$ then $$E[U_{k,j} + D_{k,j}] \leq
\frac{(1-s)^{j-k}}{q_k} + \frac{(1-s)^{j-k+1}}{q_{k-1}} \leq
\frac{(1-s)^{j-k}}{s}\big( 1 + (1-s) \big) = \frac{(2-s)(1-s)^{j-k}}{s}.
\hphantom{xx}\Box$$

\bigskip
We will now calculate the probability that the ancestor at time $t$ has 
the opposite $B$ or $b$ allele from the ancestor at time $t-1$,
given that $X_{t-1} = k$ and $X_t = l$, where $1 \leq k \leq 2N-1$,
$1 \leq l \leq 2N$, and $|k - l| \leq 1$. All of these recombination 
probabilities are the same under $P'$ and $P$ because of the conditioning on
$X_{t-1}$ and $X_t$.  We define
\begin{align}
p_B^r(k, l) &= P(B_{t-1}(A_t^{t-1}(i)) = 0|X_{t-1} = k, X_t = l,
B_t(i) = 1), \nonumber \\
p_b^r(k, l) &= P(B_{t-1}(A_t^{t-1}(i)) = 1|X_{t-1} = k, X_t = l,
B_t(i) = 0). \nonumber
\end{align}

\begin{Lemma} We have
\begin{align}
& p_B^r(k, k-1) = p_b^r(k, k+1) = 0, \nonumber \\
& p_B^r(k, k+1) = \frac{r(2N-k)}{(k+1)(2N)}, \quad
p_b^r(k, k-1) = \frac{rk}{(2N-k+1)(2N)},  \nonumber \\
& p_B^r(k,k) = p_b^r(k, k) = \frac{rk(2N-k)}{2N[k^2 + (2N-k)^2 + sk(2N-k)]} = \frac{r \beta_k}{2N}.
\nonumber
\end{align}
\label{osrlem}
\end{Lemma}

\begin{proof}
We will prove three of the six results; the others are similar.
If $X_{t-1} = k$ and $X_t = k+1$, then the new individual born
at time $t$ has a $B$ allele.  Therefore, if $B_t(i) = 0$
then $B_{t-1}(A_t^{t-1}(i)) = 0$, so 
$p_b^r(k, k+1) = 0$.  Suppose instead $B_t(i) = 1$.
Then, $B_{t-1}(A_t^{t-1}(i)) = 0$ if and only if $I_{t,1} = i$
(meaning that the $i$th individual is the new one born),
$I_{t,5} = 1$ (meaning that there is recombination),
and $B_{t-1}(I_{t,3}) = 0$ (meaning that the new individual
gets its allele at the neutral site from the member of the
$b$ population).  Conditional on 
$X_{t-1} = k$, $X_t = k+1$, and $B_t(i) = 1$, the probabilities of
$I_{t,1} = i$, $I_{t,5} = 1$, and 
$B_{t-1}(I_{t,3}) = 0$ are $1/(k+1)$, $r$, and $(2N-k)/2N$
respectively.  Multiplying them
gives the expression for $p_B^r(k, k+1)$.
To calculate $p_B^r(k, k)$ we use the fact that, conditional on
$X_{t-1} = k$ and $X_t = k$, the probability that
$B_{t-1}(I_{t,1}) = B_{t-1}(I_{t,2}) = 1$ is
$k^2/[k^2 + (2N-k)^2 + sk(2N-k)]$. Multiplying by $1/k$, $r$,
and $(2N-k)/2N$ gives $p_B^r(k,k)$.  
\end{proof}

\subsection{Bounding the probability of two recombinations}

We now begin working towards a proof of Proposition \ref{2rprop},
which shows that it is unlikely that a lineage will go from
the $B$ population to the $b$ population and then back to
the $B$ population because of two recombination events.
We begin by proving two simple lemmas.  Lemma \ref{rkslem} bounds
the probability that the number of
individuals with the $B$ allele is $k$ at the recombination time $R(i)$.
Lemma \ref{aNlem} is a useful deterministic result,
which can be proved easily by splitting the sum
into terms with $j \leq N/2$ and $j > N/2$.

\begin{Lemma}
We have $P(X_{R(i)} = k) \leq r/ks$.
\label{rkslem}
\end{Lemma}

\mn
{\it Proof.}
Considering the cases $X_{R(i)+1} = k+1$ and
$X_{R(i)+1} = k$ and using Lemmas \ref{udhlem} and \ref{osrlem},
\begin{align}
\hphantom{xxxxxxxxxxx}
P(X_{R(i)} = k) &\leq
p_B^r(k, k+1) E[U_k] + p_B^r(k, k) E[H_k] \nonumber \\
&\leq \frac{r(2N - k)}{(k+1)(2Ns)} + \frac{r}{2Ns} \leq
\frac{r(2N - k) + rk}{2Nks} = \frac{r}{ks}. 
\hphantom{xxxxxxxxxxx}\Box \nonumber
\end{align}

\begin{Lemma}
If $a > 1$, there is a $C$ depending on $a$ but not on $N$ so that
$\sum_{j=1}^N a^j/j \leq C a^N/N$.
\label{aNlem}
\end{Lemma}

\mn{\it Proof of Proposition \ref{2rprop}.}  Denote the time
of the second recombination event by
$R_2(i) = \sup\{t \leq R(i): B_t(A^t_{\tau}(i)) = 1\}$,
where $\sup\emptyset=-\infty$.
Our goal is to show $P(R_2(i) \geq 0) \leq C/(\log N)^2$.  
Note that by symmetry, the conditional distribution of $(X_t)_{t=0}^{\tau - 1}$
given $X_{\tau} = 2N$ is the same as the conditional distribution of
$(2N - X_{\tau - t})_{t=1}^{\tau}$ given $X_{\tau} = 2N$.  It follows from
this fact and the strong Markov property that
\begin{gather}
E[\#\{t < R(i): X_t = k \mbox{ and }X_{t+1} = k+1|X_{R(i)} = j\}] =
E[U_{2N-k-1, 2N-j}], \nonumber \\
E[\#\{t < R(i): X_t = k \mbox{ and }X_{t+1} = k-1|X_{R(i)} = j\}] =
E[D_{2N-k+1, 2N-j}], \nonumber \\ 
E[\#\{t < R(i): X_t = k \mbox{ and }X_{t+1} = k|X_{R(i)} = j\}] =
E[H_{2N-k, 2N-j}]. \nonumber
\end{gather}
Therefore, by Lemmas \ref{udhlem} and \ref{osrlem},
\begin{align}
P(X_{R_2(i)} = k|X_{R(i)} = j) &\leq
p_b^r(k, k-1) E[D_{2N-k+1, 2N-j}] + p_b^r(k,k) E[H_{2N-k, 2N-j}] \nonumber \\
&\leq \frac{rk}{(2N-k+1)(2Ns)} \min\{(1-s)^{k-j}, 1\} +
\frac{r}{2Ns} \min\{(1-s)^{k-j}, 1\} \nonumber \\
&\le \frac{r}{(2N-k)s} \min\{ (1 - s)^{k-j}, 1\}. \nonumber
\end{align}
Using Lemma \ref{rkslem},
\begin{align}
P(R_2(i) \geq 0) &\leq \sum_{j=1}^{2N-1} \frac{r}{js} \bigg( \sum_{k=1}^{2N-1}
\frac{r \min\{ (1 - s)^{k-j}, 1\}}{(2N-k)s} \bigg) \nonumber \\
&= \frac{r^2}{s^2} \sum_{j=1}^{2N-1} \frac{1}{j} \bigg(
\sum_{k=j}^{2N-1} \frac{(1-s)^{k-j}}{2N-k} +
\sum_{k=1}^{j-1} \frac{1}{2N-k} \bigg).
\label{jk2terms}
\end{align}
Since $r^2/s^2 \leq C/(\log N)^2$, it suffices to show that the sum on the
right-hand side of (\ref{jk2terms}) is bounded as
$N \rightarrow \infty$.  We will handle the two terms separately.
For the first term, we change variables $\ell = 2N-k$ and use Lemma \ref{aNlem} to get the bound
\begin{align}
\sum_{j=1}^{2N-1} \frac{1}{j} \bigg( \sum_{k=j}^{2N-1}
\frac{(1-s)^{k-j}}{2N-k} \bigg) 
&= \sum_{j=1}^{2N-1} \frac{(1-s)^{2N-j}}{j}
\bigg( \sum_{\ell =1}^{2N-j} \bigg( \frac{1}{1-s} \bigg)^\ell \frac{1}{\ell} \bigg)
\nonumber \\
&\le C \sum_{j=1}^{2N-1} \frac{1}{j(2N-j)} \leq \frac{2C}{N} \sum_{j=1}^{N} \frac{1}{j} \leq \frac{2C(1 + \log N)}{N}.
\label{logNNb}
\end{align}
The second term in the sum on the right-hand side of (\ref{jk2terms}) can be bounded by
$$
\hphantom{xxxxxxxxxxx}
\sum_{j=1}^N \frac{1}{j} \bigg( \frac{j}{N} \bigg) +
\sum_{j=N+1}^{2N-1} \frac{1}{N} \bigg( \sum_{\ell=2N-j}^{2N-1}
\frac{1}{\ell} \bigg) = 1 + \frac{1}{N} \sum_{\ell=1}^{2N-1} \frac{1}{\ell}
\sum_{j=2N-l}^{2N-1} (1) \leq 3.
\hphantom{xxxxxxxxxxx}\Box
$$

\subsection{Estimating the recombination probability}

Our next goal is to prove Proposition \ref{1linprop},
which gives an approximation for $P(R(i) \geq \tau_J)$.
The idea behind the proof is that every time there is a
change in the population, there is some probability that a
lineage will escape the selective sweep at that time,
given that it has not previously escaped.  Since the individual
probabilities are small, if they sum to $S$, we will be able
to approximate the probability that the lineage never escapes
by $e^{-S}$.  It will be easier to work with conditional
escape probabilities given $X$, so to justify the approximation
it will necessary be to show that the sum of the conditional
probabilities has low variance.

For $1 \leq t \leq \tau$, let $\theta_t = p_B^r(X_{t-1}, X_t)$.
Now, $\theta_t$ is the conditional probability, given $X$,
that a lineage escapes at time $t$ if it has not previously
escaped, so we have
\begin{equation}
P(R(i) \geq \tau_J|X) = 1 - \prod_{t=\tau_J+1}^{\tau}
[1 - p^r_B(X_{t-1}, X_t)] = 1 - \prod_{t=\tau_J+1}^{\tau} (1 - \theta_t).
\label{defHt}
\end{equation}
To estimate the probability that a lineage escapes after time
$\tau_J$, we will consider the sum of these conditional probabilities,
which we denote by $\eta_J = \sum_{t = \tau_J+1}^{\tau} \theta_t$.
The next lemma shows that to 
estimate $P(R(i) \geq \tau_J)$ to within an error of $O((\log N)^{-2})$,
it suffices to calculate $E[e^{-\eta_J}]$.

\begin{Lemma}
For all $J$, we
have $\big| P(R(i) \geq \tau_J) - (1 - E[e^{-\eta_J}]) \big| \leq
C/(\log N)^2.$
\label{thetalemnew}
\end{Lemma}

\begin{proof}
It follows from the Poisson approximation on p. 140 of Durrett (1996) that
\begin{equation}
|P(R(i) \geq \tau_J|X) - (1 - e^{-\eta_J})|
\leq \sum_{t = \tau_J + 1}^{\tau} \theta_t^2.
\label{poisapp}
\end{equation}
By taking expectations, we get
$$\bigg| P(R(i) \geq \tau_J) - (1 - E[e^{-\eta_J}]) \bigg|  \leq
E \bigg[ \sum_{t=\tau_J+1}^{\tau} \theta_t^2 \bigg].$$
It now remains to bound
$E \big[ \sum_{t=1}^{\tau} \theta_t^2 \big].$
By Lemma \ref{osrlem},
$$\sum_{t=1}^{\tau} \theta_t^2 = 
\sum_{k=1}^{2N-1} \bigg( U_k \frac{r^2(2N-k)^2}{(k+1)^2(2N)^2}
+ H_k \frac{r^2 \beta_k^2}{(2N)^2} \bigg).$$
Therefore, by Lemma \ref{udhlem},
\begin{align}
E \bigg[ \sum_{t=1}^{\tau} \theta_t^2 \bigg] &\leq
\sum_{k=1}^{2N} \bigg( \frac{r^2(2N-k)^2}{s(k+1)^2(2N)^2} +
\frac{r^2 \beta_k}{(2N)^2 s} \bigg) \nonumber \\
&\leq \frac{r^2}{s} \sum_{k=1}^{2N} \bigg( \frac{1}{(k+1)^2} +
\frac{1}{(2N)^2} \bigg) \leq Cr^2 \leq \frac{C}{(\log N)^2},
\label{thetatexp}
\end{align}
which completes the proof.
\end{proof}

The next result will allow is to work with a truncated version
of the sum.

\begin{Lemma} If $\eta_J' = \sum_{t = \tau_J+1}^{\tau}
\theta_t 1_{\{X_{t-1} \geq J\}}$, then $E[\eta_J - \eta_J'] \leq C/J(\log N)$.
\label{etaJ}
\end{Lemma}

\mn{\it Proof.}
Using Lemmas \ref{udhlem}, \ref{osrlem}, and \ref{aNlem}, we get
\begin{align}
\hphantom{xxxxxxxxxxxxxx}
E[\eta_J - \eta_J'] &= \sum_{k=1}^{J-1}
\big(p_B^r(k, k+1) E[U_{k,J}] + p_B^r(k,k) E[H_{k,J}] \big) \nonumber \\
& \le \sum_{k=1}^{J-1} \frac{r(2N-k)}{(k+1)(2N)} \cdot \frac{(1-s)^{J-k}}{s}
+ \frac{r \beta_k}{2N} \cdot \frac{(1-s)^{J-k}}{s \beta_k} \nonumber \\
&\leq \sum_{k=1}^{J-1} (1-s)^{J-k} \bigg( \frac{r}{ks} \bigg) =
\frac{r}{s} \sum_{k=1}^{J-1} \frac{1}{k} \bigg( \frac{1}{1-s} \bigg)^{k-J}
\leq \frac{Cr}{sJ}. \hphantom{xxxxxxxxx}\Box \nonumber
\end{align}

We will work with
$\eta_J'$ rather than $\eta_J$ because we can obtain a rather precise
estimate of its expected value, which is given in the
next lemma.  We will also be able to obtain a bound on its variance,
which will enable us to approximate $E[e^{-\eta_J'}]$ by $e^{-E[\eta_J']}$.

\begin{Lemma}
$E[\eta_J'] = \frac{r}{s} \sum_{k=J+1}^{2N} \frac{1}{k} +
O \left( \frac{1}{N} + \frac{(1-s)^J}{J \log N} \right).$
\label{ewjlem}
\end{Lemma}

\begin{proof}
It follows from Lemma \ref{udhlem} and a straightforward calculation that
$$E[H_k] = \bigg( \frac{1}{q_k} + \frac{1}{q_{k-1}} - 1 \bigg)
\bigg( \frac{1}{\beta_k (2-s)} \bigg) = \frac{1}{s \beta_k}
\bigg( \frac{(1 - (1-s)^k)(1 - (1-s)^{2N-k})}{1 - (1-s)^{2N}} \bigg).$$
Therefore,
\begin{align}
E[\eta_J'&] = \sum_{k=J}^{2N-1} \bigg( \frac{r(2N-k)}{(k+1)(2N)} E[U_k]
+ \frac{r \beta_k}{2N} E[H_k] \bigg) \nonumber \\
&= \sum_{k=J}^{2N-1}
\bigg( \frac{r(2N-k)(1 - (1-s)^{k+1})(1 - (1-s)^{2N-k})}{(k+1)
(2Ns)(1 - (1-s)^{2N})} + \frac{r(1 - (1-s)^k)(1 - (1-s)^{2N-k})}{(2Ns)
(1 - (1-s)^{2N})} \bigg) \nonumber \\
&= \frac{r}{s} \sum_{k=J}^{2N-1} \bigg( \frac{1 - (1-s)^{2N-k}}{1 -
(1-s)^{2N}} \bigg) \bigg( \frac{(2N-k)(1 - (1-s)^{k+1})}{(2N)(k+1)} +
\frac{1 - (1-s)^k}{2N} \bigg). \nonumber
\end{align}
Now
\begin{align}
\sum_{k=J}^{2N-1} &\bigg( 1 - \frac{1 - (1-s)^{2N-k}}{1 -
(1-s)^{2N}} \bigg) \bigg( \frac{(2N-k)(1 - (1-s)^{k+1})}{(2N)(k+1)} +
\frac{1 - (1-s)^k}{2N} \bigg) \nonumber \\
&\leq \sum_{k=J}^{2N-1} (1-s)^{2N-k} \bigg( \frac{1}{k} + \frac{1}{2N} \bigg)
\leq \sum_{k=1}^N (1-s)^N \bigg( \frac{2}{k} \bigg) +
\sum_{k=N+1}^{2N} (1-s)^{2N-k} \bigg( \frac{2}{N} \bigg) \nonumber \\
&\leq 2(1 + \log N)(1 - s)^N + \frac{2}{Ns} \leq \frac{C}{N}. \nonumber
\end{align}
Therefore,
$$E[\eta_J'] = \frac{r}{s} \sum_{k=J}^{2N-1} \bigg( \frac{(2N-k)
(1 - (1-s)^{k+1})}{(2N)(k+1)} + \frac{1 - (1-s)^k}{2N} \bigg) +
O \bigg( \frac{1}{N} \bigg).$$
Also, note that $$\sum_{k=J}^{2N-1} \frac{(1 - (1-s)^{k+1}) -
(1 - (1-s)^k)}{2N} = \sum_{k=J}^{2N} \frac{(1-s)^ks}{2N} \leq
\frac{1-s}{2N}.$$  Therefore, since $r \leq C_0 /\log N$,
\begin{align}
E[\eta_J'] &= \frac{r}{s} \sum_{k=J}^{2N-1} \bigg( \frac{2N-k}{(2N)(k+1)} +
\frac{1}{2N} \bigg) (1 - (1-s)^{k+1})+ O \bigg( \frac{1}{N} \bigg) \nonumber \\
&= \frac{r}{s} \left( \frac{2N+1}{2N} \right) \sum_{k=J}^{2N-1}
\frac{1 - (1-s)^{k+1}}{k+1} + O \bigg( \frac{1}{N} \bigg) = \frac{r}{s}
\sum_{k=J+1}^{2N} \frac{1 - (1-s)^k}{k} + O \bigg( \frac{1}{N} \bigg).
\label{laster}
\end{align}
Since $$\frac{r}{s} \sum_{k=J+1}^{2N} \frac{(1-s)^{k}}{k} \leq
\frac{r}{s(J+1)} \sum_{k=J+1}^{\infty} (1-s)^k = 
\frac{r(1-s)^{J+1}}{s^2(J+1)},$$
the desired result follows from (\ref{laster}).
\end{proof}

The key remaining step is to bound $\Var(\eta_J')$.  The necessary
bound is given in Lemma \ref{varwj}.  The proof uses
Lemma \ref{varcov}, which can easily be proved by conditioning
on $M$ and $N$.

\begin{Lemma}
Suppose $(X_i)_{i=1}^{\infty}$ and $(Y_i)_{i=1}^{\infty}$ are independent
i.i.d. sequences such that $E[X_1] = \mu$ and $E[Y_1] = \gamma$.
Suppose $M$ and $N$ are integer-valued random variables that are independent
of these sequences.  Then
$\Cov(X_1 + \dots + X_M, Y_1 + \dots + Y_N) = \mu \gamma \Cov(M, N)$.
\label{varcov}
\end{Lemma}

\begin{Lemma}
There exists a constant $C$ such that $\Var(\eta_J') \leq C/J(\log N)^2.$
\label{varwj}
\end{Lemma}

\proof Let 
\begin{align}
&a_k = \frac{2N-k}{(k+1)(2N)} \leq \frac{1}{k}, \label{bak} \\
& b_k = \frac{k(2N-k)}{2N[k^2 + (2N-k)^2 + sk(2N-k)]} \leq
\frac{k(2N-k)}{2N^3}. \label{bbk}
\end{align}
Then
$\eta_J' = \sum_{t = \tau_J + 1}^{\tau} \theta_t 1_{\{X_{t-1} \geq J\}} =
r \sum_{k=J}^{2N-1} (a_k U_k + b_k H_k).$
For any random variables $X$ and $Y$,
\begin{align}
\Var(X + Y) &= \Var(X) + \Var(Y) + 2\Cov(X, Y) \nonumber \\
&\leq \Var(X) + \Var(Y) + 2 \sqrt{\Var(X) \Var(Y)} \leq
4 \max\{\Var(X), \Var(Y)\}. \nonumber 
\end{align}
Therefore,
\begin{equation}
\Var(\eta_J') \leq 4r^2 \max \bigg\{ \Var
\bigg( \sum_{k=J}^{2N-1} a_k U_k \bigg),
\Var \bigg( \sum_{k=J}^{2N-1} b_k H_k \bigg) \bigg\}.
\label{vwj0}
\end{equation}
We will bound $\Var(\sum_{k=J}^{2N-1} a_k U_k)$ and
$\Var(\sum_{k=J}^{2N-1} b_k H_k)$ by $C/J$, which will prove the lemma.

To bound $\Var(\sum_{k=J}^{2N-1} a_k U_k)$, we will need to
bound $\Cov(U_k, U_l)$.  To do this, we will break up $U_l$
into jumps from $l$ to $l+1$ that occur before the last visit to $k$
and those that occur after the last visit to $k$.  More formally,
let $\zeta_k = \sup\{t: X_t = k\}$.  If $k \leq l$, then
$U_l = U_{k,l}' + {\bar U}_{k,l}$, where
\begin{gather}
U_{k,l}' = \#\{t \geq \zeta_k: X_t = l \mbox{ and }X_{t+1} = l+1\}, \nonumber \\
{\bar U}_{k,l} = \#\{t < \zeta_k: X_t = l \mbox{ and }X_{t+1} = l+1\}.
\nonumber
\end{gather}
The processes $(X_t)_{0 \leq t \leq \zeta_k}$ and
$(X_t)_{\zeta_k \leq t \leq \tau}$ are independent.  Therefore,
$U_k$ and $U_{k,l}'$ are independent, and ${\bar U}_{k,l}$ and
$U_{k,l}'$ are independent.  As observed in the proof of Lemma \ref{udhlem},
$U_l$ has a Geometric($q_l$) distribution.  Likewise, note that
$P'(X_s > l \mbox{ for all }s > t|X_t = l+1) = p(l, 2N, l+1)$ and
$P'(X_s > k \mbox{ for all }s > t|X_t = l+1) =
p(k, 2N, l+1)$.  Therefore, $$P(X_s > l \mbox{ for all }s > t| X_t = l+1,
X_s > k \mbox{ for all }s > t) = \frac{p(l, 2N, l+1)}{p(k, 2N, l+1)}.$$
It follows that if we let $v_{k,l} = p(l, 2N, l+1)/p(k, 2N, l+1)$,
then $U_{k,l}'$ has a Geometric($v_{k,l}$) distribution.
Using Lemmas 3.1 and 3.2 and
the fact that $q_l = p(l,2N,l+1)/p(0,2N,l+1)$, we have
\begin{align}
\frac{1}{q_l} - \frac{1}{v_{k,l}}
&= \frac{1 - (1-s)^{2N-l}}{s} \bigg(
\frac{1 - (1-s)^{l+1}}{1 - (1-s)^{2N}} - \frac{1 - (1-s)^{l+1-k}}{1 - (1-s)^{2N-k}} \bigg) \nonumber \\
&\leq \frac{1}{s} (1 - (1 - (1-s)^{l+1-k})) = \frac{(1-s)^{l+1-k}}{s}.
\label{vwj1}
\end{align}
Also, $\Var(U_l) = \Var(U_{k,l}') + \Var({\bar U}_{k,l})$ because
${\bar U}_{k,l}$ and $U_{k,l}'$ are independent.  Therefore, if
$J \leq k \leq l < 2N$, then by the formula for the variance of a geometric
distribution,
\begin{align}
\Var({\bar U}_{k,l}) &= \Var(U_l) - \Var(U_{k,l}') =
\frac{1 - q_l}{q_l^2} - \frac{1 - v_{k,l}}{v_{k,l}^2} 
\nonumber \\
&= \bigg( \frac{1}{q_l} + \frac{1}{v_{k,l}} - 1 \bigg)
\bigg( \frac{1}{q_l} - \frac{1}{v_{k,l}} \bigg) \leq
\frac{2}{s} \cdot \frac{(1-s)^{l-k}}{s},
\label{vwj2}
\end{align}
where the inequality uses (\ref{vwj1}) and the facts that
$q_l \geq s$ and $v_{k,l} \geq s$.  Also,
\begin{equation}
\Var(U_l) = \frac{1 - q_l}{q_l^2} \leq \frac{1}{s^2}.
\label{vwj3}
\end{equation}
Since $U_k$ and $U_{k,l}'$ are independent, it follows from
(\ref{vwj2}) and (\ref{vwj3}) that if $k \leq l$, then
\begin{align}
\Cov(U_k, U_l) &= \Cov(U_k, U_{k,l}' + {\bar U}_{k,l}) =
\Cov(U_k, {\bar U}_{k,l}) \nonumber \\
&\leq \sqrt{\Var(U_k) \Var({\bar U}_{k,l})} \leq
\frac{\sqrt{2}}{s^2} (1-s)^{(l-k)/2}.
\label{vwj4}
\end{align}
Using (\ref{vwj4}) and (\ref{bak}), we calculate
\begin{align}
\Var \bigg( \sum_{k=J}^{2N-1} a_k U_k \bigg) &=
\sum_{k=J}^{2N-1} \sum_{l=J}^{2N-1} a_k a_l \Cov(U_k, U_l) \leq
\frac{2 \sqrt{2}}{s^2} \sum_{k=J}^{2N-1} \sum_{l=k}^{2N-1}
\frac{1}{kl} (1-s)^{(l-k)/2} \nonumber \\
&\leq \frac{2 \sqrt{2}}{s^2} \sum_{k=J}^{2N-1} \frac{1}{k^2}
\bigg( \sum_{l=k}^{2N-1} (1-s)^{(l-k)/2} \bigg) \leq
C \sum_{k=J}^{2N-1} \frac{1}{k^2} \leq \frac{C}{J}.
\label{vwj5}
\end{align}

It remains to bound $\Var(\sum_{k=J}^{2N-1} b_k H_k)$.
Recall from the proof of Lemma \ref{udhlem} that
$$
p_k = P(X_t \neq X_{t-1}|X_{t-1} = k) = \frac{k(2N-k)(2-s)}{(2N)^2}
$$
and that $D_k + U_k = U_{k-1} - 1 + U_k$, using the convention that $U_0 = 1$.
Therefore, we can write
$H_k = G_1 + G_2 + \ldots + G_{U_k + U_{k-1} - 1}$, where
$(G_i)_{i=1}^{\infty}$ is an i.i.d. sequence of random variables such
that $G_i + 1$ has a Geometric($p_k$) distribution for all $i$.
Thus, $E[G_i] = p_k^{-1} - 1$.  If $k \leq l$, then by Lemma \ref{varcov}, 
\begin{align}
\Cov(H_k, H_l) &= \bigg( \frac{1}{p_k} - 1 \bigg) \bigg( \frac{1}{p_l} - 1
\bigg) \Cov( U_k + U_{k-1} - 1, U_l + U_{l-1} - 1) \nonumber \\
&\leq \frac{1}{p_k p_l} \Cov(U_k + U_{k-1}, U_l + U_{l-1}) \nonumber \\
&\leq \frac{4 \sqrt{2}}{s^2 p_kp_l} (1-s)^{(l-k-1)/2} \leq
\frac{C}{p_k p_l} (1-s)^{(l-k)/2}. \nonumber
\end{align}
Note that (\ref{bbk}) implies
$$
\frac{b_k}{p_k} \leq \frac{k(2N-k)}{2N^3} \cdot
\frac{(2N)^2}{k(2N-k)(2-s)} = \frac{2}{(2-s)N}
\leq \frac{2}{N}.
$$  
Therefore, 
\begin{align}
\Var \bigg( \sum_{k=J}^{2N-1} b_k H_k \bigg) &=
\sum_{k=J}^{2N-1} \sum_{l=J}^{2N-1} b_k b_l \Cov(H_k, H_l) \leq
C \sum_{k=J}^{2N-1} \sum_{l=k}^{2N-1} \frac{b_kb_l}{p_kp_l} (1-s)^{(l-k)/2}
\nonumber \\
&\leq \frac{C}{N^2} \sum_{k=J}^{2N-1} \sum_{l=k}^{2N-1} (1-s)^{(l-k)/2}
\leq \frac{C}{N^2} \sum_{k=J}^{2N-1} \frac{1}{1 - \sqrt{1-s}} \leq
\frac{C}{N}.
\label{vwj6}
\end{align}
The lemma follows from (\ref{vwj0}), (\ref{vwj5}), and (\ref{vwj6}). \qed

\begin{proof}[Proof of Proposition \ref{1linprop}] 
Lemma \ref{thetalemnew} gives
$$
\bigg| P(R(i) \geq \tau_J) - (1 - E[e^{-\eta_J}]) \bigg|  \leq
E \bigg[ \sum_{t=\tau_J+1}^{\tau} \theta_t^2 \bigg] \leq
\frac{C}{(\log N)^2}.
$$
Since $|\frac{d}{dx} e^{-x}| \leq 1$ for $x \geq 0$, Lemma \ref{etaJ} gives
$$
E[e^{-\eta_J'} - e^{-\eta_J}] \leq E[\eta_J - \eta_J'] \leq
\frac{C}{J(\log N)}.
$$
Using Jensen's inequality and Lemma \ref{varwj},
$$
|E[e^{-\eta_J'}] - e^{-E[\eta_J']}|  \leq E|e^{-\eta_J'} - e^{-E[\eta_J']}| \nonumber \\
\leq E|\eta_J' - E[\eta_J']| \leq \Var(\eta_J')^{1/2} \leq
\frac{C}{\sqrt{J} (\log N)}.
$$
Furthermore, it follows from Lemma \ref{ewjlem} that
$$
1 - e^{-E[\eta_J']} = q_J + O \bigg( \frac{1}{N} + \frac{(1-s)^J}{J \log N}
\bigg).
$$
Combining the last four equations gives the proposition.
\end{proof}

\section{Coalescence of two lineages}

In this section, we prove Propositions \ref{bcoalprop}, \ref{crprop},
and \ref{coaltJprop}, all of which pertain to the probabilities that two
lineages in the sample coalesce.  We begin by computing the following 
coalescence probabilities for integers $k$ and $l$ such that $1 \leq k \leq 2N-1$,
$1 \leq l \leq 2N$, and $|k - l| \leq 1$:
\begin{align}
p_{BB}^c(k,l) &= P(A_t^{t-1}(i) = A_t^{t-1}(j)|
X_{t-1} = k, X_t = l, B_t(i) = 1, B_t(j) = 1), \nonumber \\
p_{bb}^c(k,l) &= P(A_t^{t-1}(i) = A_t^{t-1}(j)|
X_{t-1} = k, X_t = l, B_t(i) = 0, B_t(j) = 0), \nonumber \\
p_{Bb}^c(k,l) &= P(A_t^{t-1}(i) = A_t^{t-1}(j)|
X_{t-1} = k, X_t = l, B_t(i) = 1, B_t(j) = 0). \nonumber
\end{align}
As with the recombination probabilities in the previous section, the Markov property implies that
the coalescence probabilities are the same under $P'$ as under $P$. 

\begin{Lemma}
We have
\begin{align}
p_{BB}^c(k, k &- 1) = p_{bb}^c(k, k+1) = 0, \nonumber \\
p_{BB}^c(k, k+1) & = \frac{2}{k(k+1)} \bigg( 1 - \frac{r(2N-k)}{2N} \bigg), \nonumber \\
p_{bb}^c(k, k-1) = & \frac{2}{(2N-k)(2N-k+1)}
\bigg( 1 - \frac{rk}{2N} \bigg), \nonumber \\
p_{bb}^c(k, k)  = \frac{2\beta_k}{k(2N-k)}
\bigg( 1 - \frac{rk}{2N} \bigg), & \quad
p_{BB}^c(k,k)  = \frac{2\beta_k}{k(2N-k)}
\bigg( 1 - \frac{r(2N-k)}{2N} \bigg),\nonumber \\
p_{Bb}^c(k,k) = \frac{r\beta_k}{k(2N-k)}, \quad
p_{Bb}^c(k,k+1) & = \frac{r}{2N(k+1)}, \quad
p_{Bb}^c(k,k-1) = \frac{r}{2N(2N - k + 1)}. \nonumber
\end{align}
\label{osclem}
\end{Lemma}

\begin{proof}
This result follows from a series of straightforward 
calculations, similar to those used to prove Lemma \ref{osrlem}.
We explain the idea behind some of these calculations.
When $X_{t-1} = k$ and $X_t = k-1$,
the new individual born at time $t$ has the $b$ allele.  Therefore, two
$B$ lineages can not coalesce at this time, so $p_{BB}^c(k, k-1) = 0$.
By the same reasoning, $p_{bb}^c(k, k+1) = 0$.
When $X_{t-1} = k$ and $X_t = k+1$, the new individual born at time $t$
has the $B$ allele.  With probability $r(2N-k)/2N$, this individual
inherits its allele at the neutral site from a member of the $b$
population because of recombination.  If this does not happen, then
two of the $B$ individuals get their allele at the neutral site
from the same parent.  Thus, conditional on $B_t(i) = B_t(j) = 1$, the
probability that the $i$th and $j$th individuals get their allele
at the neutral site from the same parent is $2/[k(k+1)]$, which
implies the formula for $p_{BB}^c(k, k+1)$.  The calculation of
$p_{bb}^c(k, k-1)$ is similar.

Now suppose $X_{t-1} = X_t = k$.  Conditional on this event,
a $B$ replaces a $B$ with probability $k^2/[k^2 + (2N-k)^2 + sk(2N-k)]$.
If the new $B$ gets its allele at the neutral site from a member
of the $B$ population, which has probability $1 - r(2N-k)/2N$, 
and if $B_t(i) = B_t(j) = 1$, then
the probability that the $i$th and $j$th lineages coalesce
is $2/k^2$, as there are $k$ possibilities both for the individual
who dies and the parent of the new individual.  The formula
for $p_{BB}^c(k,k)$ follows, and $p_{bb}^c(k,k)$ can be calculated
in the same way.  Next, to find $p_{Bb}^c(k,k)$,
note that if a $B$ replaces a $B$, and $B_t(i) = 1$ and
$B_t(j) = 0$, then the probability of coalescence is $r/(2kN)$,
as there must be recombination, and there are $k$ choices for
the $B$ individual that is just born and $2N$ choices for the
parent from which it gets its allele at the neutral site.
If instead a $b$ replaces a $b$, which happens with probability
$(2N-k)^2/[k^2 + (2N-k)^2 + sk(2N-k)]$ conditional on
$X_{t-1} = X_t = k$, the probability of coalescence is $r/[(2N-k)(2N)]$.
Adding the probabilities for the two cases
gives the formula for $p_{Bb}^c(k, k)$.

Finally, to calculate $p_{Bb}^c(k, k+1)$ and $p_{Bb}^c(k, k-1)$,
note that when a $B$ replaces a $b$, the probability that a $B$
lineage coalesces with a $b$ lineage is $r/[(k+1)(2N)]$, as there
must be recombination, and there are $k+1$ choices for the $B$
individual that was just born and $2N$ choices for its parent.
Likewise, the coalescence probability is $r/[(2N-k+1)(2N)]$
when a $b$ replaces a $B$.
\end{proof}

\begin{proof}[Proof of Proposition \ref{bcoalprop}]
We consider first the case in which the $j$th lineage is descended from a
member of the $B$ population at the time of coalescence.  Summing over the
possible values $k$ for $X_{G(i,j)}$ and applying Lemmas \ref{udhlem} and \ref{osclem}, we get
\begin{align}
P(G(i,j) &\geq 0, B_{G(i,j) + 1}(A_{\tau}^{G(i,j) + 1}(i)) = 0,
\mbox{ and }B_{G(i,j) + 1}(A_{\tau}^{G(i,j) + 1}(j)) = 1) \nonumber \\
&\leq \sum_{k=1}^{2N-1} \big( p_{Bb}^c(k, k+1) E[U_k] +
p_{Bb}^c(k, k-1) E[D_k] + p_{Bb}^c(k,k) E[H_k] \big) \nonumber \\
&\leq \sum_{k=1}^{2N-1} \bigg( \frac{r}{2N(k+1)s} + \frac{r}{2N(2N - k + 1)s}
+ \frac{r}{sk(2N-k)} \bigg) \nonumber \\
&\leq \frac{r}{2Ns} \sum_{k=1}^{2N-1}  \bigg( \frac{1}{k} +
\frac{1}{2N-k} + \frac{2N}{k(2N-k)} \bigg) \nonumber \\
& = \frac{2r}{s} \sum_{k=1}^{2N-1} \frac{1}{k(2N-k)} 
\leq \frac{4r}{Ns} \sum_{k=1}^N \frac{1}{k} 
\leq \frac{4r(1 + \log N)}{Ns} \leq \frac{C}{N}. \nonumber
\end{align}

It remains to consider the case in which the $i$th and $j$th lineages are
both descended from a member of the $b$ population at the coalescence time.
By summing over the possible values of $X_{R(i)}$ and $X_{G(i,j)}$,
we see that it suffices to show
\begin{align}
\sum_{\ell=1}^{2N-1} \sum_{k=1}^{2N-1} P(X_{R(i)} = \ell)
P&\bigg(X_{G(i,j)} = k, B_{G(i,j) + 1}(A_{\tau}^{G(i,j) + 1}(i)) = 0,
\mbox{ and } \nonumber \\
&B_{G(i,j) + 1}(A_{\tau}^{G(i,j) + 1}(j)) = 0 \bigg| X_{R(i)} = \ell \bigg)
\leq \frac{C (\log N)}{N}.
\label{ntscoal}
\end{align}
If $B_{G(i,j) + 1}(A_{\tau}^{G(i,j) + 1}(i)) = 0$, then
$G(i,j) \leq R(i)$.  Therefore, it follows from 
Lemmas \ref{udhlem} and \ref{osclem} and the time-reversal argument in the
proof of Proposition \ref{2rprop} that
\begin{align}
P(X_{G(i,j)} &= k \mbox{ and }B_{G(i,j) + 1}(A_{\tau}^{G(i,j) + 1}(i)) =
B_{G(i,j) + 1}(A_{\tau}^{G(i,j) + 1}(j)) = 0|X_{R(i)} = \ell) \nonumber \\
&\leq p_{bb}^c(k, k-1) E[D_{2N-k+1, 2N-\ell}] + 
p_{bb}^c(k,k) E[H_{2N-k, 2N-\ell}] \nonumber \\
&\leq \frac{2}{(2N-k)(2N-k+1)s} \min\{(1-s)^{k-\ell}, 1\} +
\frac{2}{sk(2N-k)} \min\{(1-s)^{k-\ell}, 1\} \nonumber \\
&\leq \bigg( \frac{2k + 2(2N - k)}{sk(2N-k)^2} \bigg) \min\{(1-s)^{k-\ell}, 1\}
= \frac{4N \min\{(1-s)^{k-l}, 1\}}{sk(2N-k)^2}. \nonumber
\end{align}
Combining this result with Lemma \ref{rkslem}, we get that the left-hand side
of (\ref{ntscoal}) is at most
\begin{align}
&\sum_{\ell=1}^{2N-1} \frac{r}{\ell s} \bigg( \sum_{k=1}^{2N-1}
\frac{4N \min\{(1-s)^{k-\ell}, 1\}}{sk(2N-k)^2} \bigg) \nonumber \\
&\qquad \leq \frac{4r}{s^2} \sum_{\ell=1}^{2N-1} \frac{1}{\ell}
\bigg( \sum_{k=\ell}^{2N-1} \frac{N(1-s)^{k-\ell}}{k(2N-k)^2} +
\sum_{k=1}^{\ell-1} \frac{N}{k(2N-k)^2} \bigg).
\label{2terms}
\end{align}
Using (\ref{logNNb}) and the fact that $N/[k(2N-k)] \leq 1$ for $1 \leq k \leq
2N-1$, we get
\begin{equation}
\frac{4r}{s^2} \sum_{\ell=1}^{2N-1} \frac{1}{\ell} \bigg( \sum_{k=\ell}^{2N-1}
\frac{N(1-s)^{k-\ell}}{k(2N-k)^2} \bigg) \leq \frac{4r}{s^2}
\bigg( \frac{2C(1 + \log N)}{N} \bigg) \leq \frac{C}{N}.
\label{ftm}
\end{equation}
For the second term in (\ref{2terms}), we have
\begin{align}
\frac{4r}{s^2} \sum_{\ell=1}^{2N-1} \frac{1}{\ell}
\bigg( \sum_{k=1}^{\ell-1} \frac{N}{k(2N-k)^2} \bigg)
&\leq \frac{4r}{s^2} \sum_{\ell=1}^N \frac{1}{\ell}
\bigg( \sum_{k=1}^{\ell} \frac{N}{kN^2} \bigg)
+ \frac{4r}{s^2} \sum_{l=N+1}^{2N-1} \frac{1}{N}
\bigg( \sum_{k=1}^{\ell} \frac{N}{k(2N-k)^2} \bigg) \nonumber \\
&\leq \frac{4r}{Ns^2} \bigg( \sum_{\ell=1}^N \frac{1}{\ell} \bigg)^2 +
\frac{4r}{s^2} \sum_{k=1}^{2N-1} \sum_{\ell=k}^{2N-1} \frac{1}{k(2N-k)^2}
\nonumber \\
&\leq \frac{4r(1 + \log N)^2}{Ns^2} + 
\frac{4r}{s^2} \cdot 2 \sum_{k=1}^{N} \frac{1}{k(2N-k)} \leq
\frac{C (\log N)}{N}.
\label{stm}
\end{align}
Using (\ref{ftm}) and (\ref{stm}) in (\ref{2terms}) proves (\ref{ntscoal}).
\end{proof}

The next lemma, which bounds the probability that there are $k$ individuals
with the $B$ allele at the time the $i$th and $j$th lineages coalesce,
will be used in the proofs of Propositions \ref{crprop} and \ref{coaltJprop}.

\begin{Lemma}
We have
\begin{equation}
P(X_{G(i,j)} = k \mbox{ and }
B_{G(i,j) + 1}(A_{\tau}^{G(i,j) + 1}(i)) =
B_{G(i,j) + 1}(A_{\tau}^{G(i,j) + 1}(j)) = 1) \leq \frac{4N}{sk^2(2N-k)}.
\label{coalkeq}
\end{equation}
\label{coalklem}
\end{Lemma}

\upat{\it Proof.}
By Lemmas \ref{udhlem} and \ref{osclem}, the probability on the
left-hand side of (\ref{coalkeq}) is at most
\begin{align}
\hphantom{xxxxxxx}
p_{BB}^c(k, k+1) E[U_k] + p_{BB}^c(k,k) E[H_k]
&\leq \frac{2}{sk(k+1)} + \frac{2}{sk(2N-k)} \nonumber \\
&\leq \frac{2(2N-k) + 2k}{sk^2(2N-k)} = \frac{4N}{sk^2(2N-k)}. 
\hphantom{xxxxxxx}\Box
\nonumber
\end{align}

\begin{proof}[Proof of Proposition \ref{crprop}]
By Proposition \ref{bcoalprop}, it suffices to show that
$$P(0 \leq R(i) \leq G(i,j) \mbox{ and }
B_{G(i,j) + 1}(A_{\tau}^{G(i,j) + 1}(i)) =
B_{G(i,j) + 1}(A_{\tau}^{G(i,j) + 1}(j)) = 1) \leq \frac{C}{\log N}.$$
By Lemmas \ref{udhlem} and \ref{osrlem} and the time-reversal argument in
the proof of Proposition \ref{2rprop},
\begin{align}
P(&X_{R(i)} = \ell \mbox{ and } 0 \leq R(i) \leq G(i,j)|X_{G(i,j)} = k)
\nonumber \\
&\leq p_B^r(\ell, \ell+1) E[U_{2N-\ell-1, 2N-k}] + p_B^r(\ell,\ell) E[H_{2N-\ell, 2N-k}]
\nonumber \\
&\leq \frac{r(2N-\ell)}{(\ell+1)(2Ns)} \min\{(1-s)^{\ell+1-k}, 1\} 
+ \frac{r}{2Ns} \min\{(1-s)^{\ell-k}, 1\}
\leq \frac{r}{\ell s} \min\{(1-s)^{\ell-k}, 1\}. \nonumber
\end{align}
Combining this result with (\ref{coalkeq}), we get
\begin{align}
P(0 \leq R(i) &\leq G(i,j) \mbox{ and }
B_{G(i,j) + 1}(A_{\tau}^{G(i,j) + 1}(i)) =
B_{G(i,j) + 1}(A_{\tau}^{G(i,j) + 1}(j)) = 1) \nonumber \\
&\leq \sum_{k=1}^{2N-1} \frac{4N}{sk^2(2N-k)} \bigg(
\sum_{\ell=1}^{2N-1} \frac{r}{\ell s} \min\{(1-s)^{\ell-k}, 1\} \bigg) \nonumber \\
&\leq \frac{4r}{s^2} \sum_{k=1}^{2N-1} \frac{N}{k^2(2N-k)} 
\bigg( \sum_{\ell=k}^{2N-1} \frac{(1-s)^{\ell-k}}{\ell} + \sum_{\ell=1}^{k-1} \frac{1}{\ell} \bigg).
\label{logNbd} 
\end{align}
The first term in the sum on the right-hand side of (\ref{logNbd}) is at most
$$
\sum_{k=1}^{2N-1} \frac{N}{k^3(2N-k)} \bigg( \sum_{\ell=k}^{2N-1}
(1-s)^{\ell-k} \bigg)
\leq \bigg( \frac{1}{s} \bigg) \bigg( \sum_{k=1}^N \frac{1}{k^3}
+ \sum_{k=N+1}^{2N-1} \frac{1}{N^2(2N-k)} \bigg),
$$
which is bounded by a constant.  The other term in the sum in (\ref{logNbd})
is at most
$$
\sum_{k=1}^{2N-1} \frac{N(1 + \log k)}{k^2(2N-k)} 
\leq \sum_{k=1}^N \frac{1 + \log k}{k^2} + \sum_{k=N+1}^{2N-1}
\frac{1 + \log(2N)}{N(2N - k)},
$$
which is also bounded by a constant.  Since $4r/s^2 \leq C/(\log N)$,
the proposition follows. \end{proof}

\begin{proof}[Proof of Proposition \ref{coaltJprop}]
By reasoning similar to that used to prove Lemma \ref{coalklem}, we have
\begin{align}
P(G(i,j) &\geq \tau_J \mbox{ and } B_{G(i,j) + 1}(A_{\tau}^{G(i,j) + 1}(i)) =
B_{G(i,j) + 1}(A_{\tau}^{G(i,j) + 1}(i)) = 1) \nonumber \\
&\leq \sum_{k=1}^{2N-1} \big( p_{BB}^c(k, k+1) E[U_{k,J}]
+ p_{BB}^c(k,k) E[H_{k,J}] \big).
\label{48eq}
\end{align}
However this time we keep the factor $\min\{(1-s)^{J-k}, 1\}$
from Lemma \ref{udhlem} to bound the right-hand side of (\ref{48eq}) by
\begin{equation}
\sum_{k=1}^J (1-s)^{J-k} \frac{4N}{sk^2(2N-k)} +
\sum_{k=J+1}^{2N-1} \frac{4N}{sk^2(2N-k)}.
\label{BJcoal}
\end{equation}
Using the fact that $N/[k(2N-k)] \leq 1$ for $1 \leq k \leq 2N-1$
and then Lemma \ref{aNlem}, we have
$$
\sum_{k=1}^J (1-s)^{J-k} \frac{4N}{sk^2(2N-k)} \leq
\frac{4}{s} (1-s)^J \sum_{k=1}^J \bigg(\frac{1}{1-s} \bigg)^k \frac{1}{k}
\leq \frac{C}{J}.
$$
For the second term in (\ref{BJcoal}), we observe
$$
\sum_{k=J+1}^{2N-1} \frac{4N}{sk^2(2N-k)} \le
\sum_{k=J+1}^{N-1} \frac{4}{sk^2} +  \sum_{k=N}^{2N-1} \frac{4}{sN(2N-k)}
\leq \frac{4}{sJ} + \frac{4(1 + \log N)}{Ns}.
$$
Since $J \leq C'N/(\log N)$, the bounds in the last two equations add up to
$C/J$, and the desired result follows from
these bounds and Proposition \ref{bcoalprop}.
\end{proof}

\section{Approximate independence of ${\bf n}$ lineages}

In this section, we prove Proposition \ref{indprop}.  We first establish
a lemma that involves the coupling of two $\{0, 1, \dots, n\}$-valued
random variables.

\begin{Lemma}
Let $V$ and $V'$ be $\{0, 1, \dots, n\}$-valued random variables such
that $E[V] = E[V']$.  Then, there exist random variables ${\tilde V}$
and ${\tilde V}'$ on some probability space such that
$V$ and ${\tilde V}$ have the same distribution, $V'$ and
${\tilde V}'$ have the same distribution, and
$$P({\tilde V} \neq {\tilde V}') \leq n \max\{ P({\tilde V} \geq 2),
P({\tilde V}' \geq 2)\}.$$
\label{coup}
\end{Lemma}

\upat
{\it Proof.} It is clear that ${\tilde V}$ and ${\tilde V}'$ can be constructed such that
they have the same distributions as $V$ and $V'$ respectively and
$P({\tilde V} = {\tilde V}') \geq \min\{P(V = 0), P(V' = 0)\} +
\min\{P(V = 1), P(V' = 1)\}.$
Note that $P(V = 0) \geq 1 - E[V]$.  Since $E[V] = E[V']$, it follows that
$\min\{P(V = 0), P(V' = 0)\} \geq 1 - E[V]$.  Also,
$P(V = 1) = E[V] - \sum_{k=2}^n k P(V = k),$
so $P(V = 1) \geq E[V] - n P(V \geq 2)$.  Likewise,
$P(V' = 1) \geq E[V] - n P(V' \geq 2)$.  It follows that
$$
\hphantom{xxxxxxxxxxxxxxxxxx}
P({\tilde V} = {\tilde V}') \geq 1 - n \max\{ P({\tilde V} \geq 2),
P({\tilde V}' \geq 2)\}.
\hphantom{xxxxxxxxxxxxxxxxxx}\Box
$$

Recall that $K_t = \# \{i \in \{1, \dots, n\}: R(i) \geq t\}$ 
for $0 \leq t \leq \tau$.  Define $\theta_t = p_B^r(X_{t-1}, X_t)$ as in
section 3, and define $\eta_J = \sum_{t=\tau_J+1}^{\tau} \theta_t$
and $\eta_J' = \sum_{t=\tau_J+1}^{\tau} \theta_t 1_{\{X_{t-1} \geq J\}}$
as in Lemma \ref{etaJ}.  Finally, let $F_J = P(R(i) \geq \tau_J|X)$,
which is shown in (\ref{defHt}) to be equal to
$1 - \prod_{t = \tau_J + 1}^{\tau} (1 - \theta_t)$.

\begin{Lemma}
If $J \le C'N/(\log N)$ then for all $d \in \{0, 1, \dots, n\}$,
$$
\bigg| P(K_{\tau_J} = d) - \binom{n}{d} E[F_J^d (1 - F_J)^{n-d}] \bigg|
\leq \min \bigg\{ \frac{C}{\log N}, \frac{C}{J} \bigg\} +
\frac{C}{(\log N)^2}.
$$
\label{indlem}
\end{Lemma}

\begin{proof}
Note that $K_{\tau} = 0$.  Also, $K_{t-1} - K_t \in
\{0, 1, \dots, n\}$ for all $1 \leq t \leq \tau$, and 
$$E[K_{t-1} - K_t|X, (K_u)_{u=t}^{\tau}] = (n - K_t)\theta_t.$$
Define another process $(K_t')_{t=0}^{\tau}$ such that
$K_{\tau}' = 0$ and the conditional distribution of
$K_{t-1}' - K_t'$ given $X$ and $(K_u')_{u=t}^{\tau}$ is
binomial($n - K_t', \theta_t$).  Note that
$E[K'_{t-1} - K'_t|X, (K'_u)_{u=t}^{\tau}] = (n - K'_t)\theta_t$.
We will show that the processes $(K_t)_{t=0}^{\tau}$ and
$(K_t')_{t=0}^{\tau}$ can be coupled so that
\begin{equation}
P(K_t \neq K'_t \mbox{ for some }t \geq \tau_J) \leq
\min \bigg\{ \frac{C}{\log N}, \frac{C}{J} \bigg\} + \frac{C}{(\log N)^2}.
\label{Ktcoup}
\end{equation}
Equation (\ref{Ktcoup}) implies the lemma because the conditional
distribution of $K_{\tau_J}'$ given $X$ is binomial with parameters $n$
and $1 - \prod_{t = \tau_J + 1}^{\tau} (1 - \theta_t) = F_J$.

By applying Lemma \ref{coup} with $V = K_{t-1} - K_t$ and
$V' = K_{t-1}' - K_t'$, we can construct the process
$(K_t')_{t=0}^{\tau}$ on the same probability space as
$(K_t)_{t=0}^{\tau}$ such that
\begin{align}
P(K_t \neq K'_t \mbox{ for some }t \geq \tau_J|X) &\leq n \sum_{t=\tau_J+1}^{\tau}  P(K_{t-1} - K_t \geq 2|X, (K_u)_{u=t}^{\tau}) \nonumber \\
&\hspace{.2in} + n \sum_{t=\tau_J+1}^{\tau}
P(K_{t-1}' - K_t' \geq 2|X, (K'_u)_{u=t}^{\tau}).
\label{indeq1}
\end{align}
If $K_{t-1} - K_t \geq 2$ for some $t \geq \tau_J$, then 
$\tau_J \leq R(i) \leq G(i,j)$ for some $i$ and $j$.  We have
$P(\tau_J \leq R(i) \leq G(i,j)) \leq C/(\log N)$ for all $J$ by
Proposition \ref{crprop} and $P(\tau_J \leq R(i) \leq G(i,j)) \leq C/J$
for all $J \leq C'N/(\log N)$ by Proposition 2.5.
Therefore, for $J \leq C'N/(\log N)$,
\begin{align}
E \bigg[\sum_{t=\tau_J+1}^{\tau}  P(K_{t-1} - K_t \geq 2|X, (K_u)_{u=t}^{\tau})
\bigg] &\leq \sum_{t=1}^{\tau} P(K_{t-1} - K_t \geq 2 \mbox{ and }
t \geq \tau_J) \nonumber \\
&\leq \frac{n}{2} P(K_{t-1} - K_t \geq 2 \mbox{ for some }t \geq \tau_J)
\nonumber \\
&\leq \min \bigg\{ \frac{C}{\log N}, \frac{C}{J} \bigg\}.
\label{indeq2}
\end{align}
Now a binomial random variable will be at least 2
if and only if there is some pair of successful trials, so 
$P(K_{t-1}' - K_t' \geq 2|X, (K'_u)_{u=t}^{\tau}) \leq
\binom{n}{2} \theta_t^2$, and
\begin{equation}
\sum_{t=\tau_J+1}^{\tau}  P(K'_{t-1} - K'_t \geq 2|X, (K'_u)_{u=t}^{\tau})
\leq \binom{n}{2} \sum_{t=\tau_J+1}^{\tau} \theta_t^2.
\label{indeq3}
\end{equation}
By taking expectations in (\ref{indeq1}) and applying (\ref{indeq2}),
(\ref{indeq3}), and (\ref{thetatexp}), we get (\ref{Ktcoup}),
which completes the proof.
\end{proof}

\begin{proof}[Proof of Proposition \ref{indprop}]
In view of Lemma \ref{indlem}, it suffices to show that
\begin{equation}
\big|E[F_J^d(1 - F_J)^{n-d}] - q_J^d (1 - q_J)^{n-d} \big| \leq
\min \bigg\{ \frac{C}{\log N}, \frac{C}{J} \bigg\} +
\frac{C}{(\log N)^2}
\label{Hqeq}
\end{equation}
for all $d \in \{0, 1, \dots, n\}$.    If
$0 \leq a_1, \dots, a_n \leq 1$ and $0 \leq b_1, \dots, b_n \leq 1$,
then $|a_1 \dots a_n - b_1 \dots b_n| \leq \sum_{i=1}^n |a_i - b_i|$,
as shown in Lemma 4.3 of chapter 2 of Durrett (1996).  Therefore,
$$\big| E[F_J^d (1 - F_J)^{n-d}] - q_J^d (1 - q_J)^{n-d} \big|
\leq E[ d|F_J - q_J| + (n-d) |(1 - F_J) - (1 - q_J)| ] = n E[|F_J - q_J|].$$
Note that
\begin{equation}
|F_J - q_J| \leq |F_J - (1 - e^{-\eta_J})| + |e^{-\eta_J'} - e^{-\eta_J}|
+ |e^{-\eta_J'} - e^{-E[\eta_J']}| + |(1 - e^{-E[\eta_J']}) - q_J|.
\label{indineq}
\end{equation}
It follows from (\ref{poisapp}) and (\ref{thetatexp}) that
$E[|F_J - (1 - e^{-\eta_J})|] \leq C/(\log N)^2$.  The expectations of the
second, third, and fourth terms on the right-hand side of (\ref{indineq})
can be bounded as in the conclusion of the proof of Proposition \ref{1linprop} 
at the end of Section 3. All of those error estimates are smaller than
the right-hand side of (\ref{Hqeq}) so the desired result follows.
\end{proof}

\section{A branching process approximation}

In this section, we will show how the evolution of the
individuals with the $B$ allele during the first stage of the selective
sweep can be approximated by a supercritical branching process.
This will lead to a
proof of Proposition \ref{branchprop}.  Recall that the first stage of
the sweep consists of the times $0 \leq t \leq \tau_J$, where
$J = \lfloor (\log N)^a \rfloor$ for some fixed constant $a > 4$.
We will assume throughout this section that $N$ is large enough that
$J \leq N$.  In subsection 6.1, we explain the coupling between the
branching process and the population model.  In subsection 6.2, 
we consider the lineages in the branching process with an infinite
line of descent.  Proposition \ref{branchprop} is proved using
these ideas in subsection 6.3.

\subsection{Coupling the population model with a branching process}

We begin by constructing a multi-type branching process with the
properties mentioned in Proposition \ref{branchprop}.  That is,
the process will start with one individual at time zero, and each
individual will give birth at rate one and die at rate $1-s$.
Each new individual has the same type as its parent with probability $1-r$
and a new type, different from all other types, with probability $r$.
We now explain how to construct this branching process so that
until the number of individuals reaches $J$, the branching process
will be coupled with the population process $(M_t)_{t=0}^{\infty}$
with high probability.

Define random variables $0 = \xi_0 < \xi_1 < \dots$ such that
$(\xi_i - \xi_{i-1})_{i=1}^{\infty}$ is an i.i.d.~sequence of random
variables, each having an exponential distribution with mean $1/2N$.
The branching process will start with one individual at time zero.
Until the population size reaches $J$,
there will be no births during the intervals $(\xi_{t-1}, \xi_t)$,
but births and deaths can occur at the times $\xi_1, \xi_2, \dots$.
This branching process will be coupled with $(M_t)_{t=0}^{\infty}$
so that, with high probability, the number of individuals with
the $B$ allele at time $t$ will be the same as the number of individuals
in the branching process at time $\xi_t$.  To facilitate this coupling, we will
also assign to each individual in the branching process a label
such that all the individuals alive at a given time have distinct labels.
We denote by $L_t$ the set of 
all $i$ such that there is an individual labeled $i$ in the population
at time $\xi_t$.  When $L_t = \{i: B_t(i) = 1\}$, meaning that the labels
are the same as the individuals in the population model with the
$B$ allele at time $t$, we say the coupling holds at time $t$.
The label of the individual at time zero will be $U$,
where $U$ is the random variable with a uniform
distribution on $\{1, \dots, 2N\}$ defined at the beginning of section 2.
We have $B_0(U) = 1$, so the coupling holds at time zero.

For the branching process to have the desired properties, each
individual must have probability $1/2N$ of giving birth at time
$\xi_t$ and probability $(1-s)/2N$ of dying at time $\xi_t$.
Also, at most one birth or death event can occur at a time.
Suppose the coupling holds at time $\xi_{t-1}$ and $i \in L_{t-1}$.
Also, assume $X_{t-1} = k$.
In the population model, the number of $B$'s increases by one at
time $t$, with $i$ being the parent of the new individual, if
$I_{t,2} = i$ and $B_{t-1}(I_{t,1}) = 0$, which has probability
$(2N - k)/(2N)^2$.  Also, the $i$th individual in the population
dies at time $t$, causing the $B$ population to decrease in size
by one, if $I_{t,1} = i$, $B_{t-1}(I_{t,2}) = 0$, and $I_{t,4} = 1$,
which has probability $(2N-k)(1-s)/(2N)^2$.  Consequently, we can
define the branching process such that the individual labeled
$i$ gives birth at time $\xi_t$ if and only if $I_{t,2} = i$, which
has probability $1/2N$.  We give the new individual the label
$I_{t,1}$, unless one of the other individuals already has this
label.  As a result, the coupling will hold at time $t$ if
$B_{t-1}(I_{t,1}) = 0$ but not if $B_{t-1}(I_{t,1}) = 1$.
The individual labeled $i$ will die
with probability $(1-s)/2N$, and will die whenever
$I_{t,1} = i$, $B_{t-1}(I_{t,2}) = 0$, and $I_{t,4} = 1$.  Then,
the probability that the coupling fails to hold at time $t$ is
\begin{equation}
k \bigg( \frac{1}{2N} - \frac{2N-k}{(2N)^2} \bigg) +
k \bigg( \frac{(1-s)}{2N} - \frac{(2N-k)(1-s)}{(2N)^2} \bigg) =
\frac{k^2(2-s)}{(2N)^2}.
\label{coupprob}
\end{equation}
If a new individual in the branching process is born at time $t$,
we say that it has a new type whenever $I_{t,5} = 1$, which has
probability $r$.  This means that births of individuals with new
types correspond to recombinations in the population model.

Fix a positive integer $m$.  On the event that the branching process has
at least $J$ individuals at some time, we define a random marked partition
${\tilde \Psi}_m$ as follows.  Define $\kappa$ such that $\xi_{\kappa}$
is the first time at which there are $J$ individuals.
Define a random injective map ${\tilde \sigma}:  \{1, \ldots, m\} \rightarrow
L_{\kappa}$ such that all $(J)_m$ possible maps are equally likely.
Then say that $i \sim_{{\tilde \Psi}_m} j$ if and only if the individuals
labeled ${\tilde \sigma}(i)$ and ${\tilde \sigma}(j)$ are of the same type.
Mark the block of ${\tilde \Psi}_m$ consisting of all $i$ such that
the individual labeled ${\tilde \sigma}(i)$ has the same type as the individual
at time zero.  Furthermore, we can define ${\tilde \sigma}$ such that $\sigma = {\tilde \sigma}$ on the
event that $\kappa = \tau_J$ and $L_{\tau_J} = \{i: B_{\tau_J}(i) = 1\}$, where
$\sigma: \{1, \ldots, m\} \rightarrow \{i: B_{\tau_J}(i) = 1\}$ is the
map defined in the section 2 that is used in the construction of
the random marked partition $\Psi_m$.  Recall that $i \sim_{\Psi_m} j$
if and only if $A_{\tau_J}^0(\sigma(i)) = A_{\tau_J}^0(\sigma(j))$,
and the block $\{i: B_0(A_{\tau_J}^0(\sigma(i))) = 1\}$ is marked.

Suppose $X_t = J$ for some $t$ and
the coupling holds for all $t \leq \tau_J$, so $\kappa = \tau_J$.
Then, the genealogy of the branching process is the same as the
genealogy of the $B$'s in the population up to time $\tau_J$.
Furthermore, groups of individuals in the branching process with
the same type correspond to groups of lineages in the population
that escape the selective sweep at the
same time, and therefore get their allele at the neutral site
from the same ancestor.  Therefore, we will have ${\tilde \Psi}_m = \Psi_m$
unless one of the following happens to a sampled lineage
during the first stage of the selective sweep:
\begin{enumerate}
\item One of the $B$ lineages experiences recombination, but the
allele at the neutral site comes from another $B$ individual.

\item Two recombinations cause a lineage to go from the $B$ population
to the $b$ population, and then back into the $B$ population.

\item There is a coalescence event involving at least one lineage
in the $b$ population.
\end{enumerate}
More formally, the lemma below is a consequence of our
construction.  Note that the events $\Lambda_3^c$, $\Lambda_4^c$,
and $\Lambda_5^c$ correspond to the three possibilities mentioned above.

\begin{Lemma}
Let $R_J(i) = \sup\{t \geq 0: B_t(A^t_{\tau_J}(i)) = 0\}$ and
$G_J(i,j) = \sup\{t \ge 0: A^t_{\tau_J}(i) = A^t_{\tau_J}(j)\}$. 
We have $\Psi_m = {\tilde \Psi}_m$ on the event $\Lambda_1 \cap
\dots \cap \Lambda_5$, where 
\begin{description}
\item []$\Lambda_1$ is the event that $X_t = J$ for some $t$,

\item []$\Lambda_2$ is the event that the coupling holds for all
$t \leq \tau_J$,

\item []$\Lambda_3$ is the event that for all $t \leq \tau_J$ for which
$B_{t-1}(I_{t,2}) = 1$, we have $B_{t-1}(I_{t,3}) = 0$,

\item []$\Lambda_4$ is the event that for
$i \in \{1, \dots, m\}$, we have $B_t(A_{\tau_J}^t(\sigma(i))) = 0$
for all $t \leq R_J(i)$, and

\item []$\Lambda_5$ is the event that for all
$i,j \in \{1, \dots, m\}$ with $G_J(\sigma(i), \sigma(j)) \geq 0$, we have
\vspace{-.1in}
$$B_{G_J(\sigma(i),\sigma(j))+1}(A_{\tau_J}^{G_J(\sigma(i), \sigma(j)) + 1}
(\sigma(i))) = B_{G_J(\sigma(i),\sigma(j))+1}(A_{\tau_J}^{G_J(\sigma(i),
\sigma(j)) + 1} (\sigma(j))) = 1.$$
\end{description}
\label{same}
\end{Lemma}

\begin{proof}
We have seen that when $\Lambda_1$ and $\Lambda_2$ occur, we have
$L_{\tau_J} = \{i: B_{\tau_J}(i) = 1\}$ and $\sigma = {\tilde \sigma}$.
For integers $u \leq t$ and $i \in L_t$, let ${\tilde A}_t^u(i)$ be the label
of the individual in the branching process at time $\xi_u$ that is the
ancestor of the individual labeled $i$ at time $\xi_t$, unless the
ancestor is of a different type then the individual labeled $i$
at time $t$, in which case we define ${\tilde A}_t^u(i) = 0$.
Note that when $\Lambda_1$ and $\Lambda_2$ occur,
we have $i \sim_{{\tilde \Psi}_m} j$ if and only if
${\tilde A}^t_{\tau_J}({\tilde \sigma}(i)) =
{\tilde A}^t_{\tau_J}({\tilde \sigma}(j)) \neq 0$ for some $t$.

Since $\sigma = {\tilde \sigma}$ when $\Lambda_1$ and $\Lambda_2$ occur,
we have $i \sim_{{\tilde \Psi}_m} j$ if and only if
${\tilde A}^t_{\tau_J}(\sigma(i)) = {\tilde A}^t_{\tau_J}(\sigma(j)) \neq 0$
for some $t$.
Suppose $j \in L_t$.  It follows from the constructions that
$A_t^{t-1}(j) = {\tilde A}_t^{t-1}(j)$ unless $j = I_{t,1}$
and $I_{t,5} = 1$.  In this case,
${\tilde A}_t^{t-1}(j) = 0$, and if $\Lambda_3$ occurs then
$B_{t-1}(A_t^{t-1}(j)) = 0$.  It follows that if $\Lambda_4$ also occurs, then 
${\tilde A}_{\tau_J}^t(\sigma(i)) = \tilde A_{\tau_J}^t(\sigma(j)) \neq 0$
if and only if we have both $A_{\tau_J}^t(\sigma(i)) = A_{\tau_J}^t(\sigma(j))$
and $B_t(A_{\tau_J}^t(\sigma(i))) = B_t(A_{\tau_J}^t(\sigma(j))) = 1$.
Furthermore, when $\Lambda_5$ occurs, we have both
$A_{\tau_J}^t(\sigma(i)) = A_{\tau_J}^t(\sigma(j))$
and $B_t(A_{\tau_J}^t(\sigma(i))) = B_t(A_{\tau_J}^t(\sigma(j))) = 1$ for some
$t$ if and only if $A_{\tau_J}^0(\sigma(i)) = A_{\tau_J}^0(\sigma(j))$,
which is exactly the condition for $i \sim_{\Psi_m} j$.  Thus, when
$\Lambda_1, \dots, \Lambda_5$ all occur, we have
$i \sim_{\Psi_m} j$ if and only if $i \sim_{{\tilde \Psi}_m} j$.

It remains only to show that the marked blocks of 
$\Psi_m$ and ${\tilde \Psi}_m$ are the same.  Note that $i$ is in the
marked block of ${\tilde \Psi}_m$ if and only if ${\tilde \sigma}(i)
= \sigma(i)$ has the same type as the individual at time zero or,
equivalently, if and only if ${\tilde A}_{\tau_J}^0(\sigma(i)) \neq 0$.
The fact that this condition is equivalent to
$B_0(A_{\tau_J}^0(\sigma(i))) = 1$ follows from the coupling
and conditions $\Lambda_3$ and $\Lambda_4$.
\end{proof}

We now use this coupling to show that the partition
${\tilde \Psi}_m$ conditioned on the survival of the branching
process has almost the same distribution as $\Psi_m$.

\begin{Lemma}
Let $\pi$ be a partition of $\{1, \dots, m\}$.  Then, there exists a
constant $C$ such that
$$|P'({\tilde \Psi}_m = \pi|\#L_t > 0 \mbox{ for all }t \in \N) -
P(\Psi_m = \pi)| \leq C/(\log N)^2.$$
\label{lem20}
\end{Lemma}

\vspace{-.2in}
\begin{proof}

We will show that if $\Lambda_1$ occurs, then
$\Lambda_2 \cap \dots \cap \Lambda_5$ occurs with high probability.
Conditional on the event that $X_{t-1} = k$ and that the coupling holds
at time $t-1$, it follows from (\ref{coupprob}) that
the probability that the coupling fails to hold at time $t$ is
$k^2(2-s)/(2N)^2$.  Likewise, conditional on these same events,
the probability that $B_{t-1}(I_{t,2}) = B_{t-1}(I_{t,3}) = 1$
is $(k/2N)^2$.  Thus, if $D_t$ is the event that $t$ is the
first integer such that either the coupling fails at time $t$ or
$B_{t-1}(I_{t,2}) = B_{t-1}(I_{t,3}) = 1$, then
$P'(D_t|X_t = k) \leq (3-s)k^2/(2N)^2$, where we use $P'$ because
we are not conditioning on the event that $X_t = 2N$ for some $t$.
Therefore,
\begin{align}
P'(\Lambda_1 \cap (\Lambda_2^c \cup \Lambda_3^c)) &\leq
\sum_{t=1}^{\infty} P'(D_t \cap \{t \leq \tau_J < \infty\})
= \sum_{t=1}^{\infty} E'[P'(D_t \cap \{t \leq \tau_J < \infty\}|X_{t-1})]
\nonumber \\
&\leq \sum_{t=1}^{\infty} E' \bigg[ \frac{(3-s)X_{t-1}^2}{(2N)^2}
1_{\{X_{t-1} \leq J\}} \bigg] = \frac{3-s}{(2N)^2}
\sum_{t=1}^{\infty} E'[X_{t-1}^2 1_{\{X_{t-1} \leq J\}}] \nonumber \\
&\leq \frac{3-s}{(2N)^2} \sum_{k=1}^J k^2 E'[T_k]. \nonumber
\end{align}
Since $P'(X_t \neq X_{t-1}|X_{t-1} = k) = P(X_t \neq X_{t-1}|X_{t-1} = k) =
p_k = k(2N-k)(2-s)/(2N)^2$ and $E'[U_k + D_k] \leq C$, it follows that
\begin{align}
P'(\Lambda_1 \cap (\Lambda_2^c \cup \Lambda_3^c)) &\leq 
\frac{3-s}{(2N)^2} \sum_{k=1}^J k^2 \frac{E'[U_k + D_k]}{p_k}
\nonumber \\
&\leq \frac{C}{N^2} \sum_{k=1}^J \frac{k^2 (2N)^2}{k(2N-k)}
\leq C \sum_{k=1}^J \frac{k}{2N-k} \leq \frac{CJ^2}{N}. \nonumber
\end{align}
To handle $\Lambda_4$ and $\Lambda_5$, note that
\begin{equation}
P'(X_{\tau} = 2N|\Lambda_1) = p(0, 2N, J) =
\frac{1 - (1-s)^J}{1 - (1-s)^{2N}} \geq 1 - (1-s)^J.
\label{condiff}
\end{equation}
It follows from (\ref{condiff}) and the proof of
Proposition \ref{2rprop} that $P'(\Lambda_1 \cap \Lambda_4^c)
\leq C/(\log N)^2$.  Likewise, it follows from (\ref{condiff}) and the proof
of Proposition \ref{bcoalprop} that $P'(\Lambda_1 \cap \Lambda_5^c)
\leq C(\log N)/N$.

Since $P'(\Lambda_1) = s/(1 - (1-s)^J)$ by Lemma \ref{survive}, it follows
from the above calculations that 
$|P'(\Lambda_1 \cap \dots \cap \Lambda_5) - s| \leq C/(\log N)^2$.  Recall
that $P'(X_{\tau} = 2N) = s/(1 - (1-s)^{2N})$ by Lemma \ref{survive}.  
Since $\{ \# L_t > 0 \mbox{ for all }t \in \N\}$ is the event that the
branching process survives, it is well-known that
$P'(\# L_t > 0 \mbox{ for all }t \in \N) = s$.  Furthermore, if
$\Lambda_1 \cap \dots \cap \Lambda_5$ occurs, then $X_t = J$ for some $t$
and $\# L_t = J$ for some $t$.  Note that $P'(X_{\tau} = 2N|X_t = J
\mbox{ for some } t) \geq 1 - (1-s)^J$ as in (\ref{condiff}) and
$P'(\# L_t > 0 \mbox{ for all }t|\# L_t = J \mbox{ for some }t) =
1 - (1-s)^J$.  Thus, the events $\Lambda_1 \cap \dots \cap \Lambda_5$,
$\{X_{\tau} = 2N\}$, and $\{\# L_t = 0 \mbox{ for all }t\}$ agree closely
enough that the probability, under $P'$, that either all or none of these
three events occurs is at least $1 - C/(\log N)^2$.  It follows from
this observation, Lemma \ref{same}, and the fact that $P$ is the conditional
probability measure of $P'$ given $X_{\tau} = 2N$ that
\begin{align}
P'({\tilde \Psi}_m = \pi| \#L_t > 0 \mbox{ for all }t \in \N) &=
P'({\tilde \Psi}_m = \pi| \Lambda_1 \cap \dots \cap \Lambda_5) +
O((\log N)^{-2}) \nonumber \\
&= P'(\Psi_m = \pi| \Lambda_1 \cap \dots \cap \Lambda_5) + O((\log N)^{-2})
\nonumber \\
&= P'(\Psi_m = \pi| X_{\tau} = 2N) + O((\log N)^{-2}) \nonumber \\
&= P(\Psi_m = \pi) + O((\log N)^{-2}), \nonumber 
\end{align}
which proves the lemma.
\end{proof}

\subsection{Infinite lines of descent}

Consider a continuous-time branching process in which each
individual gives birth at rate $1$ and dies at rate $1-s$.  Equivalently,
each individual lives for an exponentially distributed time with mean
$1/(2-s)$, and then has some number of offspring, which is $0$ with
probability $(1-s)/(2-s)$ and $2$ with probability $1/(2-s)$.
Say that an individual at time $t$ has an infinite line of descent if
it has a descendant in the population at time $u$ for all $u > t$.
Otherwise, say that the individual has a finite line of descent.

Define the process
$(Y^{(1)}_t, Y^{(2)}_t)_{t \geq 0}$ such that $Y^{(1)}_t$ is the
number of individuals at time $t$ having an infinite line of descent
and $Y^{(2)}_t$ is the number of individuals having a finite line of
descent.  Gadag and Rajarshi (1992) show that this process is
a two-type Markov branching process.  They also show that the behavior
of the process can be described as follows.  Let $p_k$ be the probability 
that an individual has $k$ offspring and let
$f(x) = \sum_{k=0}^{\infty} p_kx^k$ be the generating function of the
offspring distribution.  Let $u(x) = b[f(x) - x]$, where $b^{-1}$ is the mean
lifetime of an individual.  Let
$f^{(1)}(x,y) = \sum_{j=0}^{\infty} \sum_{k=0}^{\infty} p^{(1)}_{jk}x^jy^k,$
where $p^{(1)}_{jk}$ is the probability that an individual with an
infinite line of descent has $j$ offspring with an infinite line of
descent and $k$ offspring with a finite line of descent.  Let 
$f^{(2)}(x,y) = \sum_{j=0}^{\infty} \sum_{k=0}^{\infty} p^{(2)}_{jk}x^jy^k,$
where $p^{(2)}_{jk}$ is the probability that an individual with a
finite line of descent has $j$ offspring with an infinite line of
descent and $k$ offspring with a finite line of descent.  
Let $u^{(1)}(x,y) = b[f^{(1)}(x,y) - x]$, and let
$u^{(2)}(x,y) = b[f^{(2)}(x,y) - y]$.  Let $q$ be the smallest
nonnegative solution of the equation $u(x) = 0$, which is also the probability
that the branching process dies out.  Then, by
equation (4) of Gadag and Rajarshi (1992),
$$
u^{(1)}(x,y) = \frac{u(x(1-q) + yq) - u(yq)}{1-q}, \quad\hbox{and}\quad
u^{(2)}(x,y) = \frac{u(yq)}{q}.
$$

In the case of interest to us, we have $f(x) = \frac{1-s}{2-s} +
\frac{1}{2-s} x^2,$ and therefore 
$$u(x) = (2-s)[f(x) - x] = (1-s) + x^2 - (2-s)x.$$
Since $u(x) = x$ if and only if $x \in \{1-s, 1\}$, we have $q = 1-s$.
It follows that
\begin{align}
u^{(1)}(x,y) &= \frac{[xs + y(1-s)]^2 -
(2-s)[xs+ y(1-s)] - [y(1-s)]^2 + (2-s)[y(1-s)]}{s} \nonumber \\
&= sx^2 + 2(1-s)xy - (2-s)x. \nonumber
\end{align}
Thus, an individual with an infinite line of descent lives for an
exponentially distributed time with mean $1/(2-s)$.  It is replaced
by two individuals with infinite lines of descent at rate $s$, and it
is replaced by one individual with an infinite line of descent and another
individual with a finite line of descent at rate $2(1-s)$. 

Now, consider the process $(Y_t^{(1)}, Y_t^{(2)})$ started with one individual
and conditioned to survive forever, which is equivalent to assuming that
$Y^{(1)}_0 = 1$ and $Y^{(2)}_0 = 0$.
Assume, as in Proposition \ref{branchprop},
that the individuals are assigned types,
and that each new individual born is the same type as its parent with
probability $1-r$ and is a new type with probability $r$. Define
$\lambda^* = \inf\{t: Y_t^{(1)} = \lfloor Js \rfloor\}$. Let
$\lambda_k = \inf\{t: Y_t^{(1)} + Y_t^{(2)} = k\}$. 
Let $J_1 = \big\lfloor J \big( 1 + s^{-1} \sqrt{(\log J)/J} \big)^{-1}
\big\rfloor$
and $J_2 = \big\lceil J \big(1 - s^{-1} \sqrt{(\log J)/J} \big)^{-1}
\big\rceil$.  

\begin{Lemma}
We have $1 - P(\lambda_{J_1} \leq \lambda^* \leq \lambda_{J_2}) \leq
C/(\log N)^8$.
\label{J1J2}
\end{Lemma}

\begin{proof} If $S$ has a
binomial$(n,p)$ distribution and $p < c < 1$, then we have the large
deviations result that $P(S \geq cn) \leq e^{-2n(c-p)^2}$ (see Johnson,
Kotz, and Kemp (1992)).

Let $S_1$ have a binomial$(J_1, s)$ distribution, and let
$S_2$ have a binomial$(J_2, s)$ distribution.  Let $c = s +
\sqrt{(\log J)/J}$.  Then $J_1 = \lfloor Js/c \rfloor$, so
$c J_1 \leq Js$ and therefore
$$
P(\lambda^* \leq \lambda_{J_1}) = P(S_1 \geq \lfloor Js \rfloor|S_1 > 0) =
\frac{P(S_1 \geq \lfloor Js \rfloor)}{P(S_1 > 0)} \leq
\frac{P(S_1 \geq \lfloor cJ_1 \rfloor)}{1-(1-s)^{J_1}} \leq
\frac{P(S_1 \geq (c - \frac{1}{J_1}) J_1)}{1-(1-s)^{J_1}}.
$$
Recalling $J = \lfloor (\log N)^a \rfloor$ with $a>4$, it follows that if 
$\epsilon > 0$ is small, then for large $N$
$$
P(\lambda^* \leq \lambda_{J_1}) \leq 2
e^{-2 J_1(\sqrt{(\log J)/J} - J_1^{-1})^2} \leq Ce^{-2(J_1/J) \log J} 
\leq CJ^{-(2-\epsilon)} \leq C/(\log N)^8.
$$
Likewise, if $d = (1-s) + \sqrt{(\log J)/J}$, then
$J_2 = \lceil Js/(1-d) \rceil$, so $(1-d) J_2 \geq Js$ and thus
\begin{align}
P(\lambda^* > \lambda_{J_2}) &= P(S_2 < \lfloor Js \rfloor|S_2 > 0) \leq
P(S_2 < \lfloor Js \rfloor) \nonumber \\
&= P(J_2 - S_2 > J_2 - \lfloor Js \rfloor)
\leq P(J_2 - S_2 \geq dJ_2). \nonumber
\end{align}
Therefore, $P(\lambda^* > \lambda_{J_2}) \leq
e^{-2(J_2/J) \log J} \leq J^{-2} \leq C/(\log N)^8$,
and the lemma follows.
\end{proof}

\subsection{Proof of Proposition \ref{branchprop}}

We now prove Proposition \ref{branchprop}.
Recall that $\Upsilon_m$ is the marked partition obtained by sampling $m$
of the $\lfloor Js \rfloor$ individuals at time $\lambda^*$ having an infinite
line of descent and then declaring $i$ and $j$ to be in the same block of
$\Upsilon_m$ if and only if the $i$th and $j$th individuals in the sample
have the same type.  The marked block of $\Upsilon_m$ consists of the
individuals in the sample with the same type as the individual at time zero.
We now define three other random marked partitions $\Upsilon_m^{(1)}$,
$\Upsilon_m^{(2)}$, and $\Upsilon_m^{(3)}$ in the same way, except that the
sample of $m$ individuals is taken differently for each partition.  Namely,
to obtain $\Upsilon_m^{(1)}$, we sample $m$ of the individuals at time
$\lambda_J$.  To get $\Upsilon_m^{(2)}$, we sample $m$ of the individuals
at time $\lambda_{J_2}$.  To get $\Upsilon_m^{(3)}$, we sample $m$
of the individuals at time $\lambda_{J_2}$ that have an infinite line of
descent, assuming that $m$ such individuals exist (otherwise, sample from
all individuals at time $\lambda_{J_2}$).

Since the branching process has been conditioned to survive forever,
$\Upsilon_m^{(1)}$ has the same distribution as the conditional distribution of
${\tilde \Psi}_m$ given $\# L_t > 0$ for all $t \in \N$.  Thus,
by Lemma \ref{lem20}, it suffices to show that for all marked partitions
$\pi \in {\cal P}_m$, we have
$$
|P(\Upsilon_m^{(1)} = \pi) - P(\Upsilon_m = \pi)| \leq \frac{C}{(\log N)^2}.
$$
Note also that $\Upsilon_m^{(2)}$ and $\Upsilon_m^{(3)}$ have the same
distribution by the strong Markov property.

We can couple $\Upsilon_m^{(1)}$ and $\Upsilon_m^{(2)}$ such that the sample
at time $\lambda_J$ used to construct $\Upsilon_m^{(1)}$ includes all of the
the individuals in the sample at time $\lambda_{J_2}$ that were born before
time $\lambda_J$.  If there are fewer than $m$ such individuals, the rest of
the sample at time $\lambda_J$ can be picked from the remaining individuals.
By the strong Markov property, this way of picking the sample at time
$\lambda_J$ does not change the distribution of $\Upsilon_m^{(1)}$.
Therefore, $\Upsilon_m^{(1)} = \Upsilon_m^{(2)}$ if the $m$ individuals
sampled when constructing $\Upsilon_m^{(2)}$ were all born before time
$\lambda_J$.  Likewise, we can couple the partitions
$\Upsilon_m$ and $\Upsilon_m^{(3)}$ such that on the event
$\lambda^* \leq \lambda_{J_2}$, all of the individuals sampled at time
$\lambda_{J_2}$ that were born before time $\lambda^*$ are part of the sample
at time $\lambda^*$ used to construct $\Upsilon_m$.
Note that $\lambda^*$ is a stopping time with respect to the process
$(Y_t^{(1)}, Y_t^{(2)})_{t \geq 0}$, so the strong Markov property implies
that, conditional on $(Y_t^{(1)}, Y_t^{(2)})_{0 \leq t \leq \lambda^*}$, all
$\binom{\lfloor Js \rfloor}{m}$ $m$-tuples of individuals with an infinite
line of descent at time $\lambda^*$ are equally likely to form the sample
used to construct $\Upsilon_m$.  With this coupling,
$\Upsilon_m^{(3)} = \Upsilon_m$
if $\lambda^* \leq \lambda_{J_2}$ and all individuals sampled when constructing
$\Upsilon_m^{(3)}$ were born before time $\lambda^*$.

Since $\Upsilon_m^{(2)} =_d \Upsilon_m^{(3)}$, 
Proposition \ref{branchprop} will be proved
if the couplings described in the previous paragraph work
well enough that $P(\Upsilon_m^{(3)} \neq \Upsilon_m)$ and
$P(\Upsilon_m^{(1)} \neq \Upsilon_m^{(2)})$ can both be bounded by
$C/(\log N)^2$.
These bounds follow from Lemma \ref{J1J2}, and Lemma \ref{age} below.

\begin{Lemma}
Let $(\xi_t')_{t=0}^{\infty}$ be a random walk on $\Z$ such that $\xi_0' = 1$
and, for all $k$, $P(\xi_{t+1}' = k+1|\xi'_t = k) = 1/(2-s)$ and
$P(\xi_{t+1}' = k-1|\xi'_t = k) = (1-s)/(2-s)$.
Let $\xi = (\xi_t)_{t=0}^{\infty}$
be the Markov process whose law is the same as the conditional law of
$(\xi_t')_{t=0}^{\infty}$ given $\xi_t' \geq 1$ for all $t$.
Let $\kappa_n = \inf\{t: \xi_t = n\}$.  For all positive integers
$n$, we have $E[\kappa_{n+1} - \kappa_n] \leq (2-s)/s$.
\label{kaplem}
\end{Lemma}

\begin{proof}
Note that $\kappa_1 = 0$ and $\kappa_2 = 1$.  Therefore,
$E[\kappa_2 - \kappa_1] = 1$.  Suppose $E[\kappa_n - \kappa_{n-1}] \leq
(2-s)/s$.  Let $D_n = \#\{t: \kappa_n \leq t < \kappa_{n+1},
\xi_t = n, \mbox{ and }\xi_{t+1} = n-1\}$ be the number of times that
$\xi$ goes from $n$ to $n-1$ before hitting $n+1$.  Since
$l_n = P(\xi_t = n+1|\xi_{t-1} = n) \geq 1/(2-s)$, we have that
$D_n + 1$ follows a geometric distribution with parameter
$l_n \geq 1/(2-s)$.  Therefore, $E[D_n] = (1/l_n) - 1 \le 1-s$.
Note that each time that $\xi$ goes from $n$ to $n-1$, it must
eventually return to $n$, which takes expected time
$E[\kappa_n - \kappa_{n-1}]$.  Thus,
$E[\kappa_{n+1} - \kappa_n] = 1 + E[D_n](1 + E[\kappa_n - \kappa_{n-1}])
\leq 1 + (1-s)[1 + (2-s)/s] = (2-s)/s$.  The lemma now follows by induction.
\end{proof}

\begin{Lemma}
The probability that an individual chosen at random at time
$\lambda_{J_2}$ was born after $\lambda_{J_1}$ is at most $C/(\log N)^2$.
\label{age}
\end{Lemma}

\begin{proof}
Define $({\tilde Y}_t)_{t=0}^{\infty}$ such that if
$0 = \tau_0 < \tau_1 < \dots$ are the jump times of
$(Y_t^{(1)} + Y_t^{(2)})_{t \geq 0}$,
then ${\tilde Y}_t = Y_{\tau_t}^{(1)} + Y_{\tau_t}^{(2)}$.  Let
${\tilde \lambda}_k = \inf\{t: {\tilde Y}_t = k\}$.  The number of births
between $\lambda_{J_1}$ and $\lambda_{J_2}$ is at most
${\tilde \lambda}_{J_2} - {\tilde \lambda}_{J_1}$.  We have
$E[{\tilde \lambda}_{J_2} - {\tilde \lambda}_{J_1}] \leq
[(2-s)/s](J_2 - J_1)$ by Lemma \ref{kaplem}.  Note that
$$
\frac{J_2 - J_1}{J_2} \leq \frac{J(1 - s^{-1}\sqrt{(\log J)/J})^{-1} -
J(1 + s^{-1}\sqrt{(\log J)/J})^{-1} + 2}{J (1 + s^{-1}\sqrt{(\log J)/J})^{-1}}
\leq C \sqrt{ \frac{\log J}{J}},
$$
so the probability that a randomly-chosen individual at
time $\lambda_{J_2}$ was born after $\lambda_{J_1}$ is at most
$$\bigg( \frac{2-s}{s} \bigg) \bigg( \frac{J_2 - J_1}{J_2} \bigg) 
\leq C \sqrt{\frac{\log J}{J}} \leq \frac{C}{(\log N)^2},$$
where the last inequality holds because $J = \lfloor (\log N)^a \rfloor$
for some $a > 4$.
\end{proof}

\section{Approximating the distribution of $\Theta$}

In this section, we complete the proof of Theorem \ref{log2thm} by
proving Propositions \ref{stage2prop1}, \ref{stage2prop2}, and
\ref{stage1prop}.  We will use the notation $W_k$, $\zeta_k$,
$Y_k$, and $Z_i$ introduced before the statement of Theorem 1.2
in the introduction.  Recall also that $L = \lfloor 2Ns \rfloor$.

In subsection 7.1, we prove Propositions \ref{stage2prop1}
and \ref{stage2prop2}, which pertain to the random variables
$Z_i$ introduced in the paintbox construction given in the
introduction.  The rest of the section is devoted to the
proof of Proposition \ref{stage1prop}.  
In subsection 7.2, we introduce random variables $Z_i'$ 
using the branching process.  In subsection 7.3, we state some
lemmas comparing the $Z_i$ and $Z_i'$, and explain how these
lemmas imply Proposition \ref{stage1prop}.  In subsection 7.4,
we present some results related to Polya urns that are needed
to prove these lemmas, and finally the lemmas are proved in
subsection 7.5.

\subsection{Proofs of Propositions \ref{stage2prop1} and \ref{stage2prop2}}

\vspace{-.1in}
\mn{\it Proof of Proposition \ref{stage2prop1}.}
Since $P(Z_1 = Z_2 = k|V_k) \leq V_k^2$, we have
$P(Z_1 = Z_2 = k) \leq E[V_k^2] = E[\zeta_k^2 W_k^2] =
E[\zeta_k^2] E[W_k^2]$.  Since $E[\zeta_k^2] = E[\zeta_k] = r/s$
and $E[W_k^2] = 2/k(k+1)$, it follows that
$$
\hphantom{xxxxxxxxxxx}
P(Z_1 = Z_2 > \lfloor Js \rfloor) \leq
\sum_{k = \lfloor Js \rfloor + 1}^{L} \frac{2r}{sk(k+1)} \leq
\frac{2r}{s \lfloor Js \rfloor} \leq \frac{C}{(\log N)^{1 + a}}.
\hphantom{xxxxxxxxxxx}\Box
$$

We next prove Proposition \ref{stage2prop2}, which says that the
distribution of the number of $i$ such that $Z_i > \lfloor Js \rfloor$
is approximately binomial.  We begin with a lemma which gives an
approximation to $P(Z_i > \lfloor Js \rfloor)$.

\begin{Lemma}
$P(Z_i > \lfloor Js \rfloor) = q_J + O \left( 1/(\log N)^5 \right).$
\label{ziqj}
\end{Lemma}

\mn{\it Proof.}
By the construction in the introduction,
$P(Z_i = k|Z_i \leq k) = E[V_k] = E[\zeta_k]E[W_k] = r/sk.$
Therefore, 
$P(Z_i \leq \lfloor Js \rfloor) = \prod_{k=\lfloor Js \rfloor + 1}^L
( 1 -  r/sk).$
This is the same as the probability that none of the events
$A_{\lfloor Js \rfloor + 1}, \dots, A_L$ occurs if the events are
independent and $P(A_k) = r/sk$.  Since 
$$
\sum_{k = \lfloor Js \rfloor + 1}^L \bigg( \frac{r}{sk} \bigg)^2 \leq
\frac{r^2}{s^2  \lfloor Js \rfloor} \leq \frac{C}{(\log N)^6},
$$ 
it follows from the Poisson approximation
result on p. 140 of Durrett (1996) that
$$
P(Z_i > \lfloor Js \rfloor) =
1 - \exp \bigg(- \sum_{k = \lfloor Js \rfloor + 1}^L
\frac{r}{sk} \bigg) + O \bigg( \frac{1}{(\log N)^6} \bigg).
$$
If $1 \leq y_1 < y_2$, then
$0 \leq \sum_{k = \lfloor y_1 \rfloor}^{\lfloor y_2 \rfloor}
\frac{1}{k} - \log \big( \frac{y_2}{y_1} \big) \leq
2/\lfloor y_1 \rfloor.$  Therefore,
\begin{align}
\bigg| \sum_{k=J+1}^{2N} \frac{1}{k} - \sum_{k = \lfloor Js \rfloor + 1}^
{\lfloor 2Ns \rfloor} \frac{1}{k} \bigg| &\leq \frac{1}{J} +
\bigg| \sum_{k=J}^{2N} \frac{1}{k} - \log \bigg( \frac{2N}{J} \bigg) \bigg|
+ \bigg| \log \bigg( \frac{2Ns}{Js} \bigg) - 
\sum_{k = \lfloor Js \rfloor}^{\lfloor 2Ns \rfloor} \frac{1}{k} \bigg| 
\nonumber \\
&\leq \frac{3}{J} + \frac{2}{\lfloor Js \rfloor} \leq \frac{C}{(\log N)^a}.
\nonumber
\end{align}
It follows that 
$$
\hphantom{xxxxxxx}
P(Z_i > \lfloor Js \rfloor) =
1 - \exp \bigg(- \sum_{k = J+1}^{2N}
\frac{r}{sk} \bigg) + O \bigg( \frac{1}{(\log N)^5} \bigg)
= q_J + O \bigg( \frac{1}{(\log N)^5} \bigg).
\hphantom{xxxxxxx}\Box
$$

\begin{proof}[Proof of Proposition \ref{stage2prop2}]
Let $\eta_k = \#\{i: Z_i = k\}$.  Then $D = \eta_{\lfloor Js \rfloor + 1} +
\dots + \eta_L$.
Define the sequence $({\tilde \eta}_k)_{k=\lfloor Js \rfloor + 1}^L$
such that ${\tilde \eta}_L$ has a Binomial$(n, r/sL)$
distribution and, conditional on ${\tilde \eta}_{k+1}, \dots, 
{\tilde \eta}_L$, the distribution of ${\tilde \eta}_k$ is Binomial with
parameters $n - {\tilde \eta}_{k+1} - \dots - {\tilde \eta}_L$ and $r/sk$.
Thinking of flipping $n$ coins and continuing to flip those that don't show tails,
it is easy to see that $\tilde D = {\tilde \eta}_{\lfloor Js \rfloor + 1}
+ \dots + {\tilde \eta}_L$ has a binomial distribution with parameters $n$
and $\gamma$, where $\gamma = P(Z_i > \lfloor Js \rfloor)$.
To compare $D$ and $\tilde D$ we note that 
$$
P(\eta_k \geq 2|\eta_{k+1}, \dots, \eta_L) \leq
\binom{n}{2} E[V_k^2] = \binom{n}{2} E[\zeta_k] E[W_k^2] =
\binom{n}{2} \frac{2r}{sk(k+1)}
$$ 
and $P({\tilde \eta}_k \geq 2|{\tilde \eta}_{k+1} \dots, {\tilde \eta}_L)
\leq \binom{n}{2} (r/sk)^2.$
By Lemma \ref{coup}, we can couple the $\eta_k$ and ${\tilde \eta}_k$ such that
$P(\eta_k \neq {\tilde \eta}_k|\eta_l = {\tilde \eta}_l
\mbox{ for }l = k+1, \dots, L) \leq Cr/k^2$ for all $k$.
Therefore, $$P(\eta_k \neq {\tilde \eta}_k \mbox{ for some }k > 
\lfloor Js \rfloor)
\leq \sum_{k = \lfloor Js \rfloor + 1}^L \frac{Cr}{k^2} \leq
\frac{Cr}{\lfloor Js \rfloor} \leq \frac{C}{(\log N)^5}.$$
This result, combined with Lemma \ref{ziqj}, gives the proposition.
\end{proof}

\subsection{Random variables $Z_i'$ from the branching process}

It remains only to prove Proposition \ref{stage1prop}, which requires
considerably more work.  For convenience, let $H = \lfloor Js \rfloor$.
>From this point forward, $Z_1, \dots, Z_n$ will be random variables
defined as in the introduction but with $L = H$,
so that the associated marked partition $\Pi$ has the distribution 
$Q_{r,s, H}$.  Our goal is to describe the
distribution of the marked partition $\Upsilon_n$
from Propositions \ref{branchprop} and \ref{stage1prop} using random variables
$Z_1', \dots, Z_n'$, where $Z_i'$ will be the number of
individuals with an infinite line of descent at the time when the type
of the $i$th individual first appeared.
We will then prove Proposition \ref{stage1prop}
by comparing the distribution of $(Z_1', \dots, Z_n')$ to the distribution
of $(Z_1, \dots, Z_n)$.

Define times $0 = \gamma_1 < \gamma_2 < \dots < \gamma_H$ such that
$\gamma_j = \inf\{t: Y_t^{(1)} = j\}$ is the
first time that the branching process has $j$
individuals with an infinite line of descent.  Note that
$(\gamma_{j+1} - \gamma_j)_{i=1}^{H-1}$ is a sequence of independent random
variables, and the distribution of $\gamma_{j+1} - \gamma_j$ is exponential
with rate $js$.  Whenever a new individual with an infinite line of
descent is born, it has a new type with probability $r$.  Also,
each individual with an infinite line of descent is giving birth to
a new individual with a finite line of descent at rate $2(1-s)$.
Since a new individual has a new type with probability $r$,
between times $\gamma_j$ and $\gamma_{j+1}$, births of individuals with
new types occur at rate $2jr(1-s)$.  Whenever such a birth occurs,
the type of the individual with an infinite line of descent changes
with probability $1/2$.  Thus, between times $\gamma_j$ and
$\gamma_{j+1}$, we can view the branching process as consisting of
$j$ lineages with infinite lines of descent, and their types are
changing at rate $r(1-s)$.  It follows that if,
for some $j \geq 1$, we choose
at random one of the $j$ individuals at time $\gamma_{j+1}-$ with an
infinite line of descent, the probability that its ancestor at time
$\gamma_j$ is not of the same type is
\begin{equation}
\frac{r(1-s)}{r(1-s) + js}.
\label{bptype}
\end{equation}
Furthermore, for $j \geq 2$, the probability that its ancestor at
time $\gamma_j$ is not of the same type as its ancestor at time
$\gamma_j-$ is $r/j$ because, with probability $r$, exactly one of
the individuals at time $\gamma_j$ is of a type that did not exist
at time $\gamma_j-$.  It follows that for $j \geq 2$, the probability
that the individual sampled at time $\gamma_{j+1}-$ has a different
type from its ancestor at time $\gamma_j-$ is
\begin{equation}
\frac{r(1-s)}{r(1-s) + js} + \frac{js}{r(1-s)+js} \bigg( \frac{r}{j} \bigg)
= \frac{r}{r(1-s) + js} \leq \frac{r}{js}.
\label{rjs}
\end{equation}
Likewise, the probability that at least one of the $j$ individuals
with an infinite line of descent at time $\gamma_{j+1}-$ has a different
ancestor at time $\gamma_j-$ is
$$\frac{r(1-s)}{r(1-s) + s} + \frac{s}{r(1-s) + s} (r) = 
\frac{r}{r(1-s) + s}.$$

Let $\sigma'(1), \dots, \sigma'(n)$ represent $n$ individuals sampled
at random from those with an infinite line of descent at time $\gamma_H$.
Then we can take the partition $\Upsilon_n$ to be defined such that
$i \sim_{\Upsilon_n} j$ if and only if $\sigma'(i)$ and $\sigma'(j)$
have the same type, and the marked block is
$\{i: \sigma'(i)$ has the same type as the individual at time $0\}.$
Now define $Z'_1, \dots, Z'_n$ as follows.  Let $Z'_i = 1$ if the ancestor
at time $0$ of $\sigma'(i)$ has the same type as $\sigma'(i)$.
Otherwise, define $$Z'_i = \max\{k: \sigma'(i) \mbox{ has a different type
from its ancestor at time }\gamma_k-\}.$$  

If $Z'_i \neq Z'_j$, then since each new type is different
from all types previously in the population, $\sigma'(i)$ and
$\sigma'(j)$ have different types.  If $Z'_i = Z'_j$, then
$\sigma'(i)$ and $\sigma'(j)$ have the same type unless
$\sigma'(i)$ and $\sigma'(j)$
have different ancestors at time $\gamma_{Z'_i + 1}-$ because they both have
the same type as their ancestor at time $\gamma_{Z'_i + 1}-$.
We will show in Lemma \ref{Ri1} below that the probability that
$Z'_i = Z'_j$ and $\sigma'(i)$ and $\sigma'(j)$
have different ancestors at time $\gamma_{Z'_i + 1}-$ is
$O((\log N)^{-2})$.  Therefore, the probability
that, for some $i$ and $j$, we have $Z'_i = Z'_j$ but
$\sigma'(i)$ and $\sigma'(j)$ have different types is $O((\log N)^{-2})$.
Furthermore, it follows from (\ref{bptype}) that the individuals
$\{\sigma'(i): Z'_i = 1\}$ have the same type as the individual
at time $0$ with probability $s/(r(1-s) + s)$.  Define the marked partition
$\Upsilon'_n$ of $\{1, \dots, n\}$ such that $i \sim_{\Upsilon'_n} j$
if and only if $Z'_i = Z'_j$, and independently with probability
$s/(r(1-s) + s)$, mark the block $\{i: Z'_i = 1\}$.  The
preceding discussion implies that
\begin{equation}
|P(\Upsilon_n = \pi) - P(\Upsilon'_n = \pi)| \leq \frac{C}{(\log N)^2}
\label{mainupeq}
\end{equation}
for all $\pi \in {\cal P}_n$.  Thus, for proving Proposition \ref{stage1prop},
we may consider $\Upsilon'_n$ instead of $\Upsilon_n$.
This will be convenient because $\Upsilon'_n$ is defined
from $Z'_1, \dots, Z'_n$ in the same way that $\Pi$ is defined
from $Z_1, \dots, Z_n$.  Consequently, once we establish Lemma
\ref{Ri1} below, the remainder of the proof of Proposition
\ref{stage1prop} will just involve comparing the $Z_i$ and $Z'_i$.

\begin{Lemma}
If $i \neq j$ then
\begin{equation}
P(Z'_i = Z'_j \mbox{ and }\sigma'(i) \mbox{ and }\sigma'(j)
\mbox{ have different ancestors at time }\gamma_{Z'_i + 1}-) \leq
\frac{C}{(\log N)^2}.
\end{equation}
\label{Ri1}
\end{Lemma}

\vspace{-.1in}
\begin{proof}
First note that if $Z_i' = Z_j' = k$, then $\sigma'(i)$ and $\sigma'(j)$ have
the same type as their ancestor at time $\gamma_{k+1}-$.  If they have
different ancestors at time $\gamma_{k+1}-$, there must be a
$\gamma \in (\gamma_k, \gamma_{k+1})$ such that either $\sigma'(i)$ or
$\sigma'(j)$ has an ancestor of a different
type at time $\gamma-$ but not at time $\gamma$.  The other of $\sigma'(i)$
and $\sigma'(j)$ must have an ancestor of a different type at time
$\gamma_k-$ than at time $\gamma-$.  Given that $\sigma'(i)$ and
$\sigma'(j)$ have different ancestors at time $\gamma_{k+1}-$, the probability
that both of these things happen if $k \geq 2$ is
$$
\bigg( \frac{2r(1-s)}{2r(1-s) + ks}\bigg)\bigg( \frac{r}{r(1-s) + ks} \bigg) 
\leq \frac{2r^2}{k^2 s^2}.
$$ 
The first factor is the probability that $\sigma'(i)$ or
$\sigma'(j)$ has an ancestor of a different type at some time $\gamma-$,
while the second factor is the probability from (\ref{rjs})
that the other of $\sigma'(i)$
and $\sigma'(j)$ has an ancestor of a different type at time
$\gamma_k-$ than at time $\gamma-$.
If $k = 1$, then this conditional probability becomes
$$\bigg( \frac{2r(1-s)}{2r(1-s) + ks} \bigg)
\bigg( \frac{r(1-s)}{r(1-s) + ks} \bigg) \leq \frac{2r^2}{k^2 s^2}$$
by (\ref{bptype}).  Therefore, if $i \neq j$,
the probability that $Z'_i = Z'_j$ and $\sigma'(i)$ and $\sigma'(j)$
have different ancestors at time $\gamma_{Z'_i + 1}-$ is at most
$\sum_{k=1}^H \frac{2r^2}{k^2s^2} \leq C/(\log N)^2,$
as claimed.
\end{proof}

\subsection{Comparison of the $Z_i$ and $Z_i'$, and proof of
Proposition \ref{stage1prop}}

We first prove two fairly straightforward lemmas, one for the
$Z_i$ and one for the $Z_i'$.  Lemma \ref{Rlem} allows us
to disregard the possibility that the $Z_i'$ may take more
than two distinct values greater than one, as well as the
possibility that there may be two distinct values greater than
one, with multiple occurrences of the higher value.  
Lemma \ref{Ulem2} rules out the same possibilities for the $Z_i$.

\begin{Lemma}
\begin{align}
& P(Z'_1 = j, Z'_2 = k, Z'_3 = l \mbox{ for some }
2 \leq j < k < l) \leq \frac{C(\log (\log N))^3}{(\log N)^3},
\label{re1} \\
& P(Z'_1 = j, Z'_2 = Z'_3 = k \mbox{ for some }2 \leq j < k)
\leq \frac{C}{(\log N)^2}.
\label{re2}
\end{align}
\label{Rlem}
\end{Lemma}

\upat{\it Proof.}
>From (\ref{rjs}), we get $P(Z'_3 = l) \leq r/sl$,
$P(Z'_2 = k|Z'_3 = l) \leq r/sk$, and
$P(Z'_1 = j|Z'_2 = k, Z'_3 = l)\leq r/sj$.  Thus, the probability on the
left-hand side of (\ref{re1}) is at most
$$
\sum_{j=1}^H \sum_{k=j}^H \sum_{l=k}^H \bigg( \frac{r}{ls} \bigg)
\bigg( \frac{r}{ks} \bigg) \bigg( \frac{r}{js} \bigg)
\leq \frac{C(\log (\log N))^3}{(\log N)^3}.
$$

Conditional on the event that $\sigma'(2)$ and
$\sigma'(3)$ have different ancestors at time $\gamma_{m+1}-$, the
probability that they have the same ancestor at time $\gamma_m-$ is
$\binom{m}{2}^{-1} = 2/m(m-1)$.  Therefore, the probability that
$\sigma'(2)$ and $\sigma'(3)$ have the same ancestor at time
$\gamma_{k+1}-$ is at most $\sum_{m=k+1}^H 2/m(m-1) \leq 2/k$.  
The probability that $Z'_2 = Z'_3 = k$ given that
$\sigma'(2)$ and $\sigma'(3)$ have the same ancestor at time $\gamma_{k+1}-$ is
at most $r/ks$. Also, for $j < k$, we have $P(Z'_1 = j|Z'_2 = Z'_3 = k) \leq r/js$.
Combining these results with Lemma \ref{Ri1}, we can bound the probability
on the left-hand side of (\ref{re2}) by
$$
\hphantom{xx}
\frac{C}{(\log N)^2} + \sum_{j=1}^H \sum_{k=j+1}^H \bigg( \frac{r}{js} \bigg)
\bigg( \frac{r}{ks} \bigg) \bigg( \frac{2}{k} \bigg) 
\leq \frac{C}{(\log N)^2} +
\frac{2r^2}{s^2} \sum_{j=1}^H \sum_{k=j+1}^H \frac{1}{jk^2}
\leq \frac{C}{(\log N)^2}.
\hphantom{xxx}\Box
$$

\begin{Lemma}
\begin{align}
& P(Z_1 = j, Z_2 = k, Z_3 = l \mbox{ for some }
2 \leq j < k < l) \leq \frac{C(\log (\log N))^3}{(\log N)^3}, \nonumber \\
& P(Z_1 = j, Z_2 = Z_3 = k \mbox{ for some }2 \leq j < k)
\leq \frac{C}{(\log N)^2}.
\label{ue2}
\end{align}
\label{Ulem2}
\end{Lemma}

\begin{proof}
Fix $j,k,l$ such that $2 \leq j < k < l \leq H$.  We have
$P(Z_3 = l|Z_3 \leq l) = \frac{r}{sl}$,
$P(Z_2 = k|Z_3 = l, Z_2 \leq k) = \frac{r}{sk}$,
$P(Z_1 = j|Z_2 = k, Z_3 = l, Z_1 \leq j) = \frac{r}{sj}$, and hence
$$
P(Z_1 = j, Z_2 = k, Z_3 = l) \leq \bigg( \frac{r}{sj} \bigg) \bigg( \frac{r}{sk} \bigg)
\bigg( \frac{r}{sl} \bigg). 
$$
Summing as in the proof of
Lemma \ref{Rlem} gives the first result. To prove (\ref{ue2}), first note that
$$
P(Z_2 = Z_3 = k) \leq E[V_k^2] = E[\zeta_k^2] E[W_k^2]
= \frac{2r}{sk(k+1)}
$$
and $P(Z_1 = j|Z_2 = Z_3 = k) \leq r/sj$, then compute as in
the proof of Lemma \ref{Rlem}.
\end{proof}

Throughout the rest of this section, we will use the notation 
$$
q_{k,a,n} = \frac{(k-1)a!(n-a+k-2)!}{(n+k-1)!}.
$$
We now state four more lemmas related to the $Z_i$ and $Z'_i$.  Their
proofs will be given after we explain how they imply Proposition
\ref{stage1prop}.

\begin{Lemma}
Suppose $1 \leq a \leq n-1$.  Then
\begin{align}
P(Z'_1 = l, Z'_2 = \dots = Z'_{a+1} = k, &Z'_{a+2} = \dots = Z'_n = 1
\mbox{ for some } 2 \leq k < l) \nonumber \\
&= \frac{r^2}{s^2} \sum_{k=2}^H
\sum_{l=k+1}^H \frac{q_{k,a,n-1}}{l} + O \bigg( \frac{1}{(\log N)^2} \bigg).
\nonumber
\end{align} 
\label{3brlem}
\end{Lemma}

\begin{Lemma}
Suppose $1 \leq a \leq n-1$.  Then
\begin{align}
P(Z_1 = l, Z_2 = \dots = Z_{a+1} = k, &Z_{a+2} = \dots = Z_n = 1
\mbox{ for some } 2 \leq k < l) \nonumber \\
&= \frac{r^2}{s^2} \sum_{k=2}^H
\sum_{l=k+1}^H \frac{q_{k,a,n-1}}{l} + O \bigg( \frac{1}{(\log N)^2} \bigg).
\nonumber
\end{align} 
\label{3bulem}
\end{Lemma}

\begin{Lemma}
If $2 \leq a \leq n$, then
\begin{align}
P(&Z'_1 = \dots = Z'_a = k \mbox{ and }Z'_{a+1} = \dots = Z'_n = 1 
\mbox{ for some }k \geq 2) \nonumber \\
&= \frac{r}{s} \sum_{k=2}^H q_{k,a,n} - \frac{nr^2}{s^2}
\sum_{k=2}^H \sum_{l=k+1}^H \frac{q_{k,a,n}}{l} +
O \bigg( \frac{1}{(\log N)^2} \bigg).
\label{ag2eq} \\
P(&Z'_1 = k \mbox{ and }Z'_2 = \dots = Z'_n = 1 
\mbox{ for some }k \geq 2) = \frac{r}{s} \sum_{k=2}^H q_{k,1,n}\nonumber \\
& - \frac{nr^2}{s^2}
\sum_{k=2}^H \sum_{l=k+1}^H \frac{q_{k,1,n}}{l}
- \frac{(n-1)r^2}{s^2} \sum_{k=2}^H \sum_{l=2}^{k-1} \frac{1}{k(n+l-2)}
+ O \bigg( \frac{1}{(\log N)^2} \bigg).
\label{ae1eq}
\end{align}
\label{r2blem}
\end{Lemma}

\begin{Lemma}
If $2 \leq a \leq n$, then
\begin{align}
P(& Z_1 = \dots = Z_a = k \mbox{ and }Z_{a+1} = \dots = Z_n = 1 
\mbox{ for some }k \geq 2) \nonumber \\
&= \frac{r}{s} \sum_{k=2}^H q_{k,a,n} - \frac{nr^2}{s^2}
\sum_{k=2}^H \sum_{l=k+1}^H \frac{q_{k,a,n}}{l} +
O \bigg( \frac{1}{(\log N)^2} \bigg).
\label{uag2eq} \\
P(&Z_1 = k \mbox{ and }Z_2 = \dots = Z_n = 1 
\mbox{ for some }k \geq 2) = \frac{r}{s} \sum_{k=2}^H q_{k,1,n} \nonumber \\
& - \frac{nr^2}{s^2}
\sum_{k=2}^H \sum_{l=k+1}^H \frac{q_{k,1,n}}{l}
- \frac{(n-1)r^2}{s^2} \sum_{k=2}^H \sum_{l=2}^{k-1} \frac{1}{k(n+l-2)}
+ O \bigg( \frac{1}{(\log N)^2} \bigg).
\label{uae1eq}
\end{align}
\label{u2blem}
\end{Lemma}

\begin{proof}[Proof of Proposition \ref{stage1prop}]
Let $\pi \in {\cal P}_n$.  If $\pi$ has four or more blocks, or three
blocks of size at least two, then $P(\Upsilon'_n = \pi) \leq C/(\log N)^2$
by Lemma \ref{Rlem} and $Q_{r,s, H}(\pi) \leq C/(\log N)^2$ by
Lemma \ref{Ulem2}.  If $\pi$ has three blocks, at least one
containing just one integer,
then the fact that $|P(\Upsilon'_n = \pi) - Q_{r,s, H}(\pi)| \leq C/(\log N)^2$
follows from Lemmas \ref{Rlem}, \ref{Ulem2}, \ref{3brlem}, and \ref{3bulem}, as
well as the fact that the probabilities that the blocks
$\{i: Z_i = 1\}$ and $\{i: Z'_i = 1\}$ are marked in the two partitions
are both $s/(r(1-s) + s)$.  If $\pi$ has just two blocks, then 
$|P(\Upsilon'_n = \pi) - Q_{r,s, H}(\pi)| \leq C/(\log N)^2$ follows from
Lemmas \ref{r2blem} and \ref{u2blem}, Lemmas \ref{3brlem} and \ref{3bulem}
with $a = n-1$, and equations (\ref{re2}) and (\ref{ue2}).
Finally, when $\pi$ has just one block,
$|P(\Upsilon'_n = \pi) - Q_{r,s, H}(\pi)| \leq C/(\log N)^2$ follows
from Lemmas \ref{r2blem} and \ref{u2blem} with $a = n$, and the fact that
$P(Z_1 = \dots = Z_n = 1)$ and $P(Z_1' = \dots = Z_n' = 1)$ can be
obtained by subtracting from one the remaining possibilities.
Proposition \ref{stage1prop} now follows from these
results and (\ref{mainupeq}).
\end{proof}

\subsection{Polya urn facts}

It remains to prove Lemmas \ref{3brlem}, \ref{3bulem}, \ref{r2blem},
and \ref{u2blem}.  In this subsection, we establish three lemmas
that are related to Polya urns.  The first two lemmas are standard
and straightforward, and their proofs are omitted.

\begin{Lemma}
Suppose $X$ has a beta distribution with parameters $1$ and $k-1$,
where $k$ is an integer.
Let $U_1, \dots, U_n$ be i.i.d. random variables with a uniform distribution
on $[0,1]$.  Then
$$
P(U_i \leq X \mbox{ for }i = 1, \dots, a \mbox{ and }
U_i > X \mbox{ for }i = a+1, \dots, n) = q_{k,a,n}.
$$
\label{betalem}
\end{Lemma}

\begin{Lemma}
Consider an urn with one red ball and $k-1$ black balls.  Suppose that $n$ new
balls are added to the urn one at a time.  Each new ball is either red or
black, and the probability that a given ball is red is equal to the fraction
of red balls currently in the urn.  Let $S$ be any $a$-element subset of
$\{1, \dots, n\}$.  The probability that the $i$th ball added is red for
$i \in S$ and black for $i \notin S$ is $q_{k,a,n}$. Note that this implies
the sequence of draws is exchangeable.
\label{urnlem}
\end{Lemma}

\begin{Lemma}
In the setting of Lemma \ref{urnlem}, suppose instead
$l-k$ new balls are added to the urn.  Then, suppose we sample $n$ of the
$l$ balls at random.  Let $p_{k,l,a,n}$ be the probability that the
first $a$ balls sampled are red and the next $n-a$ are black.
If $a \geq 1$, then
there exists a constant $C$, which may depend on $a$ and $n$, such
that $|p_{k,l,a,n} - q_{k,a,n}| \leq C/kl$ for all $k$ and $l$.
\label{urnlem2}
\end{Lemma}

\mn{\it Proof.}
It follows from Lemma \ref{urnlem} that, conditional on the event
that none of the original $k$ balls is in the sample of $n$, the probability
that the first $a$ balls sampled are red and the next $n-a$ are black is
exactly $q_{k,a,n}$.  The probability that the sample of $n$ balls 
contains exactly $j$ of the original $k$ balls, an event we call $D_{j,k}$, is
\begin{equation}
\frac{\binom{k}{j} \binom{l-k}{n-j}}{\binom{l}{n}} \le
\bigg( \frac{k^j}{j!} \bigg) \bigg( \frac{(l-k)^{n-j}}{(n-j)!} \bigg)
\bigg( \frac{n!(l-n)!}{l!} \bigg) \leq
\binom{n}{j} \frac{k^j l^{n-j} (l-n)!}{l!} \leq C \bigg( \frac{k}{l} \bigg)^j,
\label{Djk}
\end{equation}
since $n$ is a constant and thus so are $a \leq n$ and $j \leq n$.

Conditional on the event $D_{j,k}$, we can calculate the probability that we sample 
$a$ red balls and $n-a$ black balls.  The probability that the original red ball is 
in the sample is $j/k$.  If it is, then by Lemma \ref{urnlem} the probability that $a-1$ of
the other balls in the sample are red is $\binom{n-j}{a-1} q_{k, a-1, n-j}$.  
Likewise, conditional on the event that the original red ball is not in
the sample, the probability that $a$ of the other balls in the sample
are red is $\binom{n-j}{a} q_{k, a, n-j}$.  Thus, conditional on $D_{j,k}$, 
the probability that we sample $a$ red balls and $n-a$ black balls is
$$
\frac{j}{k} \frac{(n-j)!(k-1)(n-j-a+k-1)!}{(n-j-a+1)!(n-j+k-1)!}
+ \frac{k-j}{k} \frac{(n-j)!(k-1)(n-j-a+k-2)!}{(n-j-a)!(n-j+k-1)!}. \nonumber 
$$
Our next step is to bring $\binom{n}{a} q_{k,a,n}$ out in front.
Using that $(m-j)! = m!/(m)_j$ for integers $1 \leq j \leq m$, we get,
for $k \geq 3$,
\begin{align}
&\binom{n}{a} \cdot a! \cdot \frac{(k-1)(n-a+k-2)!}{(n+k-1)!} \nonumber \\
&\times 
\bigg[  \frac{(n-a)_{j-1}}{(n)_j} \frac{j}{k} \frac{(n+k-1)_j}{(n-a+k-2)_{j-1}}
+  \frac{(n-a)_j}{(n)_j} \frac{k-j}{k}  \frac{(n+k-1)_j}{(n-a+k-2)_j}
\bigg]. 
\label{ag2}
\end{align}
Consider the expression in brackets.  Each term can be written as a ratio
of two polynomials in $k$ of the same degree.  Since $a \leq n$ and
$j \leq n$, if $k \rightarrow \infty$ with $n$ fixed, the expression in
brackets is bounded by a constant.  Now, suppose $a = 1$.
The bracketed expression becomes
\begin{align}
&\frac{j(n+k-1)(n+k-2)}{nk(n+k-j-1)} + \frac{(n-j)(k-j)(n+k-1)(n+k-2)}{nk
(n+k-j-1)(n+k-j-2)} \nonumber \\
&= \frac{j(n+k-1)(n+k-2)(n+k-j-2) + (n-j)(k-j)(n+k-1)(n+k-2)}{nk
(n+k-j-1)(n+k-j-2)}. \nonumber
\end{align}
Both the numerator and denominator of this fraction can be written
as third-degree polynomials in $k$ whose leading term is $nk^3$.
Consequently, this fraction minus $1$ can be written as a second-degree
polynomial in $k$ divided by a third-degree polynomial in $k$, which
can be bounded by $Ck^{-1}$ for some constant $C$.

Note that
\begin{equation}
q_{k,a,n} = \frac{(k-1)a!(n-a+k-2)!}{(n+k-1)!} \leq
\frac{a! (n-a+k-2)!}{(n+k-2)!} = \frac{a!}{(n+k-2)_a} \leq \frac{C}{k^a}.
\label{qkan}
\end{equation}
To compare $p_{k,l,a,n}$ and $q_{k,a,n}$ when $a \geq 2$, we will
break up the probability $p_{k,l,a,n}$ by conditioning on the number
of the original $k$ balls that were sampled.  Conditional on sampling
$j \geq 1$ of the original $k$ balls, the probability that the first $a$
balls sampled are red and the next $n-a$ are black is $\binom{n}{a}^{-1}$
times the probability in (\ref{ag2}), which 
can be bounded by $C q_{k,a,n}$.  The probability of sampling 
$j$ of the original $k$ balls is at most $C(k/l)^j$ by (\ref{Djk}), so
$$|p_{k,l,a,n} - q_{k,a,n}| \leq C \sum_{j=1}^n \bigg( \frac{k}{l} \bigg)^j
q_{k,a,n} \leq Ck^{-a} \sum_{j=1}^n \bigg( \frac{k}{l} \bigg)^j
\leq Ck^{-a} n \cdot \frac{k}{l} \leq \frac{C}{kl}.$$
Finally, when $a = 1$, we have
$$
\hphantom{xxxxxxxxxxx}
|p_{k,l,a,n} - q_{k,a,n}| \leq C \sum_{j=1}^n \bigg( \frac{k}{l} \bigg)^j
q_{k,a,n} \frac{C}{k} \leq C \sum_{j=1}^n \bigg( \frac{k}{l} \bigg)^j k^{-2}
\leq \frac{C}{kl}.
\hphantom{xxxxxxxxxxx}\Box
$$

\subsection{Proofs of Lemmas \ref{3brlem}, \ref{3bulem}, \ref{r2blem},
and \ref{u2blem}}

\noindent {\it Proof of Lemma \ref{3brlem}.} For $2 \leq k \leq l$,
let $A_1^{k,l}$ be the event that $\sigma'(1), \dots, \sigma'(n)$ all
have distinct ancestors at time $\gamma_{l+1}-$.  Let $A_2^{k,l}$ be the event
that the ancestor of $\sigma'(1)$ at time $\gamma_{l}-$ has a different type
from the ancestor of $\sigma'(1)$ at time $\gamma_{l+1}-$.
Let $A_3^{k,l}$ be the event
that one of the $k$ individuals at time $\gamma_{k+1}-$ is the ancestor
of $\sigma'(2), \dots, \sigma'(a+1)$ but not $\sigma'(a+2), \dots, \sigma'(n)$,
and let $A_4^{k,l}$ be the event that the ancestor of this individual at time
$\gamma_{k}-$ has a different type.  We claim that
\begin{align}
P(Z'_1 = l, Z'_2 = \dots = Z'_{a+1} = k, &Z'_{a+2} = \dots = Z'_n = 1
\mbox{ for some } 2 \leq k < l) \nonumber \\
&= P \bigg( \bigcup_{2 \leq k < l} A_1^{k,l} \cap A_2^{k,l} \cap A_3^{k,l}
\cap A_4^{k,l} \bigg) + O \bigg( \frac{1}{(\log N)^2} \bigg).
\label{a1a4}
\end{align}

First consider the probability that
$Z'_1 = l, Z'_2 = \dots = Z'_{a+1} = k \mbox{ and } Z'_{a+2} = \dots = Z'_n
= 1 \mbox{ for some } 2 \leq k < l$ but that not all of
$A_1^{k,l}$, $A_2^{k,l}$, $A_3^{k,l}$,
and $A_4^{k,l}$ occur for any $k$ and $l$.
Note that this can only happen in two ways.
One way would be for $A_1^{k,l}$ not to hold, which would mean
$\sigma'(1), \dots, \sigma'(n)$ do not all have distinct ancestors at
time $\gamma_{l+1}-$.  However, it follows from the argument used to prove
(\ref{re2}) that
$P((A_1^{k,l})^c \cap \{Z'_1 = l\} \cap \{Z'_2 = k\} \mbox{ for some }
2 \leq k < l)$ is $O((\log N)^{-2})$.  The second way would be for
$A_1^{k,l}$ to hold but for $\sigma'(2), \dots, \sigma'(a+1)$ not all to have
the same ancestor at time $\gamma_{k+1}-$.  It follows from
Lemma \ref{Ri1}
that this possibility also has probability $O((\log N)^{-2})$.

Next, we consider the probability that $A_1^{k,l}$, $A_2^{k,l}$,
$A_3^{k,l}$, and $A_4^{k,l}$
all hold, but we do not have 
$Z'_1 = l, Z'_2 = \dots = Z'_{a+1} = k$, and $Z'_{a+2} = \dots = Z'_n = 1$.
This is only possible if there is a third time $\gamma$, other than the
times between $\gamma_l$ and $\gamma_{l+1}$ and between
$\gamma_k$ and $\gamma_{k+1}$, such that the type of the
ancestor of one of the individuals
$\sigma'(1), \dots, \sigma'(n)$ at time $\gamma$ is different from the 
type of the
ancestor at time $\gamma-$.  However, it is a consequence of (\ref{re1})
that the probability that this occurs is at most
$O((\log \log N)^3/(\log N)^3)$.  It follows that (\ref{a1a4}) holds.

Recall from the proof of Lemma \ref{Rlem} that if two individuals with
an infinite line of descent are chosen at random at time $\gamma_{k+1}-$,
then the probability that they will have the same ancestor at time
$\gamma_k-$ is $2/k(k-1)$.
Since there are $\binom{n}{2}$ pairs of individuals,
we have $$P(A_1^{k,l}) \geq 1 - \binom{n}{2} \sum_{k=l+1}^H \frac{2}{k(k-1)}
\geq 1 - \binom{n}{2} \frac{2}{l} \geq 1 - \frac{C}{l}.$$
We have $P(A_2^{k,l}|A_1^{k,l}) = r/[r(1-s) + ls]$ by (\ref{rjs}).  Next,
note that if we choose at random one of the $k$ individuals between
times $\gamma_k$ and $\gamma_{k+1}$, then the probability that the individual
born at time $\gamma_{k+1}$ is a descendant of the randomly chosen individual
is $1/k$, and thereafter the probability that each new individual is
a descendant of the randomly chosen individual is the fraction of the
current individuals that are descended from the randomly chosen individual.
This is the same description as the urn problem of Lemma \ref{urnlem2},
so conditional on $A_1^{k,l}$,
the probability that $\sigma'(2), \dots, \sigma'(a+1)$ but not
$\sigma'(a+2), \dots, \sigma'(n)$ are descended from the randomly chosen
individual is $p_{k,l,a,n-1}$.  Therefore,
$P(A_3^{k,l}|A_1^{k,l} \cap A_2^{k,l}) = kp_{k,l,a,n-1}$.  By (\ref{rjs}),
we have $P(A_4^{k,l}|A_1^{k,l} \cap A_2^{k,l} \cap A_3^{k,l}) =
r/[r(1-s) + ks]$.  By the arguments used to prove (\ref{re1}), the probability
that $A_1^{k,l} \cap A_2^{k,l} \cap A_3^{k,l} \cap A_4^{k,l}$ holds for
more than one pair $(k,l)$ is at most $O((\log \log N)^3/(\log N)^3)$.  Thus,
\begin{align}
&P \bigg( \bigcup_{2 \leq k < l} A_1^{k,l} \cap A_2^{k,l} \cap A_3^{k,l}
\cap A_4^{k,l} \bigg) \nonumber \\
&= \sum_{k=2}^H \sum_{l=k+1}^H 
\bigg( \frac{r}{r(1-s) + ls} \bigg) \bigg( \frac{kr}{r(1-s) + ks} \bigg)
\big( p_{k,l,a,n-1} \big) P(A_1^{k,l})
+ O \bigg( \frac{(\log \log N)^3}{(\log N)^3}
\bigg).
\label{rap1}
\end{align}
By Lemma \ref{urnlem2}, we can write $p_{k,l,a,n-1} = q_{k,a,n-1} + \delta$,
where $|\delta| \leq C/kl$.  Also, $P(A_1) = 1 - \eta$, where
$\eta \leq C/l$.  Note that $r/[r(1-s) + ls] \leq r/ls$ and
$kr/[r(1-s) + ks] \leq r/s$.  Recall from (\ref{qkan}) that
$q_{k,a,n} \leq C/k$ for all $a \geq 1$.
To complete the proof, we will need to simplify the four factors inside the sum in
(\ref{rap1}) by obtaining four inequalities.  First, note that
$$
\bigg| \frac{r}{r(1-s) + ls} - \frac{r}{ls} \bigg| =
\frac{r^2(1-s)}{(r(1-s) + ls)(ls)} \leq \frac{r^2}{l^2 s^2}.
$$
Therefore,
\begin{equation}
\sum_{k=2}^H \sum_{l=k+1}^H \bigg| \frac{r}{r(1-s) + ls} - \frac{r}{ls} \bigg|
\bigg( \frac{r}{s} \bigg) \bigg( \frac{C}{k} \bigg) \leq 
\frac{Cr^3}{s^3} \sum_{k=2}^H \sum_{l=k+1}^H \frac{1}{kl^2} \leq
Cr^3 = O \bigg( \frac{1}{(\log N)^3} \bigg). 
\label{rap2}
\end{equation}
Also, 
$$
\bigg| \frac{kr}{r(1-s) + ks} - \frac{r}{s} \bigg| =
\frac{r^2(1-s)}{(r(1-s) + ks)s}\leq \frac{r^2}{ks^2}.
$$
Therefore,
\begin{equation}
\sum_{k=2}^H \sum_{l=k+1}^H \bigg| \frac{kr}{r(1-s) + ks} - \frac{r}{s} \bigg|
\bigg( \frac{r}{ls} \bigg) \bigg( \frac{C}{k} \bigg) \leq
\frac{Cr^3}{s^3} \sum_{k=2}^H \sum_{l=k+1}^H \frac{1}{k^2l} \leq
Cr^3 \log H = O \bigg( \frac{\log \log N}{(\log N)^3} \bigg).
\label{rap3}
\end{equation}
Also,
\begin{equation}
\sum_{k=2}^H \sum_{l=k+1}^H
\bigg( \frac{r}{ls} \bigg) \bigg( \frac{r}{s} \bigg)
\bigg( \frac{C}{kl} \bigg) \leq \frac{Cr^2}{s^2} \sum_{k=2}^H
\sum_{l=k+1}^H \frac{1}{kl^2} \leq Cr^2 = O \bigg( \frac{1}{(\log N)^2} \bigg)
\label{rap4}
\end{equation}
and
\begin{equation}
\sum_{k=2}^H \sum_{l=k+1}^H
\bigg( \frac{r}{ls} \bigg) \bigg( \frac{r}{s} \bigg)
\bigg( \frac{C}{k} \bigg) (1 - P(A_1^{k,l})) \leq 
\frac{Cr^2}{s^2} \sum_{k=2}^H \sum_{l = k+1}^H \frac{1}{kl^2} =
O \bigg( \frac{1}{(\log N)^2} \bigg).
\label{rap5}
\end{equation}
It follows from (\ref{rap1})-(\ref{rap5}) that
$$P \bigg( \bigcup_{2 \leq k < l} A_1^{k,l} \cap A_2^{k,l} \cap A_3^{k,l}
\cap A_4^{k,l} \bigg) = \sum_{k=2}^H \sum_{l=k+1}^H \bigg( \frac{r}{ls} \bigg)
\bigg( \frac{r}{s} \bigg) q_{k,a,n-1} + O \bigg( \frac{1}{(\log N)^2} \bigg),$$
which, combined with (\ref{a1a4}), implies the lemma. \qed

\bigskip
\noindent {\it Proof of Lemma \ref{3bulem}.}
Suppose $2 \leq k < l$.  Let $B_1^{k,l}$ be the event that
$Z_i \leq l$ for $i = 1, \dots, n$.  Let $B_2^{k,l}$
be the event that $Z_1 = l$ and $Z_i \neq l$ for all $2 \leq i \leq n$.
Let $B_3^{k,l}$ be the event that $Z_i \leq k$ for all
$2 \leq i \leq n$.  Let $B_4^{k,l}$ be the event that
$Z_2 = \dots = Z_{a+1} = k$ but $Z_i \neq k$ for $a+2 \leq i \leq n$.
Let $B_5^{k,l}$ be the event that $Z_{a+2} = \dots = Z_n = 1$.
Note that $Z_1 = l, Z_2 = \dots = Z_{a+1} = k$, and $Z_{a+2} = \dots = Z_n = 1
\mbox{ for some } 2 \leq k < l$ if and only if, for some $2 \leq k < l$,
the event $B_1^{k,l} \cap B_2^{k,l} \cap B_3^{k,l} \cap B_4^{k,l} \cap
B_5^{k,l}$ occurs.  Furthermore, the events $B_1^{k,l} \cap \dots \cap
B_5^{k,l}$
are disjoint for different values of $k$ and $l$, so we need to calculate
$\sum_{k=2}^H \sum_{l=k+1}^H P(B_1^{k,l} \cap B_2^{k,l} \cap B_3^{k,l}
\cap B_4^{k,l} \cap B_5^{k,l}).$ We have
\begin{equation}
P(B_1^{k,l}) = \prod_{j=l+1}^H E[(1 - V_j)^n] \geq \prod_{j=l+1}^H
E[1 - nV_j] \geq 1 - \sum_{j=l+1}^H n E[V_j] =
1 - n \sum_{j=l+1}^H \frac{r}{js}.
\label{B1}
\end{equation}
By Lemma \ref{betalem},
$$
P(B_2^{k,l}|B_1^{k,l}) = \frac{r}{s} q_{l,1,n} =
\frac{r}{s} \bigg( \frac{(l-1)(n+l-3)!}{(n+l-1)!} \bigg) =
\frac{r}{sl} \bigg( \frac{l(l-1)}{(n+l-1)(n+l-2)} \bigg) \leq \frac{r}{sl}.
$$
By the same reasoning used to get (\ref{B1}), we have
\begin{equation}
P(B_3^{k,l}|B_1^{k,l} \cap B_2^{k,l}) \geq 1 - (n-1) \sum_{j=k+1}^{l-1}
\frac{r}{js}.
\label{B3}
\end{equation}
By Lemma \ref{betalem},
$$
P(B_4^{k,l}|B_1^{k,l} \cap B_2^{k,l} \cap B_3^{k,l}) =
\frac{r}{s} q_{k,a,n-1}.
$$
Finally, by the argument used to establish (\ref{B1}) and (\ref{B3}),
\begin{equation}
P(B_5^{k,l}|B_1^{k,l} \cap B_2^{k,l} \cap B_3^{k,l} \cap B_4^{k,l}) \geq
1 - (n-a-1) \sum_{j=2}^{k-1} \frac{r}{js}.
\label{B5}
\end{equation}
Note that the product of the probabilities on the right-hand side of
(\ref{B1}), (\ref{B3}), and (\ref{B5}) is at least
$1 - n \sum_{j=1}^H \frac{r}{js} \geq 1 - \frac{C \log H}{\log N}.$
Since $q_{k,a,n-1} \leq C/k$ by (\ref{qkan}), we have
$$
\sum_{k=2}^H \sum_{l=k+1}^H \bigg( \frac{r}{sl} \bigg)
\bigg( \frac{r}{s} \bigg) q_{k,a,n-1} \bigg( \frac{C \log H}{\log N} \bigg)
\leq \frac{C}{(\log N)^3} \sum_{k=2}^H \sum_{l=k+1}^H \frac{(\log H)}{kl}
\leq \frac{C(\log H)^3}{(\log N)^3},
$$
and so
\begin{align}
\sum_{k=2}^H \sum_{l=k+1}^H P (&B_1^{k,l} \cap B_2^{k,l} \cap B_3^{k,l}
\cap B_4^{k,l} \cap B_5^{k,l} ) \nonumber \\
&= \sum_{k=2}^H \sum_{l=k+1}^H \bigg( \frac{r}{sl} \bigg)
\bigg( \frac{r}{s} \bigg) \bigg[ \frac{l(l-1)}{(n+l-1)(n+l-2)} \bigg]
q_{k,a,n-1} + O \bigg( \frac{1}{(\log N)^2} \bigg).
\label{Ulneq}
\end{align}
Finally, note that $\left| 1 - \frac{l(l-1)}{(n+l-1)(n+l-2)} \right|
\leq C/l$ for some constant $C$.  Since $q_{k,a,n-1} \leq C/k$
and $$\sum_{k=2}^H \sum_{l=k+1}^H \bigg( \frac{r}{sl} \bigg)
\bigg( \frac{r}{s} \bigg) \frac{C}{kl} \leq \frac{C}{(\log N)^2},$$
equation (\ref{Ulneq}) remains true if the term in brackets is replaced
by $1$.  The lemma follows. \qed

\bigskip
\noindent {\it Proof of Lemma \ref{r2blem}.}  Let $A_1^k$ be the event
that one of the $k$ individuals at time $\gamma_{k+1}-$ is the ancestor
of $\sigma'(1), \dots, \sigma'(a)$ but not $\sigma'(a+1), \dots, \sigma'(n)$,
and let $A_2^k$ be the event that the ancestor of this individual at time
$\gamma_{k}-$ has a different type.  It follows from
Lemma \ref{Ri1}
that the probability that, for some $k \geq 2$, we have
$Z'_1 = \dots = Z'_a = k \mbox{ and }Z'_{a+1} = \dots = Z'_n = 1$ but
the event $A_1^k \cap A_2^k$ does not occur is at most
$O((\log N)^{-2})$.  We will therefore calculate the probability that the event
$A_1^k \cap A_2^k \cap \{Z'_1 = \dots = Z'_a = k\} \cap
\{Z'_{a+1} = \dots = Z'_n = 1\}$ occurs for some $k \geq 2$.
Note that this occurs for at most one value of $k$, so we may sum the
probabilities over $k = 2, \dots, H$.

Note that $P(A_1^k) = kp_{k,H,a,n}$ and
$P(A_2^k|A_1^k) = r/(r(1-s) + ks)$ by (\ref{rjs}).  It follows that 
$P(A_1^k \cap A_2^k) = [kr/(r(1-s) + ks)] p_{k,H,a,n}$.  Note
that $kr/(r(1-s) + ks) \leq r/s$, and recall that
$|p_{k,H,a,n} - q_{k,a,n}| \leq C/kH$ by Lemma \ref{urnlem2}.  Therefore,
$$
\sum_{k=2}^H \bigg( \frac{kr}{r(1-s) + ks} \bigg) |p_{k,H,a,n} -
q_{k,a,n}| \leq \frac{Cr}{s} \sum_{k=2}^H \frac{1}{kH} \leq
\frac{Cr \log H}{H} \leq \frac{C}{(\log N)^{5}}.
$$  
It follows that
$\sum_{k=2}^H P (A_1^k \cap A_2^k) =
\sum_{k=2}^H \big( \frac{kr}{r(1-s) + ks} \big) q_{k,a,n} +
O \big( 1/(\log N)^5 \big).$
Also, $q_{k,a,n} \leq C/k$, so
$$\sum_{k=2}^H \bigg( \frac{kr}{r(1-s) + ks} - \frac{r}{s} \bigg)
q_{k,a,n} \leq \sum_{k=2}^H \bigg( \frac{r^2}{ks^2} \bigg) \frac{C}{k}
= O \bigg( \frac{1}{(\log N)^2} \bigg).$$  Thus,
\begin{equation}
\sum_{k=2}^H P (A_1^k \cap A_2^k) =
\frac{r}{s} \sum_{k=2}^H q_{k,a,n} + O \bigg( \frac{1}{(\log N)^2} \bigg).
\label{a1ka2k}
\end{equation}

If $A_1^k$ and $A_2^k$ both occur, then we will have
$Z'_1 = \dots = Z'_a = k$ and $Z'_{a+1} = \dots = Z'_n = 1$ unless
either $Z'_i = l$ for some $i = 1, \dots, n$ and $l \notin \{1, k\}$, or
$Z_i' = k$ for some $i \geq a + 1$.  By Lemma \ref{Ri1}, we have
$P(A_1^k \cap A_2^k \cap \{Z_i = k\} \mbox{ for some }k \geq 2 
\mbox{ and }i \geq a+1) \leq
C/(\log N)^2$.  Therefore, we only need to consider the possibility that
$Z'_i = l$ for some $i = 1, \dots, n$ and $l \notin \{1, k\}$.
We will treat separately the cases $l < k$ and $l > k$.  Note that
by (\ref{re1}), the probability that $A_1^k$ and $A_2^k$ both occur,
$Z'_i = l_1$, and $Z'_j = l_2$, where $l_1$ and $l_2$
are distinct integers not in $\{1, k\}$ is at most
$O((\log \log N)^3/(\log N)^3)$.

We first consider $l > k$.  By (\ref{re2}) the probability that
$A_1^k$ and $A_2^k$ both occur and $Z'_i = Z'_j = l$ for some $i \neq j$
is $O((\log N)^{-2})$.
By the same argument used to prove Lemma \ref{3brlem}, the probability
that $A_1^k \cap A_2^k$ for some $k$ but $Z'_i = l$ for some $l > k$ is
\begin{equation}
\frac{nr^2}{s^2} \sum_{k=2}^H \sum_{l=k+1}^H \frac{q_{k,a,n}}{l} +
O \bigg( \frac{1}{(\log N)^2} \bigg).
\label{lgk}
\end{equation}
There are two differences
between this formula and the result of Lemma \ref{3brlem}, which can
be explained as follows.  First, in place of the event $A_2^{k,l}$, 
we need the event that, for some $i = 1, \dots, n$, the ancestor of
$\sigma'(i)$ at time $\gamma_l-$ has a different type from the ancestor of
$\sigma'(i)$ at time $\gamma_{l+1}-$.  This is why the double summation
is multiplied by $n$.  Second, instead of $A_3^{k,l}$, we need one
of the individuals at time $\gamma_{k+1}-$ to be the ancestor of
$\sigma'(1), \dots, \sigma'(a)$ but not $\sigma'(a+1), \dots, \sigma'(n)$,
rather than $\sigma'(2), \dots, \sigma'(a+1)$ but not
$\sigma'(a+2), \dots, \sigma'(n)$.  This is why we have $q_{k,a,n}$
in the formula rather than $q_{k,a,n-1}$.  Otherwise, 
the calculation proceeds as before.

If $a \geq 2$, a consequence of (\ref{re2}) is that the probability
that $A_1^k \cap A_2^k$ for some $k$ but $Z'_i = l$ for some $l < k$
is $O((\log N)^{-2})$.  Thus, (\ref{ag2eq}) follows by subtracting
(\ref{lgk}) from (\ref{a1ka2k}).  Now, consider the case $a = 1$.
Let $S$ be an $d$-element subset of $\{2, \dots, n\}$.  By the argument
used to prove Lemma \ref{3brlem}, the probability that, for some
$2 \leq l < k$, the events $A_{1,k}$ and $A_{2,k}$ occur but
$Z'_i = l$ for $i \in S$ and $Z'_i = 1$ for $i \in \{2, \dots, n\}
\setminus S$ is
$$\frac{r^2}{s^2} \sum_{k=2}^H \sum_{l=2}^{k-1} \frac{q_{l,d,n-1}}{k}
+ O \bigg( \frac{1}{(\log N)^2} \bigg).$$
Summing this over $d = 1, \dots, n-1$ and all subsets $S$ of size $d$,
we get that the probability that $A_{1,k}$ and $A_{2,k}$ occur but
$Z'_i = l$ for $i \in S$ and $Z'_i = 1$ for $i \in \{2, \dots, n\}
\setminus S$ for some nonempty $S \subset \{2, \dots, n\}$ is
\begin{equation}
\frac{r^2}{s^2} \sum_{k=2}^H \sum_{l=2}^{k-1} \frac{1}{k}
\bigg( \sum_{d=1}^{n-1} \binom{n-1}{d} q_{l,d,n-1} \bigg)
+ O \bigg( \frac{1}{(\log N)^2} \bigg).
\label{7285}
\end{equation}
Using the probabilistic interpretation of the $q_{l, d, n-1}$ as in Lemma
\ref{urnlem}, we have
$$\sum_{d=1}^{n-1} \binom{n-1}{d} q_{l,d,n-1} =
1 - q_{l,0,n-1} = 1 - \frac{(l-1)((n-1)+l-2)!}{((n-1)+l-1)!} =
1 - \frac{l-1}{n+l-2} = \frac{n-1}{n+l-2}.$$  Thus, (\ref{7285}) becomes
\begin{equation}
\frac{(n-1) r^2}{s^2} \sum_{k=2}^H \sum_{l=2}^{k-1} \frac{1}{k(n+l-2)}
+ O \bigg( \frac{1}{(\log N)^2} \bigg).
\label{llk}
\end{equation}
We get (\ref{ae1eq}) by subtracting (\ref{llk}) and (\ref{lgk})
from (\ref{a1ka2k}). \qed

\begin{Lemma}
Let $\delta_1, \dots, \delta_N \in (0,1)$.  Assume that
$\delta = \delta_1 + \dots + \delta_n \in (0,1)$.  Then
$$\delta(1 - \delta) \leq 1 - \prod_{n=1}^N (1 - \delta_n) \leq \delta.$$
\label{dlem}
\end{Lemma}

\proof The second inequality follows from $|\prod_{n=1}^N 1 - \prod_{n=1}^N (1 - \delta_n)| \le
\sum_{n=1}^N \delta_n$. To prove the first inequality using the second, note that
\begin{align}
1 - \prod_{n=1}^N (1 - \delta_n) &= \sum_{m=1}^N \bigg(
\prod_{n=1}^{m-1} (1 - \delta_n) - \prod_{n=1}^m (1 - \delta_n) \bigg)
\nonumber \\
&= \sum_{m=1}^N \bigg( \prod_{n=1}^{m-1} (1 - \delta_n) \bigg) \delta_m
\geq \sum_{m=1}^N (1 - \delta) \delta_m = \delta(1 - \delta). \qed \nonumber
\end{align}

\noindent {\it Proof of Lemma \ref{u2blem}.}
Let $B_1^k = \{Z_i \leq k$ for $i = 1, \dots, n \}$.
Let $B_2^k = \{Z_i = k$ for $1\le i \le a$ and
$Z_j < k$ for $ a+1 \le j \le n\}$.  Let $B_3^k = \{
Z_i = 1$ for $a+1 \leq i \leq n\}$.  We have
\begin{align}
P(B_1^k \cap B_2^k \cap B_3^k &\mbox{ for some }k \geq 2) 
= \sum_{k=2}^H P(B_1^k) P(B_2^k|B_1^k) P(B_3^k|B_1^k \cap B_2^k) \nonumber \\
&= \sum_{k=2}^H \bigg( \prod_{l=k+1}^H E[(1 - V_l)^n] \bigg)
\bigg( \frac{r}{s} q_{k,a,n} \bigg) \bigg( \prod_{l=2}^{k-1}
E[(1 - V_l)^{n-a}] \bigg).
\label{Binteq}
\end{align}
Using Lemma \ref{betalem}, $$E[(1 - V_l)^m] = \bigg(1 - \frac{r}{s} \bigg) +
\frac{r}{s} q_{l,0,m} = \bigg(1 - \frac{r}{s} \bigg) +
\frac{r}{s} \bigg( \frac{(l-1)(m+l-2)!}{(m+l-1)!} \bigg) =
1 - \frac{rm}{s(m+l-1)}.$$
Therefore, the expression on the right-hand side of (\ref{Binteq}) is
$$\frac{r}{s} \sum_{k=2}^H \bigg[ \prod_{l=k+1}^H
\bigg( 1 - \frac{nr}{s(n+l-1)} \bigg) \bigg] \bigg[ \prod_{l=2}^{k-1}
\bigg( 1 - \frac{(n-a)r}{s(n-a+l-1)} \bigg) \bigg] q_{k,a,n}.$$
Let $\delta = \frac{r}{s} \sum_{l=k+1}^H \frac{n}{n+l-1} +
\frac{r}{s} \sum_{l=2}^{k-1} \frac{n-a}{n-a+l-1}.$
Then,
$$
\delta^2 = \frac{r^2}{s^2} \bigg( \sum_{l=k+1}^H
\frac{n}{n+l-1} + \sum_{l=2}^{k-1} \frac{n-a}{n-a+l-1} \bigg)^2
\leq \frac{r^2}{s^2} \bigg(n \sum_{l=1}^H \frac{1}{l} \bigg)^2
\leq Cr^2 (\log H)^2.
$$
Since $q_{k,a,n} \leq C/k$ by (\ref{qkan}), we have
$\frac{r}{s} \sum_{k=2}^H \delta^2 q_{k,a,n} \leq
Cr^3 (\log H)^2 \sum_{k=2}^H \frac{1}{k} \leq Cr^3 (\log H)^3.$
Using Lemma \ref{dlem}, the right-hand side of (\ref{Binteq}) can be written as
\begin{equation}
\frac{r}{s} \sum_{k=2}^H \bigg( 1 - \sum_{l=k+1}^H \frac{nr}{s(n+l-1)}
- \sum_{l=2}^{k-1} \frac{(n-a)r}{s(n-a+l-1)} \bigg) q_{k,a,n} +
O \bigg( \frac{(\log \log N)^3}{(\log N)^3} \bigg).
\label{approx1}
\end{equation}
We have $\frac{1}{l} - \frac{1}{n+l-1} = \frac{n-1}{l(n+l-1)} \leq
\frac{n}{l^2}.$  Since $q_{k,a,n} \leq C/k$, it follows that
\begin{equation}
\frac{nr^2}{s^2} \sum_{k=2}^H \sum_{l=k+1}^H \bigg| \frac{1}{(n+l-1)} -
\frac{1}{l} \bigg| q_{k,a,n} \leq Cr^2 \sum_{k=2}^H \sum_{l=k+1}^H
\frac{1}{k l^2} \leq \frac{C}{(\log N)^2}.
\label{approx2}
\end{equation}
Since $q_{k,a,n} \leq C/k^a$ by (\ref{qkan}), when $a \geq 2$ we have
\begin{equation}
\frac{r^2}{s^2} \sum_{k=2}^H \sum_{l=2}^{k-1} \frac{n-a}{n-a+l-1} q_{k,a,n}
\leq Cr^2 \sum_{l=2}^H \sum_{k=l+1}^H \frac{1}{lk^2} =
O \bigg( \frac{1}{(\log N)^2} \bigg).
\label{approx3}
\end{equation}
By combining (\ref{approx1}), (\ref{approx2}), and (\ref{approx3}),
we get (\ref{uag2eq}) when $a \geq 2$.  When $a = 1$, note that
$$\frac{n-a}{n-a+l-1} q_{k,a,n} = \frac{(n-1)(k-1)}{(n+l-2)(n+k-1)(n+k-2)}.$$
Also, note that
$\big| \frac{k-1}{(n+k-1)(n+k-2)} - \frac{1}{k} \big| \leq \frac{C}{k^2}.$
It follows that, when $a = 1$, we have
\begin{equation}
\frac{r^2}{s^2} \sum_{k=2}^H \sum_{l=2}^{k-1} \frac{n-a}{n-a+l-1} q_{k,a,n}
= \frac{(n-1)r^2}{s^2} \sum_{k=2}^H \sum_{l=2}^{k-1} \frac{1}{k(n+l-2)}
+ O \bigg( \frac{1}{(\log N)^2} \bigg).
\label{approx4}
\end{equation}
Equations (\ref{approx1}), (\ref{approx2}), and (\ref{approx4})
establish (\ref{uae1eq}) when $a = 1$. \qed

\clearpage

\bigskip
\begin{center}
{\bf {\Large References}}
\end{center}

\mn K. B. Athreya and P. E. Ney (1972). {\it Branching Processes}.
Springer-Verlag, New York.

\mn N. H. Barton (1998). The effect of hitch-hiking on neutral genealogies.
{\it Genet. Res., Camb.} {\bf 72}, 123-133.

\mn N. H. Barton, A. M. Etheridge, and A. K. Sturm (2004).  Coalescence
in a random background.  {\it Ann. Appl. Probab.} {\bf 14}, 754-785.

\mn P. Donnelly and T. G. Kurtz (1999).  Genealogical processes for
Fleming-Viot models with selection and recombination.
{\it Ann. Appl. Probab.} {\bf 9}, 1091-1148.

\mn R. Durrett (1996).  {\it Probability: Theory and Examples}.
2nd ed.  Duxbury, Belmont, CA.

\mn R. Durrett (2002).  {\it Probability Models for DNA Sequence
Evolution}.  Springer-Verlag, New York.

\mn R. Durrett and J. Schweinsberg (2004a).  Approximating selective sweeps.
{\it Theor. Pop. Biol.} {\bf 66}, 129-138.


\mn R. Durrett and J. Schweinsberg (2004b).  A coalescent model for the effect of advantageous mutations on the genealogy of a population.  Preprint, available at http://front.math.ucdavis.edu/ math.PR/0411071.

\mn A. M. Etheridge, P. Pfaffelhuber, and A. Wakolbinger (2005).
An approximate sampling formula under genetic hitchhiking.
Preprint, available at http://front.math.ucdavis.edu/math.PR/ 0503485.

\mn V. G. Gadag and M. B. Rajarshi (1987).  On multitype processes
based on progeny length particles of a supercritical Galton-Watson process.
{\it J. Appl. Probab} {\bf 24}, 14-24.

\mn V. G. Gadag and M. B. Rajarshi (1992).  On processes associated
with a super-critical Markov branching process.  {\it Serdica.} {\bf 18},
173-178.

\mn J. H. Gillespie (2000).  Genetic drift in an infinite population:
the pseudohitchhiking model.  {\it Genetics}, {\bf 155}, 909-919.

\mn N. L. Johnson, S. Kotz, and A. W. Kemp (1992).  {\it Univariate discrete distributions}. 2nd. ed.  Wiley, New York.

\mn P. Joyce and S. Tavar\'e (1987).  Cycles, permutations and the
structure of the Yule process with immigration.  {\it Stoch. Proc.
Appl.} {\bf 25}, 309-314.

\mn N. L. Kaplan, R. R. Hudson, and C. H. Langley (1989).
The ``hitchhiking effect'' revisited.  {\it Genetics.} {\bf 123}, 887-899.

\mn J. F. C. Kingman (1978).  The representation of partition
structures.  {\it J. London Math. Soc.} {\bf 18}, 374-380.

\mn J. F. C. Kingman (1982).  The coalescent.
{\it Stochastic Process. Appl.} {\bf 13}, 235-248.

\mn J. Maynard Smith and J. Haigh (1974).  The hitch-hiking
effect of a favorable gene.  {\it Genet. Res.}  {\bf 23}, 23-35.

\mn P. A. P. Moran (1958).  Random processes in genetics.
{\it Proc. Cambridge Philos. Soc.} {\bf 54}, 60-71.

\mn M. M{\"o}hle and S. Sagitov (2001).
A classification of coalescent processes for haploid exchangeable
population models.  {\it Ann. Probab.} {\bf 29}, 1547-1562.

\mn N. O'Connell (1993).  Yule process approximation for the
skeleton of a branching process.  {\it J. Appl. Probab.} {\bf 30}, 725-729.

\mn J. Pitman (1999).  Coalescents with multiple collisions.
{\it Ann. Probab.} {\bf 27}, 1870-1902.  

\mn M. Przeworski (2002).  The signature of positive selection at
randomly chosen loci.  {\it Genetics.} {\bf 160}, 1179-1189.

\mn S. Sagitov (1999).  The general coalescent with asynchronous mergers of ancestral lines.  {\it J. Appl. Probab.} {\bf 36}, 1116-1125.

\mn J. Schweinsberg (2000).  Coalescents with simultaneous multiple
collisions.  {\it Electron. J. Probab.} {\bf 5}, 1-50.

\mn K. L. Simonsen, G. A. Churchill, and C. F. Aquadro (1995).
Properties of statistical tests of neutrality for DNA polymorphism data.
{\it Genetics.} {\bf 141}, 413-429.

\mn W. Stephan, T. Wiehe, and M. W. Lenz (1992).  The effect of
strongly selected substitutions on neutral polymorphism: Analytical
results based on diffusion theory.  {\it Theor. Pop. Biol.} {\bf 41}, 237-254.

\bigskip
\bigskip
\begin{tabular}{p{3in}l}
Department of Mathematics, 0112 & Department of Mathematics \\
University of California at San Diego & Malott Hall \\
9500 Gilman Drive & Cornell University \\
La Jolla, CA 92093-0112 & Ithaca, NY 14853-4201 \\
E-mail: jschwein@math.ucsd.edu & E-mail: rtd1@cornell.edu
\end{tabular}

\end{document}